%% file: stochsys24.tex
\documentclass[11pt]{article}

\usepackage[mathscr]{euscript}

\usepackage{multirow}

\usepackage[utf8]{inputenc}

\usepackage{pgfplots}
\usepackage{tikz}
 \usetikzlibrary{arrows,decorations.pathmorphing,through,positioning,fit, patterns}
 \usetikzlibrary{shapes}
 \usetikzlibrary{shapes.multipart}
 \usetikzlibrary{calc}
 \usetikzlibrary{decorations.pathreplacing}
 \usetikzlibrary{plotmarks}
 \usetikzlibrary{arrows.meta}
\definecolor{azulcito}{RGB}{0,80,150}
\definecolor{verdecito}{RGB}{80,150,0}
\definecolor{rojito}{RGB}{150,0,80}
\pgfplotsset{width=\graphwidth, height=\graphheight, grid style = {dashed, gray!50!white}, every tick label/.append style = {font=\footnotesize}, compat=1.18, legend style={font=\footnotesize}, every axis/.append style={label style={font=\footnotesize}}}

\tikzset{>=stealth}

\newcommand{\graphheight}{7cm}
\newcommand{\graphwidth}{11cm}
\newcommand{\given}{\;\middle|\;}

\usepackage[small, margin=1cm]{caption}

\usepackage{appendix}
\usepackage{color}
\definecolor{strcolor}{rgb}{0.6, 0.2, 0.6}
\definecolor{commentcolor}{rgb}{0.3125, 0.5, 0.3125}
\definecolor{keycol}{rgb}{0, 0, 1}

\usepackage{amssymb}
\usepackage{amsmath}
\usepackage{bbm}
\usepackage{amsthm}
\newtheorem{assumption}{Assumption}
\newtheorem{example}{Example}
\newtheorem{corollary}{Corollary}
\newtheorem{proposition}{Proposition}
\newtheorem{lemma}{Lemma}
\newtheorem{theorem}{Theorem}
\newtheorem{definition}{Definition}
\newtheorem*{unnumberedtheorem}{Theorem}
\newcommand{\Halmos}{}
\newcommand{\halmos}{}
\usepackage[top=2.5cm,bottom=2.5cm,left=2.3cm,right=2.3cm]{geometry}

\newcommand {\E}[1]{\mathrm{E}\left[ #1 \right]}

\newcommand {\p}{{\rm P}}
\renewcommand{\leq}{\leqslant}  
\renewcommand{\geq}{\geqslant}  

\def\blot{\quad \mbox{$\vcenter{ \vbox{ \hrule height.4pt
				\hbox{\vrule width.4pt height.9ex \kern.9ex \vrule width.4pt}
				\hrule height.4pt}}$}}

\usepackage{natbib}
\bibpunct[, ]{(}{)}{,}{a}{}{,}%
\def\BIBand{and}%

\gdef\AQ#1{}
\gdef\CQ#1{}

\begin{document}

\title{Optimal local storage policy based on stochastic intensities and its large scale behavior}

\author{Matias Carrasco, Andres Ferragut, Fernando Paganini\thanks{Research partially supported by AFOSR-US under grant \#FA9550-23-1-0350.}
}

\date{Facultad de Ingenieria, Universidad ORT Uruguay\\ \bigskip
\emph{Preprint submmitted to Stochastic Systems} \\ \bigskip
\today}

	




	
	%
	
\maketitle

\begin{abstract}
  In this paper, we analyze the optimal management of local memory systems, using the tools of stationary point processes. We provide a rigorous setting of the problem, building upon recent work, and characterize the optimal causal policy that maximizes the hit probability. We specialize the result for the case of renewal request processes and derive a suitable large scale limit as the catalog size $N\to \infty$, when a fixed fraction $c$ of items can be stored. We prove that in the limiting regime, the optimal policy amounts to comparing the stochastic intensity (observed hazard rate) of the process with a fixed \emph{threshold}, defined by a quantile of an appropriate limit distribution, and derive asymptotic performance metrics, as well as sharp estimates for the pre-limit case. Moreover, we establish a connection with optimal timer based policies for the case of monotonic hazard rates. We also present detailed validation examples of our results, including some close form expressions for the miss probability that are compared to simulations. We also use these examples to exhibit the significant superiority of the optimal policy for the case of regular traffic patterns.
\end{abstract}

\section{Introduction}\label{sec:intro}

In modern computing systems, a crucial component for improving performance at several layers is the use of a \emph{local memory}, typically referred to as a \emph{cache}. The rationale is that storing frequently used items at a readily available location can improve latency, reducing the cost of retrieval from a more remote location. This strategy is pervasive in computing systems, from local caching of instructions at the processor level, texture caching in graphics processing units, disk caching for quickly retrieving data from hard disks, content caching in web applications, geographically distributed caching in content delivery networks and cloud storage gateways keeping readily available items stored in the cloud data centers.

A \emph{local memory} consists of a certain amount of memory space that may store a subset of items locally and temporarily; the goal is to appropriately choose from a large catalog, the subset of items more likely to be requested next. All the aforementioned applications can be subsumed into this basic structure. For simplicity, we will assume that items are of homogeneous size, and thus the space constraint is reflected in the \emph{number} of items $C$ which may be locally stored, from the catalog of size $N \gg C$.  

Classical analyses of local memory systems have focused on modeling the (discrete time) sequence of item requests. Based on this empirical sequence, the system must determine which items are stored, and which must be \emph{evicted} from memory to make room for others. In a stationary regime, one natural strategy would be to store at all times the $C$ items with the highest \emph{popularity}, measured by their mean intensity of requests; this simple \emph{static policy} requires, however, popularities to be known. A practical approximation is the \emph{Least-Frequently-Used (LFU)} eviction policy: here, item request intensities are dynamically estimated from the arrival sequence, and the least popular are evicted. Under mild assumptions on stationary arrivals, LFU will eventually converge to the static policy.

A popular alternative is the \emph{Least-Recently-Used (LRU)} policy, which keeps in memory the $C$ most recently requested items. Upon receiving a request, it will store the item if not already present, and in that case evict the oldest item in the sequence of requests.  This method is better suited to handle highly \emph{correlated} demands, where a requested item is more likely to be fetched again, i.e. bursty request patterns. Smoother variants have also been analyzed, and also combinations of both approaches such as LFRU (cf. \cite{bilal2017lfru}) have been proposed. We review the relevant literature below.

A main drawback of the classical analysis is that relevant (continuous) \emph{time} information is neglected: by focusing only on the sequence of requests, the models ignore inter-request times that are important characteristics of the request processes. An alternative approach that gives time the center stage, are \emph{timer-based (TTL)} caching policies; here, items are stored upon request, and evicted only after a certain amount of time has elapsed since their last appearance in the request stream. This is a common approach in Internet-based systems such as Domain Name System queries and Web applications.

A crucial step towards incorporating time information is the seminal paper \cite{fofack2014performance}, where the incoming request stream is assumed to be a stationary point process on the real line. The mean intensities of the underlying request processes for each item capture their relative popularity, while by modeling the inter-request times we can express different types of behavior: for instance, heavy-tailed inter-request times are well-adapted to bursty arrival patterns. This approach has led to new insights in the analysis of both replacement and timer-based caching policies, which we also review below.

Within this modeling framework, a natural question to ask is: what is the \emph{optimal} memory management policy? By this we mean the one that maximizes the \emph{hit rate}, the frequency of successful retrievals from local memory. In \cite{ferragut2016caches_sigmetrics, ferragut2018optimal}, the optimality question for TTL caching policies was investigated; in particular, it was shown that if requests form a renewal process, optimality may be characterized by the {\em hazard rate} function of the inter-arrival distribution. More recently   \cite{panigrahy2022upper} address optimality in replacement policies, characterizing it for general point processes in terms of the \emph{stochastic intensity}, a notion applicable beyond the renewal case. 

In this paper we develop the latter connection more extensively, bringing to bear the machinery of stationary point processes as in \cite{bremaud2020point} to formally characterize causal memory management policies and the condition for optimality.  Subsequently, our aim is to understand the large-scale behavior of the optimal replacement policy for large $N$, under some natural assumptions: namely, that inter-arrival times for different items are drawn from a \emph{scale family}, parametrized by a distribution of intensities which has an asymptotic limit. We characterize the optimal memory management as a \emph{threshold} policy on observed stochastic intensities, and prove the threshold has a deterministic limit. We also obtain a closed expression for the optimal \emph{performance} in the large-scale limit, as a function of the fundamental parameters of the scale family and the popularity distribution. We further investigate the properties of threshold policies before the asymptotic limit, showing that under monotonicity of the process hazard rate function, they become equivalent to timer-based policies, This also yields as a result the convergence of the two types of policies (replacement or TTL) to a common optimum in the asymptotic limit. Our results are illustrated by a series of parametric examples, and their properties are further exhibited through stochastic simulations.

\subsection{Related work}

Being essential to modern computing systems, the performance analysis of caching and local memory management has a long history. The seminal work on replacement based algorithms with a fixed memory size started in \cite{king1971paging,gelenbe1973unified}. Despite its apparent simplicity, the analysis of replacement policies is complex even for a single cache: in \cite{king1971paging}, the author gives an explicit expression of the hit probability of the LRU policy, with exponential complexity. In \cite{gelenbe1973unified}, it is shown that under the so-called \emph{independent reference model (IRM)} assumption, where the requests are independent and identically distributed, similar replacement policies such as first-in-first-out achieve the same hit probability. Exact computation is however intractable, \cite{dan1990approximate} provides an approximate computational procedure. Their method has been further extended in \cite{rosensweig2010approximate} to the case of networks of cache systems. Also in \cite{gast2015list}, another approach is given for list-based replacement policies.

A related line of work that includes replacement policies are analyses based on the ``Move-To-Front'' rule, which is equivalent to LRU. A first step in this direction is \cite{fill1996limits}, which addresses the limit cost of the Move-to-Front rule. By exploiting the connection between MTF and LRU, \cite{jelenkovic1999asymptotic,jelenkovic2008persistent} provide asymptotic expressions for the hit probability under the IRM assumption and Zipf popularities with parameter $\beta$. With a fundamentally different approach, based on Laplace transforms, the same limit result is obtained in \cite{barrera2010limiting}. Our asymptotic limit assumptions are related to this latter approach. Few results on replacement policies move beyond the IRM assumptions, such as \cite{jelenkovic2003asymptotic,jelenkovic2004least,jelenkovic2006critical} where the authors show some insensitivity properties of LRU in a large scale regime and dependent requests. However, their approach is only valid for light-tailed popularities, and is related to our result below on optimal performance for the non-uniformly integrable popularity limit. Extensions to a network of caches working cooperatively to optimize performance are proposed in \cite{borst2010distributed,ioannidis2010distributed,ioannidis2016adaptive}.

Since exact analysis of LRU is difficult, the most popular technique for approximating its performance  is the so-called ``Che approximation'', introduced in \cite{che2002hierarchical} for Web caches. The crucial point is to define a \emph{characteristic time} representing the average permanence time \emph{common} to all files. This assumption is valid in a large scale regime, and in \cite{fricker2012versatile}, the authors perform second order analysis to explain why this is a good approximation. 

On the other hand, the analysis of policies based on timers is more recent, associated with the growth of Internet caches. A first contribution in this line is \cite{jung2003modeling}, with expressions for the steady state hit probabilities for web pages, later extended in \cite{bahat2005measuring} to include update delay. However, a crucial contribution is the introduction of point process theory to capture request correlations, and in particular move beyond the \emph{sequence} of requests to a continuous time model with general arrival processes. The foundations for this approach were laid down in \cite{fofack2012analysis,fofack2014performance}. In \cite{berger2014exact} the analysis is extended to a family of TTL policies with the focus of approximating LRU performance in linear and tree networks, and a different policy is proposed in \cite{berger2015maximizing} with the aim of maximizing hit ratios by variance reduction. In \cite{bianchi2013check} the Che approximation is analyzed in a more general setting and in \cite{martina2014unified}, it is extended, using TTL cache tools, to renewal arrivals, showing good accuracy in the case of small caches and Zipf popularities, making the connection between replacement and timer based policies. The fact that timer based policies decouple the analysis over independent request streams has led to more amenable generalizations to the case of cache networks \cite{panigrahy2017hit,dehghan2019utility,panigrahy2020network}.

The search for optimal policies has received more recent attention. In \cite{ferragut2016caches_sigmetrics,ferragut2018optimal}, the authors formulate the problem of characterizing the optimal timer based policy, and introduced the connection with the \emph{hazard rate function} of the inter-request times. Using tools from convex optimization, the authors show that \emph{decreasing} hazard rates lead to a non-trivial policy that achieves optimality, and characterize its large scale behavior. However, more regular traffic with \emph{increasing} hazard rates do not benefit from caching at all. Later, in \cite{panigrahy2022upper}, the authors consider fixed memory system (i.e. replacement based policies) and incorporate the notion of stochastic intensity of point processes. This notion is a generalization of the hazard rate, and allows to characterize the optimal \emph{causal} policy, i.e. one that does not make use of future arrival information. In \cite{ferragut2024prefetch}, a generalization of timer-based policies to consider \emph{pre-fetching} was introduced, leading to a non-trivial optimal policy for increasing hazard rates. In this paper, we build upon these works to establish a rigorous connection between the optimal causal policy, and its asymptotic limit, with the timer policies, and in particular we show that they become equivalent in a large scale limit.

\subsection{Main contributions}

We now summarize the main contributions of the paper. Our first result is to establish a rigorous characterization of \emph{causal} policies applied to the local memory or caching problem, in the language of point processes in the real line. Under this characterization, we identify the optimal causal policy, which will depend on the stochastic intensities of the underlying request processes.

The main result of the paper is to prove that, under some suitable additional assumptions, the optimal policy converges to a \emph{threshold} policy, under a fluid scaling appropriate for large systems. This fixed threshold value is used to decide whether an item should be stored. Armed with this result, we characterize the limit of the optimal \emph{miss rate} in large scale systems, through a tractable formula. Therefore, we provide a \emph{universal asymptotic bound} on performance for any practical policy.

A third contribution is based on further studying the behavior of threshold policies: under suitable monotonicity assumptions of the \emph{hazard rate function} of the inter-request times, it can be shown that previously analyzed timer based policies are in fact threshold policies, and in particular the \emph{optimal timer policy} for decreasing hazard rates identified in \cite{ferragut2018optimal} achieves the universal bound.

As a final contribution, we propose a dual policy of the timer-based caching policy dubbed \emph{timer-based prefetching}, which is also asymptotically optimal in the case of increasing hazard rates. Through analysis and simulations we show that indeed this policy can greatly outperform classical caching policies in this context.

\subsection{Organization of the paper}

The paper is organized as follows: in Section \ref{sec:preliminaries} we lay out the main tools of point processes and stochastic intensities required to analyze our system. We then define our local memory model, establish the framework for causal policies and characterize the optimal policy in Section \ref{sec:optimal}. Our main Theorem describing the large scale limit is presented in Section \ref{sec:asymptotic}, and the ensuing universal performance bound is presented in Section \ref{sec:universalbound}. In Section \ref{sec:threshold} we describe the connection between the optimal policy and timer-based policies. Simulations and examples are presented in Section \ref{sec:simulations} and conclusions are given in Section \ref{sec:conclusions}.

\section{Preliminaries and notation}\label{sec:preliminaries}

Throughout this paper, we will consider \emph{stationary point processes} defined on a common probability space $(\Omega,\mathcal{F},P)$. We recall now some basic concepts that will be useful in the following and introduce our notation; we refer the reader to \cite{bremaud2020point} for a thorough treatment. 

A simple and locally finite point process $\Phi$ on the real line is a random and strictly increasing sequence of points $\Phi=\{\tau_k\}_{k\in\mathbb{Z}}$ satisfying $\lim_{k\to\pm\infty}\tau_k=\pm\infty$. More formally, $\Phi$ can be cast as a \emph{random counting measure}, i.e. $\Phi = \sum_k \delta_{\tau_k}$, a measurable map from $(\Omega,\mathcal{F}) \to (M^{\#}(\mathbb{R}),\mathcal{M}^{\#}(\mathbb{R}))$. Here $M^{\#}(\mathbb{R})$  is the space of locally finite measures on $\mathbb{R}$ taking values in $\mathbb{N}\cup \{\infty\}$, and $\mathcal{M}^{\#}(\mathbb{R})$ is the smallest $\sigma$-algebra such that, for all Borel sets $B\in \mathcal{B}(\mathbb{R})$, $\Phi(B) = \sum_k \mathbf{1}_{\{\tau_k \in B\}}$ is measurable. $\Phi(B)$ is assumed to be finite for bounded $B$, and is thus a non-negative integer valued random variable.  By definition, all points $\tau_k$ are different and thus $\Phi(\{x\}) \leq 1$ $P$-a.s. for all $x\in\mathbb{R}$. In order to label the points, we follow the usual convention \cite{bremaud2020point} where $\tau_0(\Phi) \leqslant 0$ and $\tau_1(\Phi)>0$. With this convention $\tau_0=\tau_0(\Phi)$ represents the first point before the time origin of the process $\Phi$.\footnote{Note that $\tau_0$ is a proper random variable since $\{\tau_0\leqslant t\} = \{\Phi((t,0])=0\}$ and thus measurable for all $t\leq  0$. Similarly, any point $\tau_k$ is a random variable.}

Let $S_t(\Phi)$ denote the shift operator for measures in $\mathbb{R}$, i.e. $S_t(\Phi)(B) = \Phi(B+t)$. The point process $\Phi$ is stationary if $\displaystyle{S_t(\Phi)}$ has the same distribution as $\Phi$ for all $t\in\mathbb{R}$. The mean measure of the point process is $\lambda(B):=\E{\Phi(B)}$. If the process is stationary, then this measure is translation-invariant, and thus a multiple of the Lebesgue measure on $\mathbb{R}$, i.e. $\lambda(B)=\lambda m(B)$. The constant $\lambda$ is called the (average) \emph{intensity} of the stationary point process. In what follows we assume $\lambda>0$ to avoid the trivial case where the process has no points.

For a simple stationary point process $\Phi$, an important measure is the \emph{Palm probability} $P_\Phi^0$. This is a probability measure defined in $(\Omega,\mathcal{F})$ that  captures the stochastic behavior of the point process when the observer is \emph{synchronized} with it. In particular $P_\Phi^0(\tau_0=0)=1$, i.e. there is $P_\Phi^0$ a.s. a point at the origin. We refer the reader to \cite{bremaud2020point} for a formal definition.

The key relationship between the Palm probability and the stationary probability is the following \emph{inversion formula} valid for any non-negative real valued measurable function $f:M^\#(\mathbb{R})\to\mathbb{R}_{+}$:
\begin{equation}\label{eq:inversion_formula}
  \E{f(\Phi)} = \lambda{\rm E}^0_\Phi\left[\int_0^{\tau_1} f(S_t(\Phi)) dt\right]
\end{equation}
That is, in order to know the average value of a property in the stationary measure, we can integrate over one cycle of the process using the Palm measure and scale by $\lambda$.

The \emph{inter-arrival distribution} of the point process $\Phi$ is defined as $F_0(t):=P^0_\Phi(\tau_1-\tau_0\leqslant t) = P^0_\Phi(\tau_1\leqslant t)$. Since the process is simple, $F_0$ has support in $\mathbb{R}_+$. Moreover, ${\rm E}^0_\Phi[\tau_1] = 1/\lambda$, which follows from taking $f\equiv 1$ in \eqref{eq:inversion_formula}.

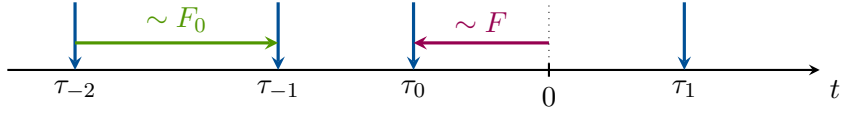
\begin{figure}
  \begin{center}
  \begin{tikzpicture}[scale=0.9]
       
      \draw[->,thick] (0,0) -- (12,0);
      \node[below right] at (12,0) {$t$};
      
      \draw[->,very thick, azulcito] (1,1) -- (1,0);
      \draw[->,very thick, azulcito] (4,1) -- (4,0);
      \draw[->,very thick, azulcito] (6,1) -- (6,0);
      \draw[->,very thick, azulcito] (10,1) -- (10,0);
      \draw[-,thick] (8,-.1) node[below] {$0$} -- (8,0.1);
    
      \node[below] at (1,0) {$\tau_{-2}$};
      \node[below] at (4,0) {$\tau_{-1}$};
      \node[below] at (6,0) {$\tau_{0}$};
      \node[below] at (10,0) {$\tau_{1}$};
       
      \draw[->, very thick, verdecito] (1,.4) -- node[midway, above, anchor=south] {$\sim F_0$} (4,.4);
  
      \draw[->, very thick, rojito] (8,.4) -- node[midway, above, anchor=south] {$\sim F$} (6,.4);

      \draw[dotted] (8,0) -- (8,1);
    \end{tikzpicture}
    \end{center}
  \caption{Inter-arrival and age distribution of a stationary point process.}\label{fig:age_distribution}
\end{figure}

A second important distribution is the \emph{age distribution}, which is the age of the current interval when the process is observed at a point $t$ not synchronized with it. Since the process is stationary, we can take without loss of generality $t=0$ and thus the age distribution is just:
\begin{equation}\label{eq:age_distribution}
  F(t):=P(-\tau_0 \leqslant t) = \lambda \int_0^t (1-F_0(s)) ds,
\end{equation}
a result that also follows  from \eqref{eq:inversion_formula} with $f(\Phi)=\mathbf{1}_{\{\tau_0(\Phi)\geqslant -t\}}$. Note that these distributions are different in general, due to the bias towards larger intervals when sampling in steady state. A depiction of this sampling effect is shown in Figure \ref{fig:age_distribution}.

\subsection{Stochastic intensity}\label{ssec:stochint}
We now introduce the concept of  \emph{stochastic intensity} \cite{bremaud2020point} which is crucial for the analysis in this paper.  First we need the following definition:
\begin{definition}\label{def:history}
    A filtration $\{\mathcal{F}_t\}_{t\in\mathbb{R}}$ in $(\Omega,\mathcal{F})$ (i.e. an increasing family of $\sigma$-algebras contained in $\mathcal{F}$)  is \emph{a history} of the simple locally finite point process $\Phi$ if $\Phi((a,b])$ is $\mathcal{F}_t$-measurable for all $a<b\leqslant t$. The \emph{natural history} of $\Phi$ is $\mathcal{F}_t=\sigma\left(\Phi((a,b]):a<b\leq t\right)$ the smallest filtration satisfying this property.
\end{definition}

We define the stochastic intensity of a process for a given history \cite[Definition 5.1.1]{bremaud2020point}.

\begin{definition}\label{def:stochastic_intensity}
  Let $\Phi$ be a simple locally finite point process and $\mathcal{F}_t$ a history of $\Phi$. If there exists a locally integrable $\mathcal{F}_t$-adapted process $\lambda(t)\geqslant 0$ satisfying:
  \begin{gather} \label{eq:stochint}
    \E{\Phi((s,t])\given \mathcal{F}_s} = \E{\int_s^t \lambda(u)du \given \mathcal{F}_s},
  \end{gather}
then $\lambda(t)$ is called an $\mathcal{F}_t-$\emph{stochastic intensity} of $\Phi$.
\end{definition}

The process $\lambda(t)$ acts as the local likelihood of a point appearing at time $t$ given past information. This notion will play a key role in our particular application. Also, directly from the definition, one can prove that $\E{\lambda(t)} = \E{\lambda(0)} = \lambda$, the average intensity of the process.

As an example, the stationary Poisson process of intensity $\lambda$ satisfies \eqref{eq:stochint} with $\lambda(t)\equiv \lambda$, a deterministic constant, and thus the likelihood of a point appearing in time is independent of the past, reflecting the total randomness property of the Poisson process.

We now highlight some key properties of the stochastic intensity that will be useful for our later analysis. The first one is related to predictability, for this we need the following:
\begin{definition}\label{def:predictable}
  Let $\{\mathcal{F}_t\}_{t\in\mathbb{R}}$ be a filtration on $(\Omega,\mathcal{F})$. The \emph{predictable} $\sigma-$algebra $\mathcal{P}(\mathcal{F}.)$ associated to $\mathcal{F}_t$ is the $\sigma-$algebra on $\mathbb{R}\times \Omega$ generated by the sets of the form:
    \begin{equation*}
      (a,b] \times A, \; a<b, \; A\in\mathcal{F}_a.
    \end{equation*}
    A stochastic process $X(t, \omega)$ taking values on a measurable space $(E,\mathcal{E})$ is $\mathcal{F}_t$-\emph{predictable} if the mapping $(t,\omega) \mapsto X(t,\omega)$ is $\mathcal{P}(\mathcal{F}.)$-measurable.
\end{definition}

In particular, any real-valued and \emph{left-}continuous stochastic process adapted to $\mathcal{F}_t$ is $\mathcal{F}_t$-predictable \cite[Example 5.1.9]{bremaud2020point}. Moreover, the stochastic intensity of a point process $\Phi$, if it exists, can always be chosen to be a.s. predictable, up to a set of Lebesgue measure $0$, and thus is essentially unique \cite[Theorem 5.1.38]{bremaud2020point}.

The following result from \cite{bremaud2020point}\footnote{This is essentially \cite[Theorem 5.11]{bremaud2020point}; in fact, this reference provides a slightly stronger version of the converse, where condition 
\eqref{eq:smoothing_formula} need only be checked for adapted, left continuous $Z(t)$.} provides a \emph{smoothing formula} for predictable processes, based on the stochastic intensity:
\begin{unnumberedtheorem}
  Let $\Phi$ be a simple point process in $\mathbb{R}$ with history $\mathcal{F}_t,$ and $\lambda(t)$ an $\mathcal{F}_t$-adapted and a.s. integrable process. Then: $\lambda(t)$ is an $\mathcal{F}_t$-stochastic intensity of $\Phi$ if and only if 
  \begin{equation}\label{eq:smoothing_formula}
    \E{\int Z(t)\Phi(dt)} = \E{\int Z(t)\lambda(t)dt}
  \end{equation}
 holds for all $\mathcal{F}_t$-predictable processes $Z(t)$.
\end{unnumberedtheorem}

\subsection{Renewal point processes}\label{ssec:renewal}

As an important special case, consider now that $\Phi$ is a \emph{stationary renewal process}, meaning that the inter-request time sequence $\{\tau_{k+1}-\tau_k\}_{k\in\mathbb{Z}}$ are $\p^0_\Phi$ independent and identically distributed random variables, with common distribution $F_0$. We assume that $F_0$ has a density $f_0$ and recall the definition of the \emph{failure rate} or \emph{hazard rate} function associated to $F_0$:
  \begin{equation}\label{eq:def_hazard_rate}
    \eta(t) = \frac{f_0(t)}{1-F_0(t)}.
  \end{equation}

Then if $\mathcal{F}_t$ is the natural history of the process, the stochastic intensity of $\Phi$ is given by \cite[Chapter 7]{daley2003introduction}:
  \begin{equation}\label{eq:stoch_intensity_renewal}
    \lambda(t) = \eta(t-\tau^-(t)),
  \end{equation}
  where $\tau^-(t) = \sup\{\tau_k:\tau_k < t\}$ is the last point before $t$. Note that $\tau^{-}(t)$ is a left continuous process.

  In particular, due to the renewal property, the local intensity of the process depends only on the \emph{age} of the current interval and the hazard rate function of the inter-arrival distribution. 
  
  When $F_0$ is the exponential distribution as in the Poisson process, $\eta(t)\equiv \lambda$ so the stochastic intensity is constant as previously stated.  An interesting parametric example is presented below in \ref{ssec:pareto}, which we will use throughout the paper to illustrate the results.

  For the purposes of our analysis, an important random variable is the stochastic intensity observed when sampling the renewal process at a fixed time, for convenience chosen to be  $t=0$. We call it the \emph{observed hazard rate} (OHR), and denote it by $X:=\lambda(0)=\eta(-\tau^-(0))$. 

To find the distribution of the OHR, there are two underlying probability measures to consider, depending on whether the time under consideration is synchronized with process arrivals. 

If time $t=0$ is \emph{not} synchronized with arrivals, then $\tau^-(0) = \tau_0$ and the distribution is: 
      \begin{equation}\label{def:ohr}
    G(x) := \p(X\leqslant x) = \p(\eta(-\tau_0)\leqslant x) = \p \left(-\tau_0 \in \eta^{-1}{([0,x])}\right) = \int_{\eta^{-1}([0,x])} F(dt),
  \end{equation}
  since $-\tau_0\sim F$, the age distribution. We note that the set $\eta^{-1}([0,x])$ will be an interval if the hazard rates are monotone.

If, instead, we are evaluating the OHR at an arrival time, we must use  the Palm probability $P_\Phi^0$, for which $\tau_0\equiv 0$ a.s., and $\tau^{-}(\tau_0) =  \tau^{-}(0) = \tau_{-1}$. Therefore: $ X =  \eta(-\tau_{-1})$, 
and  its distribution is given by:
  \begin{equation}\label{def:ohrua}
    G_0(x) := \p_\Phi^0 (X\leqslant x) = \p_\Phi^0 (\eta(-\tau_{-1}) \leqslant x) = \int_{\eta^{-1}([0,x])} F_0(dt),
  \end{equation}
  since $-\tau_{-1}\sim F_0$ under the Palm probability. A basic inequality for this distribution follows from the definition of the hazard rate:
\begin{align}\label{eq:lipschitz}
  G_0(x) & = \int_{\eta^{-1}([0,x])} f_0(t) dt =  \int_{\{t: \eta(t) \leqslant x\}} \eta(t)(1-F_0(t)) dt 
  &\leq x \int_{\mathbb{R}} (1-F_0(t)) dt = x E[\tau_1 - \tau_0] = \frac{x}{\lambda}. 
\end{align}

\subsection{Pareto inter-arrival times.} \label{ssec:pareto}

An interesting parametric example of a stationary renewal processes, considered in \cite{ferragut2016caches_sigmetrics} for the same application, is when inter-arrival times follow a heavy-tailed Pareto distribution. In this case, all the above magnitudes have explicit expressions, that we now compute.

\begin{example}[Renewal Pareto process]
  In the above setting, choose
  \begin{equation}\label{eq:def_pareto}
    F_0(t) = 1 - \left(\frac{1}{1+\gamma t}\right)^\alpha, \quad f_0(t) = \frac{\alpha \gamma}{\left(1+\gamma t\right)^{\alpha+1}}\quad (t\geq 0).
  \end{equation}
  Thus, $F_0$ is a Pareto distribution (starting at $0$) with tail parameter $\alpha>1$. The number $\gamma$ acts as a scale parameter and by direct computation, in order to have $\rm{E}^0_\Phi[\tau_1] = 1/\lambda$, it should be chosen such that $\gamma = \frac{\lambda}{\alpha-1}$.
\end{example}

From equations \eqref{eq:age_distribution} and \eqref{eq:def_hazard_rate} we can compute:
\begin{equation}\label{eq:age_pareto}
  F(t) = 1 - \left(\frac{1}{1+\gamma t}\right)^{\alpha-1}, \quad \eta(t) = \frac{\alpha \gamma}{1+\gamma t},\quad (t\geq 0).
\end{equation}

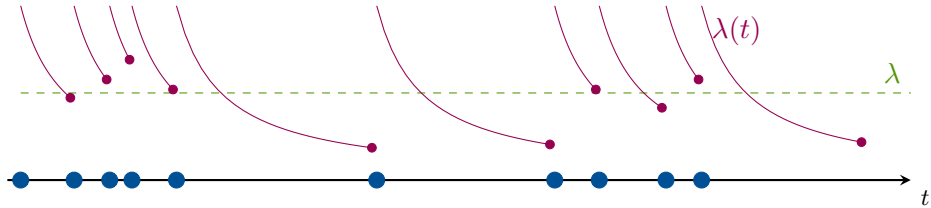
\begin{figure}
    \begin{center}
      \begin{tikzpicture}
        \begin{axis}[xlabel=$t$,
          ymin = -.3,
          ymax = 2,
          xmin = -.3,
          xmax=20,
          xlabel style = {at={(axis cs:20,0)},anchor=north west},
          y axis line style={draw=none},
          x axis line style={->, thick},
          ticks=none,
          axis x line*=middle,
          width=0.8\textwidth,
          height= 0.25\textwidth,
          ]
          \addplot[verdecito,domain=0:20, dashed] {1};
          \addplot[rojito,domain=0:1.2, -Circle] {2/(1+x)};
          \addplot[rojito,domain=1.2:2, -Circle] {2/(1+x-1.2)};
          \addplot[rojito,domain=2:2.5, -Circle] {2/(1+x-2)};
          \addplot[rojito,domain=2.5:3.5, -Circle] {2/(1+x-2.5)};
          \addplot[rojito,domain=3.5:8, -Circle] {2/(1+x-3.5)};
          \addplot[rojito,domain=8:12, -Circle] {2/(1+x-8)};
          \addplot[rojito,domain=12:13, -Circle] {2/(1+x-12)};
          \addplot[rojito,domain=13:14.5, -Circle] {2/(1+x-13)};
          \addplot[rojito,domain=14.5:15.3, -Circle] {2/(1+x-14.5)};
          \addplot[rojito,domain=15.3:19, -Circle] {2/(1+x-15.3)};
          \node[below right, rojito] at (axis cs:15.3,2) {$\lambda(t)$};
          \node[above left, verdecito] at (axis cs:20,1) {$\lambda$};    
          \addplot+[azulcito, mark options={fill=azulcito, scale=1.5}, mark=*,only marks] coordinates {
            (0,0)
            (1.2,0)
             (2,0)
             (2.5,0)
             (3.5,0)
             (8,0)
             (12,0)
             (13,0)
             (14.5,0)
             (15.3,0)
          };
        \end{axis}
      \end{tikzpicture}
      \end{center}
      \caption{Stochastic intensity of a renewal Pareto process, showing decreasing hazard rates.} \label{fig:pareto_example}
  \end{figure}

For this example, note that the hazard rate function is decreasing for any choice of the parameters. Therefore, at an arrival time $\tau_k$, the stochastic intensity \emph{increases} (the hazard rate resets to $\eta(0)$ following \eqref{eq:stoch_intensity_renewal}), and a subsequent arrival becomes more likely. This gives rise to \emph{bursty} traffic as depicted in Figure \ref{fig:pareto_example}.

For this process we can also compute the distributions of the observed hazard rate, at non-synchronized or synchronized times:
\begin{equation}\label{eq:ohr-pareto}
   G(x) = \int_{\eta^{-1}(x)}^\infty F(dt) = \int_{\frac{\alpha}{x}-\frac{1}{\gamma}}^\infty F(dt) = \begin{cases}
     \left(\frac{x}{\alpha\gamma}\right)^{\alpha-1}, &\quad 0\leqslant x \leqslant \alpha \gamma, \\
     1 &\quad x> \alpha\gamma.
   \end{cases}
 \end{equation}

 \begin{equation}\label{eq:ohrua-pareto}
  G_0(x) = \int_{\eta^{-1}(x)}^\infty F_0(dt) = \int_{\frac{\alpha}{x}-\frac{1}{\gamma}}^\infty F_0(dt) = \begin{cases}
    \left(\frac{x}{\alpha\gamma}\right)^{\alpha}, &\quad 0\leqslant x \leqslant \alpha \gamma, \\
    1 &\quad x> \alpha\gamma.
  \end{cases}
\end{equation}

\begin{figure}
  \begin{center}
  \begin{tikzpicture}
    \begin{axis}[xlabel=$t$,
      ymin = -.05,
      ymax = 1.3,
      xmin = -.1,
      xmax=4,
      xlabel style = {at={(axis cs:4,0)},anchor=north west},
      y axis line style={->,thick},
      x axis line style={->, thick},
      xtick=\empty,
      ytick={1},
      axis x line*=middle,
      axis y line*=middle,
      width=0.3\textwidth,
      height= 0.28\textwidth,
      ]
      \addplot[domain=0:4, dashed] {1};
      \addplot[azulcito, thick, domain=0:4] {1-(1/(1+x))^2};
      \node[above left, azulcito] at (axis cs:4,1) {\footnotesize $F_0(t)$};
    \end{axis}
  \end{tikzpicture} \hspace{1cm}
  \begin{tikzpicture}
    \begin{axis}[xlabel=$t$,
      ymin = -.05,
      ymax = 1.3,
      xmin = -.1,
      xmax=4,
      xlabel style = {at={(axis cs:4,0)},anchor=north west},
      y axis line style={->,thick},
      x axis line style={->, thick},
      xtick=\empty,
      ytick={1},
      axis x line*=middle,
      axis y line*=middle,
      width=0.3\textwidth,
      height= 0.28\textwidth,
      ]
      \addplot[domain=0:4, dashed] {1};
      \addplot[azulcito, thick, domain=0:4] {1-(1/(1+x))};
      \node[above left, azulcito] at (axis cs:4,1) {\footnotesize $F(t)$};
    \end{axis}
  \end{tikzpicture} \hspace{1cm}
  \begin{tikzpicture}
    \begin{axis}[xlabel=$t$,
      ymin = -.1,
      ymax = 2.3,
      xmin = -.1,
      xmax=4,
      xlabel style = {at={(axis cs:4,0)},anchor=north west},
      y axis line style={->,thick},
      x axis line style={->, thick},
      xtick=\empty,
      ytick={2},
      axis x line*=middle,
      axis y line*=middle,
      width=0.3\textwidth,
      height= 0.28\textwidth,
      ]
      \addplot[azulcito,thick, domain=0:4] {2/(1+x)};
      \node[above left, azulcito] at (axis cs:4,.5) {\footnotesize $\eta(t)$};
    \end{axis}
  \end{tikzpicture}

  \begin{tikzpicture}
    \begin{axis}[xlabel=$x$,
      ymin = -.05,
      ymax = 1.3,
      xmin = -.1,
      xmax=4,
      xlabel style = {at={(axis cs:4,0)},anchor=north west},
      y axis line style={->,thick},
      x axis line style={->, thick},
      xtick={2},
      ytick = {1},
      axis x line*=middle,
      axis y line*=middle,
      width=0.3\textwidth,
      height= 0.28\textwidth,
      ]
      \addplot[domain=0:4, dashed] {1};
      \addplot[azulcito, thick, domain=0:2] {x/2};
      \addplot[azulcito, thick,domain=2:4] {1};
      \node[above left, azulcito] at (axis cs:4,1) {\footnotesize $G(x)$};
      \draw[dashed] (axis cs:2,0) -- (axis cs:2,1);
  \end{axis}
  \end{tikzpicture}\hspace{1cm}
  \begin{tikzpicture}
    \begin{axis}[xlabel=$x$,
      ymin = -.05,
      ymax = 1.3,
      xmin = -.1,
      xmax=4,
      xlabel style = {at={(axis cs:4,0)},anchor=north west},
      y axis line style={->,thick},
      x axis line style={->, thick},
      xtick={2},
      ytick = {1},
      axis x line*=middle,
      axis y line*=middle,
      width=0.3\textwidth,
      height= 0.28\textwidth,
      ]
      \addplot[domain=0:4, dashed] {1};
      \addplot[azulcito, thick, domain=0:2] {x^2/4};
      \addplot[azulcito, thick,domain=2:4] {1};
      \node[above left, azulcito] at (axis cs:4,1) {\footnotesize $G_0(x)$};
      \draw[dashed] (axis cs:2,0) -- (axis cs:2,1);
  \end{axis}
  \end{tikzpicture}
  \end{center}
  \caption{Distribution shapes for the Pareto renewal example with $\alpha=2, \gamma=1$, and therefore $\lambda=1$.}\label{fig:pareto_example_2}
\end{figure}
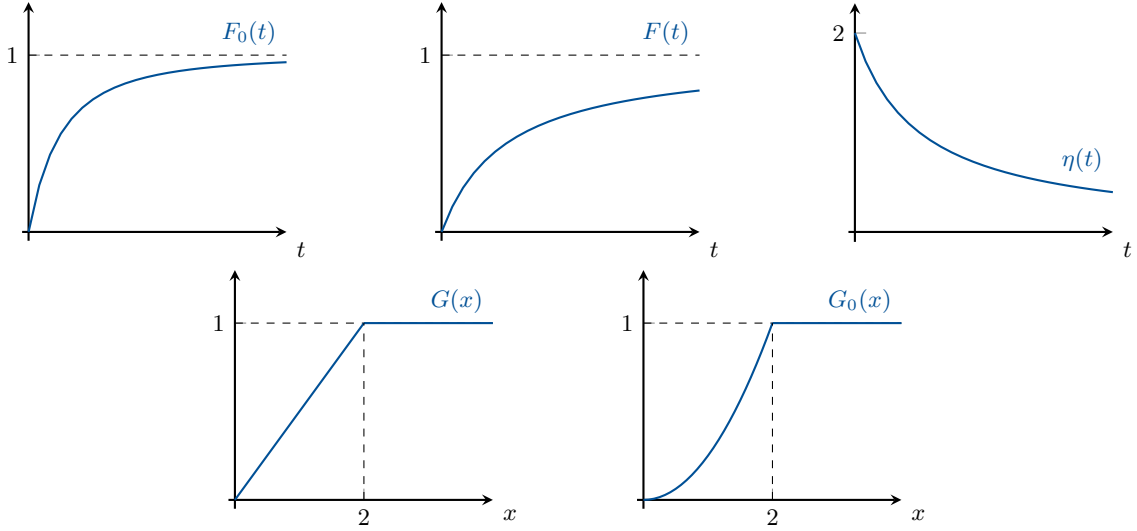

 In Figure \ref{fig:pareto_example_2} we depict the distributions $F_0$, $F$, $G$ and $G_0$, as well as the hazard-rate function $\eta$, for the case $\alpha=2$, $\gamma=1$ (with average intensity $\lambda=1$). In this particular case the non-synchronized OHR distribution is uniform in $[0,2]$.

\section{Local memory systems and optimal storage policy}\label{sec:optimal}

We start by rigorously defining our model for a \emph{local memory system}. Requests from a \emph{catalog} of $N$ (equally sized) items are received. We model item requests by independent stationary stochastic point processes $\Phi_i$, $i=1,\ldots,N$, defined on a common probability space and with finite intensity $\lambda_i>0$; a measure of the popularity of item $i$. By appropriate labelling we can choose $\lambda_1\geqslant \lambda_2 \geqslant \ldots$, i.e. the objects are ordered by decreasing popularities.

The complete request process is thus the superposition $\Phi := \sum_{i=1}^N \Phi_i$, and its intensity is $\lambda^N = \sum_{i=1}^N \lambda_i$. Let $\p^0_\Phi$ denote the Palm probability of the superposition, and $\p^0_{\Phi_i}$ the Palm probability of the $i-$th process. Then, for the superposition of independent processes, we have \cite[Examples 7.2.8 and 7.2.11]{bremaud2020point}:
\begin{subequations}\label{eq:palm_superposition}
\begin{equation}\label{eq:palm_superposition_a}
  \p^0_\Phi(A) = \sum_{i=1}^N \frac{\lambda_i}{\lambda^N}\p^0_{\Phi_i}(A),
\end{equation}
for any $A\in \mathcal{F}$. Moreover, if $\{\Gamma_i\} \subset \mathcal{M}^\#(\mathbb{R})$, then
\begin{equation}\label{eq:palm_superposition_b}
  \p^0_{\Phi_i} (\Phi_1 \in \Gamma_1,\ldots,\Phi_N \in \Gamma_N) = \p_{\Phi_i}^0\left(\Phi_i \in \Gamma_i\right) \prod_{j\neq i} \p\left(\Phi_j \in \Gamma_j\right).
\end{equation} 
\end{subequations}
The interpretation of eqs. \eqref{eq:palm_superposition} is the following: given a point in $\Phi$ occurring at time $0$, $\lambda_i/\lambda^N$ is the probability that this point comes from process $i$. Then, in order to compute Palm probabilities given that the point comes from process $i$, we must use the Palm probability for process $i$ and the stationary probability for any other $j\neq i$, and the processes remain independent.

The local memory is limited, and thus can only keep readily available a subset of size $C<N$ of the items. These items can, upon request, be served from the local memory with lower cost than retrieving them from a central repository. This formulation is quite general, and subsumes the typical notion of \emph{caching}, useful in many applications.\footnote{We specifically avoid using the term \emph{cache} for our system, because in caching normally only the recently requested objects are stored, and our model aims to generalize these policies.} Mathematically, the local memory can keep, at any point in time, a subset:
\begin{equation*}
  \mathcal{C}(t) = \{i_1,\ldots,i_k\} \subset \{1,\ldots,N\},\text{ with } 0\leq k \leq C.
\end{equation*}
We call this process the \emph{storage process} of the system. For a given storage process, we can define the \emph{hit counting process} of the local memory system as:
\begin{equation}\label{eq:hit_process}
  \Psi_\mathcal{C}(B) = \sum_{i=1}^N \int_B \mathbf{1}_{\{i\in \mathcal{C}(t)\}} \Phi_i(dt),
\end{equation}
i.e. the thinned process that counts the requests of objects only if they are stored in memory at request time. The mean intensity $h_\mathcal{C}$ of the process $\Psi_\mathcal{C}$ is called the \emph{hit rate}, and accordingly satisfies $h_\mathcal{C}\leqslant \lambda^N$. 

The natural objective in this setting is to maximize $h_\mathcal{C}$, i.e. the rate at which requests can be served directly from the local memory, by choosing an appropriate \emph{policy} that defines the storage process $\mathcal{C}(t)$. However, this policy should be \emph{causal}, i.e. it cannot make use of future information. Otherwise, the policy that only stores the next arriving item achieves maximum rate using only one unit of memory for any $N$, and the problem is trivial. This is where the predictability notion introduced in Section \ref{ssec:stochint} becomes important.

Let us denote by $\mathcal{F}_t^{(i)}$ the natural history of the $i$-th request process $\Phi_i$, and $\lambda_i(t)$ its stochastic intensity with respect to $\mathcal{F}_t^{(i)}$. Define $\mathcal{F}_t = \sigma\left(\mathcal{F}_t^{(i)};\,i=1,\ldots,N\right)$ the aggregated history. Then, due to independence, we have the following property:
\begin{lemma}\label{lemma:superposition_intensity}
If $\Phi = \sum_i \Phi_i$ is the superposition process, then $\mathcal{F}_t$ is a history of $\Phi$. Furthermore, the total stochastic intensity of $\Phi$ is simply $\lambda(t) = \sum_{i=1}^N \lambda_i(t)$. Moreover, $\lambda_i(t)$ is a stochastic intensity for process $i$ with respect to the enlarged history $\mathcal{F}_t$.
\end{lemma}
\proof{Proof.} 
First, $\Phi((a,b])$ is $\mathcal{F}_t$-measurable for all $a<b\leqslant t$ since the sum function is measurable and $\mathcal{F}_t^{(i)}\subset\mathcal{F}_t$. Second, $\lambda_i(t)$ is the stochastic intensity of $\Phi_i$ with respect to the shared history $\mathcal{F}_t$, since the conditional expectation given $\mathcal{F}^{(i)}_t$ coincides with the conditional expectation given $\mathcal{F}_t$ for any random variable independent of $\mathcal{F}_t^{(j)}$ with $j\neq i$. This fact is in turn a consequence of the independence of the filtrations $\{\mathcal{F}_t^{(i)}\}_{i=1}^N$ and that $\mathcal{F}_t$ is generated by events of the form $\bigcap_{i=1}^N A_i$ with $A_i\in\mathcal{F}_t^{(i)}$, $i=1,\ldots,N$. Finally, that $\lambda(t)$ is the stochastic intensity of $\Phi(t)$ with respect to $\mathcal{F}_t$ follows immediately by linearity of conditional expectation.
\hfill\Halmos
\endproof

The above Lemma guarantees that the stochastic intensity of process $i$ is not altered by superposition with other independent processes; there is no cross-information between the processes, a crucial assumption. As an example of what could happen without independence, consider in particular two processes $\Phi_1$, and $\Phi_2 = S_{t_0}(\Phi_1)$, i.e. a delayed version of the first, with $t_0>0$. Then their stochastic intensities with respect to their own natural histories are just delayed versions of one another: $\lambda_2(t) = \lambda_1(t-t_0)$. However, in the combined history of both processes, the stochastic intensity of process $2$ degenerates into point masses at the $t_0$-shifted points of process $1$.

The storage process defined above is quite general, we now restrict our attention to \emph{causal} policies, which are only the ones that use past information from the process.
\begin{definition}
  Consider a local memory system with independent request processes $\{\Phi_i:i=1\ldots N\}$, aggregated natural history $\mathcal{F}_t$ and capacity $C$. A \emph{causal memory management} policy is a stationary $\mathcal{F}_t$-predictable stochastic process $\mathcal{C}(t)$, with values in the space $\mathscr{P}_C(N)$ of all subsets of $\{1,\ldots, N\}$ of size at most $C$ (equipped with the discrete $\sigma$-algebra).
\end{definition}

We can now compute the stochastic intensity of the hit process for any causal policy:
\begin{lemma}\label{lemma:hit_rate_process_intensity}
Let $\mathcal{C}(t)$ be a causal policy. The stochastic intensity of its associated hit process $\Psi_\mathcal{C}(t)$ defined by \eqref{eq:hit_process} is
\begin{equation}\label{eq:hit_rate_process_intensity}
  \lambda_\mathcal{C}(t) = \sum_{i=1}^N \lambda_i(t)\mathbf{1}_{\{i\in\mathcal{C}(t)\}}.
\end{equation}
\end{lemma}

\proof{Proof.}
We apply the smoothing formula Theorem of Section \ref{ssec:stochint}. For an
$\mathcal{F}_t$-predictable process $Z(t)$ we write:
\[
\begin{aligned}
\E{\int_{\mathbb{R}}Z(t)\Psi_{\mathcal{C}}(dt)} & =\sum_{i=1}^N\E{\int_{\mathbb{R}}Z(t)\mathbf{1}_{\{i\in\mathcal{C}(t)\}}\Phi_i(dt)} =\sum_{i=1}^N\E{\int_{\mathbb{R}}Z(t)\mathbf{1}_{\{i\in\mathcal{C}(t)\}}\lambda_i(t)\,dt}\\
                                                & =\E{\int_{\mathbb{R}}Z(t)\left(\sum_{i=1}^N \mathbf{1}_{\{i\in\mathcal{C}(t)\}}\lambda_i(t)\right)\,dt} =\E{\int_{\mathbb{R}}Z(t)\lambda_{\mathcal{C}}(t)\,dt}.
\end{aligned}
\]
The second equality follows from \eqref{eq:smoothing_formula}, since by Lemma \ref{lemma:superposition_intensity} the stochastic intensity of $\Phi_i$ with respect $\mathcal{F}_t$ is still $\lambda_i(t)$, and the process $Z(t)\mathbf{1}_{\{i\in\mathcal{C}(t)\}}$ is $\mathcal{F}_t$-predictable, because $\mathcal{C}(t)$ is causal. 

Since the identity holds for \emph{any} $\mathcal{F}_t$-predictable process $Z(t)$, the converse implication in the smoothing formula establishes \eqref{eq:hit_rate_process_intensity}. \hfill \Halmos
\endproof

Define now the policy $\mathcal{C}^*(t)$ as follows: at any point in time, rank the items by decreasing stochastic intensity, and store the $C$ highest; ties may be broken arbitrarily.  Note that this policy maximizes the sum of the stochastic intensities of the items in memory:
\begin{equation}\label{eq:optimal_memory}
    \mathcal{C}^{*}(t):=\arg\max_{\mathcal{C}\in \mathscr{P}_C(N)}\sum_{i\in \mathcal{C}}\lambda_i(t)
\end{equation}
This proposal was already presented in \cite{panigrahy2022upper}. The result below provides a detailed proof of optimality within the framework of 
stationary point processes.

\begin{theorem}[Optimal causal memory management policy]\label{thm:optimal}
  For any causal policy $\mathcal{C}(t)$, its stationary hit rate $h_{\mathcal{C}}$ satisfies:
  \begin{equation*}
    h_{\mathcal{C}}= \E{\lambda_{\mathcal{C}}(0)} \leqslant \E{\lambda_{\mathcal{C}^*}(0)} = 
    h_{\mathcal{C}^*}.
  \end{equation*}
  In addition, the hit probability:
  \begin{equation*}
    H_{\mathcal{C}} := \p^0_\Phi\left( \Psi_{\mathcal{C}}(\{0\}) = 1 \right)
  \end{equation*}
  is also maximized by the policy $\mathcal{C}^*$.
\end{theorem}
\proof{Proof.}
Note from eq. \eqref{eq:hit_rate_process_intensity} and the construction of $\mathcal{C}^*$, that for any realization:
\begin{equation*}
  \lambda_{\mathcal{C}}(0) = \sum_{i=1}^N \lambda_i(0)\mathbf{1}_{\{i \in \mathcal{C}\}} \leqslant \max_{\{i_1,\ldots,i_C\}} \sum_{i\in \{i_1.\ldots,i_C\}} \lambda_i(0) = \lambda_{\mathcal{C}^*}(0). 
\end{equation*}
The first inequality follows by taking expectations on both sides with respect to the joint probability measure $\p$ of the arrival processes (non-synchronized with arrivals). 

To derive the second statement, it suffices to show that for any causal policy $h_{\mathcal{C}} = \lambda^N H_{\mathcal{C}}$, i.e. the hit rate is the total rate ``thinned'' by the hit probability. This is a natural property, but non-trivial since the thinning is not independent of the arrival process. 

To proceed, apply \eqref{eq:hit_process} to the set $B=\{0\}$ to conclude that
\begin{align}\label{eq:hit at zero}
\Psi_{\mathcal{C}}(\{0\}) = \#\left\{j: \phi_j(\{0\}) = 1, j \in \mathcal{C}(0)\right\}.
\end{align}
Apply now the Palm probability formula \eqref{eq:palm_superposition} to the event 
$A=\{\omega: \Psi_{\mathcal{C}}(\{0\}) = 1\}$:
  \begin{equation*}
      H_{\mathcal{C}} = \p^0_\Phi\left( \Psi_\mathcal{C}(\{0\}) = 1 \right) = \sum_{i=1}^N \frac{\lambda_i}{\lambda^N}\p^0_{\Phi_i}\left(\Psi_\mathcal{C}(\{0\}) = 1\right).
\end{equation*}
Now note that under $\p^0_{\Phi_i}$, the events $\{\Psi_{\mathcal{C}}(\{0\}) = 1\}$ and $\{i \in \mathcal{C}(0)\}$ coincide almost surely; this follows from \eqref{eq:hit at zero} observing that the set of the right-hand side can (at most) contain the index $j=i$, $\p^0_{\Phi_i}$-almost surely. Therefore, we arrive at 
   \begin{align}\label{eq:hit prob decomp}
      H_{\mathcal{C}} = \sum_{i=1}^N \frac{\lambda_i}{\lambda^N}\p^0_{\Phi_i}\left(i\in \mathcal{C}(0)\right).
\end{align}
Now return to  \eqref{eq:hit_rate_process_intensity} and its expectation, at $t=0$:
\begin{equation*}
  h_{\mathcal{C}} = \E{\sum_{i=1}^N \lambda_i(0)\mathbf{1}_{\{i\in\mathcal{C}(0)\}}} = \sum_{i=1}^N \E{\lambda_i(0)\mathbf{1}_{\{i\in\mathcal{C}(0)\}}}.
\end{equation*}
Since $\mathbf{1}_{\{i\in\mathcal{C}(t)\}}$ is $\mathcal{F}_t-$predictable for our causal policy, it follows from
Papangelou's Theorem \cite[Theorem 7.7.5]{bremaud2020point} that:
\begin{equation*}
  \E{\lambda_i(0)\mathbf{1}_{\{i\in\mathcal{C}(0)\}}} = \lambda_i \mathrm{E}_{\Phi_i}^0\left[\mathbf{1}_{\{i\in\mathcal{C}(0)\}}\right] = \lambda_i \p^0_{\Phi_i}\left(i\in \mathcal{C}(0) \right).
\end{equation*}
Summation over $i$ gives, together with \eqref{eq:hit prob decomp}, $ h_{\mathcal{C}} = \lambda^N H_{\mathcal{C}}$ as claimed, and the result follows.
\hfill\Halmos
\endproof

The above Theorem characterizes the structure of the optimal \emph{causal} policy for any superposition of \emph{independent} request processes, i.e. where no additional information is available about the future other than the natural history of the requests and there is no cross-information between them. 

We now analyze some examples were the above policy can be further characterized. The simplest is the Poisson case:
\begin{example}[Poisson arrivals]
  In case all processes $\Phi_i$ are Poisson with $\lambda_1\geq \ldots \geq\lambda_N$, $\lambda_i(t)\equiv\lambda_i$. By Theorem \ref{thm:optimal}, the optimal policy $\mathcal{C}^*$ is the \emph{static policy} that stores the $C$ objects with higher (average) intensities at all times. This is of course connected to the total independence property of Poisson processes and is also well known in the caching literature (see e.g. \cite{garetto2016unified}) under the name Independent Reference Model.
\end{example}

In the more general class of renewal arrival processes, by eq. \eqref{eq:stoch_intensity_renewal}, $\lambda_i(t) = \eta_i\left(t-\tau_i^-(t)\right)$. In this case, the optimal policy can be recast in the following way:

\begin{corollary}\label{cor:optimal_renewal}
If all input processes are renewal and independent, then the optimal causal policy is, at time $t$:
\begin{itemize}
  \item Rank all contents in decreasing order of the \emph{current interval hazard rates}, $\eta_i(t-\tau^{-}_i(t))$.
  \item Store in memory the first $C$ objects in the ranking. 
\end{itemize}
\end{corollary}

As we can see, the role of the hazard rates is crucial, as already identified in \cite{ferragut2016caches_sigmetrics} for timer-based caching. Different monotonicity assumptions on these hazard rates lead to completely different optimal policies, as we shall see in Section \ref{sec:threshold}.

If the processes are \emph{not} renewal or are \emph{correlated} among them, this statement will not be true in general, and we will have to keep track of the shared stochastic intensity with respect to the common history.

An issue with the optimal policy is that its main performance metric, the hit rate, cannot be computed exactly except in some special cases. We now turn our attention to an asymptotic result that characterizes the optimal policy for large scale systems and allows us to calculate asymptotic performance limits.

\section{Large-scale analysis of the optimal causal policy}\label{sec:asymptotic}

In this section we will present results concerning the asymptotic behavior of the optimal policy in a large-scale regime, where both the system load and memory size scale appropriately to infinity. For this purpose we need to introduce more structure into the problem: in particular we will assume that requests for different items come from independent renewal processes that differ in their mean intensity, but with their inter-arrival distributions sharing a common \emph{scale family}. 

The base process from which the scaled versions are derived will have \emph{unit} intensity. We reserve henceforth the notation $F_0(t)$ for the inter-arrival distribution of normalized mean $1$; $F(t)$, $\eta(t)$ are the respective age distribution and hazard rate functions, and $G_0(x)$,  $G(x)$ the corresponding distributions for observed hazard rate at, respectively, times synchronized and non-synchronized with arrivals. We now specify our scale family. 

\begin{assumption}\label{ass:renewal}
The request processes $\Phi_i$ for items $i\in \{1,\ldots, N\}$ are independent renewal processes, whose inter-arrival time distributions are given by:
 \begin{equation*}
  F_0^{(i)}(t) = F_0(\lambda_i t),
 \end{equation*}
 where $\lambda_i$ is the process intensity. Without loss of generality we will assume that the intensities are in decreasing order, i.e. $\lambda_1 \geq \dots \geq \lambda_N > 0$; equivalently, items are ordered in decreasing popularity. 
\end{assumption}
Invoking \eqref{eq:age_distribution} and \eqref{eq:def_hazard_rate}, we can express the age distributions and hazard rates of the scaled processes in terms of the base (unit mean) process:
\begin{subequations}\label{eq:scaling_equations}
  \begin{align}
     F^{(i)} (t) &= \lambda_i \int_0^t (1 - F_0(\lambda_i s))ds  = \int_0^{\lambda_i t}(1 - F_0(u))du = F(\lambda_i t);\\
    \eta^{(i)}(t) &= \frac{f_0^{(i)}( t)}{1-F_0^{(i)}(t)} = \frac{\lambda_i f_0(\lambda_i t)}{1-F_0(\lambda_i t)} =    \lambda_i \eta(\lambda_i t).
  \end{align}
\end{subequations}

We can carry out an analogous calculation for the observed hazard rate introduced in Section \ref{ssec:renewal}. Denote by $X$ the observed hazard rate for the base process, and let $G(x)$ be its distribution at a time non-synchronized with arrivals. Analogously we denote by 
$X_i$, $G^{(i)}(x)$ the corresponding objects for the scaled process with intensity $\lambda_i$. We have 
\begin{equation*}
  G^{(i)}(x) := \p\left(X_i\leqslant x\right) = 
  \p\left(\eta^{(i)}\left(-\tau_0^{(i)}\right)\leqslant x\right) = \p\left(\eta\left(-\lambda_i \tau_0^{(i)}\right) \leqslant x/\lambda_i \right).
\end{equation*}
The key observation is that  $-\lambda_i \tau_0^{(i)} \sim F$, the age distribution of the base process: indeed, 
\[
\p\left(-\lambda_i \tau_0^{(i)} \leq t\right) =  
\p\left(-\tau_0^{(i)}\leq t/\lambda_i\right) = F^{(i)}\left(t/\lambda_i\right) = F(t).
\]
Therefore, $\eta(-\lambda_i \tau_0^{(i)}) \sim G$ and we have:
\begin{equation}\label{eq:G_scaling}
  G^{(i)}(x)= G\left(x/\lambda_i\right).
\end{equation}
If, instead, we look at the observed hazard rate distribution at a time synchronized with arrivals, we can derive analogously the relationship:
\begin{equation}\label{eq:G0_scaling}
  G_0^{(i)}(x)= G_0\left(x/\lambda_i\right).
\end{equation}

\subsection{Threshold characterization of the optimal policy}

Under Assumption \ref{ass:renewal}, the optimal policy of  Corollary \ref{cor:optimal_renewal} can be recast as follows: at any given time (e.g., $t=0$), we have a sample of random variables $X_i:=\eta_i(-\tau_0^{(i)}), \, i=1,\ldots,N$, each representing the stochastic intensity (observed hazard rate for the current interval) of the $i-$th request process. Since the processes are independent, the $X_i$ are independent but \emph{non-identically} distributed, in fact $X_i\sim  G^{(i)}$. 

According to  Corollary \ref{cor:optimal_renewal}, an item will be stored in memory at time $t=0$ if and only if its observed hazard rate is one of the $C$-highest. An alternative
way to cast this optimal policy is to consider the \emph{empirical distribution} of the observed hazard rates, defined by:
\begin{equation}\label{eq:def_empirical_hazrates}
  \widehat{G}_N(x) = \frac{1}{N} \sum_{i=1}^N \mathbf{1}_{\{X_i\leq x\}}.
\end{equation}
With the above definition, object $i$ is locally stored at time $t=0$ if and only if $\widehat{G}_N(X_i) \geqslant \frac{N-C}{N} = 1-\frac{C}{N}$. 

Introducing the \emph{quantile function} (inverse of the empirical distribution) 
\begin{align}
    \widehat{Q}_N(p):=\inf\left\{x: \widehat{G}_N(x)\geqslant p\right\},\quad p\in[0,1],
\end{align}
the condition for storing item $i$ in memory may be expressed as: 
\begin{equation}\label{eq:threshold_definition}
  X_i\geqslant \widehat{\theta}_N :=  \widehat{Q}_N  \left(1-\frac{C}{N}\right).
\end{equation} 
The random variable $\widehat{\theta}_N$ acts as a \emph{threshold} in the OHR that determines the items to optimally store in local memory.
Equivalently, $\widehat{\theta}_N = Y_{N-C}$ where $\{Y_i:i=1,\ldots,N\}$ are the order statistics of the sample $\{X_i, i=1,\ldots,N\}$.

We will pursue asymptotic results for a system with a very large catalog size $N$. If we can find a suitable limit for the empirical distribution $\widehat{G}_N(x)$ as $N\to\infty$, and the memory size scales linearly with $N$, such that $\frac{C}{N}\to_N c$, then the quantile $\widehat{\theta}_N$ should approach a limit: i.e., the large-scale behavior should resemble a policy with a fixed threshold in OHR. 

Since the $X_i$ are not identically distributed, we cannot use classical Glivenko-Cantelli arguments for empirical distribution convergence. Nevertheless, we will show that a suitable limit arises  under Assumption \ref{ass:renewal} if the system load, as defined by the request intensities, also scales appropriately with $N$.

\subsection{Scaling of the request intensities}

We will construct a sequence of systems indexed by $N$, the number of arrival streams or, in other words, items in its catalog. Denote by $\{\lambda_i^{(N)}\}_{i=1}^N$ the arrival rates of the system of size $N$, with the above convention that $\lambda_1^{(N)}\geqslant \dots \geqslant \lambda_N^{(N)}>0$. We may think of these points as a discrete measure on the axis $\lambda \in \mathbb{R}_+$ of possible process intensities. Normalizing this measure to total unit mass we may write its distribution function:
\begin{equation}\label{eq:def_LN}
  L_N(\lambda) := \frac{1}{N}\sum_{i=1}^N \mathbf{1}_{\left\{\lambda_i^{(N)}\leqslant \lambda\right\}}.
\end{equation}
The total arrival rate for system of size $N$ can then be expressed as:
\begin{equation*}
  \lambda^N = \sum_i \lambda_i^{(N)} = N\int_0^\infty \lambda L_N(d\lambda).
\end{equation*}
Our limit theorems will assume that this family of discrete distributions has a limit with $N\to \infty$:
\begin{assumption}\label{ass:weak_convergence}
  There exists a fixed distribution $L$, with no atoms at $\lambda=0$, such that $L_N \Rightarrow L$ when $N\to\infty$. Here $\Rightarrow$ denotes usual weak convergence of probability distributions.
\end{assumption}
To gain some intuition on the above condition we turn to the following important example. 

\begin{example}[Scaling for Zipf popularities]\label{ex:zipf}
A widely used model of popularity among different items is the so-called Zipf distribution, where the request rates $\lambda_i\propto i^{-\beta}$, where $\beta\geq 0$ is the tail parameter of the Zipf law. Note that $\beta=0$ corresponds to a uniform distribution, whereas as for large $\beta$, relative popularities decay very fast.

In order to model this in our setting we can use the following arrival rates for the $N-$th system:
  \begin{equation*}
    \lambda_i^{(N)} = \left(\frac{N}{i}\right)^\beta,
  \end{equation*}
  with $\beta\geq 0$. Under this scaling, the least  popular object has intensity $1$ for all $N$. As $N$ grows, larger intensities are included in the mix.  Now, for any $\lambda>1$, we have:
  \begin{eqnarray*}
    1-L_N(\lambda) = \frac{1}{N}\sum_{i=1}^N \mathbf{1}_{\left\{\left( \frac{N}{i}\right)^\beta > \lambda \right\}} = \frac{1}{N}\sum_{i=1}^N \mathbf{1}_{\left\{i<\frac{N}{\lambda^{1/\beta}} \right\}} = \frac{1}{N} \left\lfloor \frac{N}{\lambda^{1/\beta}} \right\rfloor \mathop{\longrightarrow}_{N\to\infty} \lambda^{-1/\beta}.
  \end{eqnarray*}
Since the above convergence is pointwise and the limit is continuous, we have $L_N(\lambda)\Rightarrow L(\lambda)$ given by:
 \begin{equation}\label{eq:zipf_phi}
  L(\lambda) = 1 - \lambda^{-1/\beta} \quad \text{for }\lambda\geq 1.
 \end{equation}
In the limit the popularities follow a continuous distribution over $[1,\infty)$, namely a standard Pareto distribution with tail parameter $1/\beta$. 

If $\beta \geq 1$, i.e. the popularities are light-tailed, some objects are extremely more popular than others; in this case $L$ does not have a finite mean. If instead $0<\beta<1$, where popularities are heavy-tailed and thus more homogeneous, $L$ has finite mean $1/(1-\beta)$. For $\beta=0$, the system degenerates into every object having the same popularity, and thus $L$ is the step function at $\lambda=1$.

It is worth observing also that the total arrival rate of the $N$-th system satisfies:
\begin{equation*}
  \lambda^N = \sum_{i=1}^N \lambda_i^{(N)} = N^\beta \sum_{i=1}^N \frac{1}{i^\beta} =:N^\beta S_N(\beta),
\end{equation*}
where $S_N(\beta)$ is the generalized harmonic series partial sum. Using the well known equivalents for this series, we have that:
\begin{equation*}
  \lambda^N = \begin{cases}
      O(N^\beta) & \text{if } \beta>1,\\
      O(N \log N) & \text{if } \beta=1, \\
      O(N) & \text{if } \beta<1.
  \end{cases}
\end{equation*}
In particular, with our scaling, the total arrival rate $\lambda^N \to \infty$ as $N\to \infty$, albeit at different rates depending on the tail parameter $\beta$.
\end{example}

\subsection{Asymptotic behavior of the optimal causal policy}

We now return to our family of systems indexed by $N$, with renewal requests coming from a common scale family with base distribution $F_0$ (Assumption \ref{ass:renewal}), and incorporate the scaling in Assumption \ref{ass:weak_convergence} on the arrival rates  $\{\lambda_i^{(N)}\}_{i=1}^N$. 
Namely, their empirical distribution $L_N(\lambda)$ in \eqref{eq:def_LN}  has a weak limit $L(\lambda)$. 

If we look at each system at a fixed time $t=0$, we obtain a sample of the current observed hazard rates $\{X^{(N)}_i:i=1,\ldots,N\}$. Considered collectively for all $N \geqslant 1$, these random variables 
constitute a \emph{triangular array}; without loss of generality we may assume they are all defined in a common probability space $(\Omega, \mathcal{F}, \p)$.

For each $N$ we can define the random function $\widehat{G}_N$ by eq. \eqref{eq:def_empirical_hazrates}. The main result below concerns the asymptotic behavior of these empirical hazard rate distributions. 

\begin{theorem}\label{thm:main_thm}
  Consider a family of local memory systems, indexed by $N$, with request processes satisfying Assumption \ref{ass:renewal}, and with intensities $\{\lambda_i^{(N)}\}_{i=1}^N$ satisfying Assumption \ref{ass:weak_convergence}. Then:
  \begin{equation*}
    \widehat{G}_N \mathop{\Longrightarrow}_{N\to\infty} G_\infty  \quad \p-\text{a.s.}
  \end{equation*}
  where the function $G_\infty$ is given by:
  \begin{equation}\label{eq:def_Ginf}
    G_{\infty}(x):=\int_0^\infty G(x/\lambda) L(d\lambda).
  \end{equation}
  Moreover, assume that the memory size of the $N$-th system satisfies $\frac{C_N}{N} \mathop{\rightarrow}_N c$, with $0 < c \leqslant 1$. \\ Then, if $1-c$ is a continuity point of the quantile function $Q_\infty = G^{-1}_\infty$, the random threshold $\widehat{\theta}_N$ defined by eq. \eqref{eq:threshold_definition} converges $\p-$almost surely to $\theta^*= Q_\infty(1-c)$.
\end{theorem}

\proof{Proof.} 
The proof begins by computing the expected value of the random function $\widehat{G}_N$:
\begin{equation}\label{eq:gbarN}
  \overline{G}_N(x) := \E{\widehat{G}_N(x)} = \E{\frac{1}{N} \sum_{i=1}^N \mathbf{1}_{\left\{X_i^{(N)}\leq x\right\}}} = \frac{1}{N} \sum_{i=1}^N G^{(i)}(x) = \frac{1}{N} \sum_{i=1}^N G \left(\frac{x}{\lambda_i^{(N)}}\right),
\end{equation}
where we have applied the scaling in \eqref{eq:G_scaling}. $\overline{G}_N$ is a (deterministic) distribution function in $x \geq 0$. Using the definition of $L_N$, we can rewrite the above as:
\begin{equation*}
  \overline{G}_N(x) = \int_0^\infty G(x/\lambda) L_N(d\lambda).
\end{equation*}

We first show that, as distribution functions, $\overline{G}_N \Rightarrow G_\infty$ as ${N\to\infty}$. To do so, it is convenient to interpret \eqref{eq:gbarN} as follows: consider a pair $(X,\Lambda_N)$ of independent random variables in $\mathbb{R}_+$, with respective distribution functions $G(x)=P(X\leq x)$, $L_N(\lambda)=P(\Lambda_N \leq \lambda)$. Then $\overline{G}_N(x)$ is the distribution of the product $Z_N = X \Lambda_N$. Indeed, 
\begin{align*}
P(Z_N\leq x) = \sum_{n=1}^N P(X\Lambda_N  \leq x | \Lambda_N = \lambda_i^{(N)})P(\Lambda_N = \lambda_i^{(N)})  = \frac{1}{N} \sum_{n=1}^N P(X\leq x/\lambda_i^{(N)}) =  \overline{G}_N(x). 
\end{align*}
Now, consider the limit in distribution of the pair $(X,\Lambda_N)$. Due to their independence, 
and $L_N\Rightarrow L$, the limit corresponds to the distribution of $(X,\Lambda)$, a vector of independent random variables with marginal distributions $G(x)$, $L(\lambda)$. 
By continuity of the map $(x,\lambda) \mapsto x \lambda$, we conclude that  $Z_N = X \Lambda_N \stackrel{d}{\rightarrow} Z = X\Lambda$. It remains to compute the distribution function of the latter product:
\begin{align*}
    P(Z \leq x) &= \int_{\mathbb{R}_+^2} \mathbf{1}_{\{\xi \lambda \leq x\}}L(d\lambda)G(d\xi)  = \int_0^\infty L(d\lambda)  \int_0^\infty \mathbf{1}_{\{\xi\leq x/\lambda\}}G(d\xi)
      = \int_0^\infty G(x/\lambda) L(d\lambda).
\end{align*}
We conclude that $\overline{G}_N\Rightarrow G_\infty$ as $N\to\infty$. Equivalently, we have pointwise convergence  $\overline{G}_N(x) \to G_\infty(x)$ at any  continuity point $x\geqslant 0$ of $G_\infty$. 

To relate the mean function $\overline{G}_N(x)$ to the stochastic one  $\widehat{G}_N(x)$, we now resort to \cite[Theorem 2.1]{shorack1979weighted}, a generalization of the classical Glivenko-Cantelli Theorem for empirical distributions for random variables that are not identically distributed. The theorem states that, in the probability space $(\Omega,\mathcal{F},\p)$ where the triangular array $\{X^{(N)}_i: i=1,\ldots, N, N\geq 1\}$ is defined, we have: 
\begin{equation}\label{eq:shorack}
\left\|\widehat{G}_N-\overline{G}_N\right\|_\infty\mathop{\longrightarrow}_{N\to\infty}0\quad \p\text{-a.s.}
\end{equation}
In particular, with $\p$ probability one, $|\widehat{G}_N(x)-\overline{G}_N(x)|\to 0$ as $N\to\infty$ for all $x\geqslant 0$. Combining this with the previously obtained pointwise convergence of $\overline{G}_N(x)$, we conclude that, $\p$ almost surely: $\widehat{G}_N(x) \to G_\infty(x)$ at any  continuity point $x\geqslant 0$ of $G_\infty$, and thus $\widehat{G}_N \mathop{\Rightarrow} G_\infty$.

Finally, convergence of $\widehat{\theta}_N$ now follows from  the fact that convergence in distribution implies convergence of quantiles: specifically \cite[Lemma 21.2]{vaart1998asymptotic}, $\widehat{G}_N \Rightarrow G_\infty$ is equivalent to $\widehat{Q}_N(p)\to Q_\infty(p)$ for all continuity points $p\in[0,1]$ of $Q_\infty=G_\infty^{-1}$. See also Proposition \ref{prop:conv_quant} in Appendix \ref{sec:quantiles}.
\hfill\Halmos
\endproof

We have thus established that the optimal causal local memory policy converges, in the large scale limit, to a deterministic threshold policy in the stochastic intensities (observed hazard rates).  Specifically, 
when memory scales as a fraction $c$ of the catalog, a large scale system should store, 
at any given time, the items whose current OHR exceeds a threshold $\theta^*$, chosen as the $1-c$ quantile of the limit distribution $G_\infty$ in \eqref{eq:def_Ginf}.

\section{Asymptotic optimal performance}\label{sec:universalbound}

In this section we analyze the \emph{performance} achieved by the optimal policy in the large scale limit. As discussed in Section \ref{sec:optimal}, performance in local memory systems is measured by the \emph{hit probability} $H_\mathcal{C}$, or equivalently the \emph{hit rate} $h_\mathcal{C}=\lambda^N H_\mathcal{C}$. 

We will find it more convenient to derive expressions for the complementary 
\emph{miss probability} $M_\mathcal{C}=1 - H_\mathcal{C}$, and the \emph{miss rate} $m_\mathcal{C}=\lambda^N M_\mathcal{C}=\lambda^N - h_\mathcal{C}$. Referring back to eq. \eqref{eq:hit prob decomp}, we can write the following expressions for these  quantities in a system of size $N$ \footnote{We now make explicit the dependence on $N$, since we will study the limits as $M\to \infty$.}:
\begin{subequations}
     \begin{align}\label{eq:miss prob decomp}
      M_{\mathcal{C}}^{(N)} &= \sum_{i=1}^N \frac{\lambda^{(N)}_i}{\lambda^N}\p^0_{\Phi_i}\left(i\not\in \mathcal{C}(0)\right). \\
      \label{eq:miss rate}
      m_{\mathcal{C}}^{(N)} &= \sum_{i=1}^N \lambda^{(N)}_i\p^0_{\Phi_i}\left(i\not\in \mathcal{C}(0)\right).       
\end{align}
\end{subequations}

The key challenge in the evaluation of the above formulas is to compute 
$\p^0_{\Phi_i}\left(i\not\in \mathcal{C}(0)\right)$, i.e. the probability that an incoming arrival does not find its file in local memory. In the optimal policy $\mathcal{C}^*$ for fixed $N$, the determining condition is whether the OHR $X^{(N)}_i$ of the requested item finds itself among the $C$-largest. To correctly evaluate this probability, however, we must recognize an asymmetry: the competing  OHRs $X^{(N)}_j, \ j\neq i$  are being evaluated at a time \emph{not} synchronized with their arrivals, and thus follow independent distributions $G(\cdot/\lambda^{(N)}_j)$; for the item in question we must use the distribution $G_0(\cdot/\lambda^{(N)}_i)$ that applies to synchronized sampling of the OHR.

For this reason, to analyze the comparison it is convenient to define the empirical distribution 
\begin{align}\label{eq:ghat minus i}
\widehat{G}_N^{(-i)}(x)=\frac{1}{N-1}\sum_{j\neq i} \mathbf{1}_{\left\{X_j^{(N)}\leqslant x\right\}}
\end{align}
of the OHRs of non-requested items, and express the \emph{miss} condition as follows:
\begin{align*}
    i\not\in \mathcal{C}(0) &\Longleftrightarrow
\sum_{j\neq i}\mathbf{1}_{\left\{X_j^{(N)} >  X_i^{(N)}\right\}}\geqslant C \\
 &\Longleftrightarrow 
\sum_{j\neq i}\mathbf{1}_{\left\{X_j^{(N)} \leq  X_i^{(N)}\right\}}\leqslant N-1-C. \\
& \Longleftrightarrow\widehat{G}^{(-i)}_N\left(X_i^{(N)}\right)\leqslant 1 - \frac{C}{N-1} = :p_N.
\end{align*}
The above equivalence assumes there are no ties in the comparison of OHRs; to simplify the analysis to follow, we will make this a standing assumption in this section. 

\begin{assumption}\label{ass:no_ties}
For each $N\geqslant 1$ and $1\leqslant i\neq j \leqslant N$, $\p_{\Phi}^0\left(X^{(N)}_i=X^{(N)}_j\right)=0$.  
\end{assumption}

Returning to \eqref{eq:miss rate}, we may now express the performance criterion of the optimal policy as:
\begin{align}\label{eq:miss rate N}
m_{\mathcal{C}^*}^{(N)} &= \sum_{i=1}^N \lambda_i^{(N)}\p^0_{\Phi_i}\left(\widehat{G}^{(-i)}_N\left(X_i^{(N)}\right)\leqslant p_N\right) = 
\sum_{i=1}^N \lambda_i^{(N)} \p^0_{\Phi_i}\left(X_i^{(N)}\leqslant \widehat{Q}^{(-i)}_N(p_N)\right),
\end{align}
where $\widehat{Q}^{(-i)}_N$ is the inverse (quantile) function of $\widehat{G}^{(-i)}_N$.

We now outline, informally, the essence of the analysis that follows. Since the empirical distribution $\widehat{G}^{(-i)}_N$ corresponds to $N-1$ OHRs which are sampled at a non-synchronized point, its asymptotic behavior under the assumed scaling should follow the conclusions of Theorem \ref{thm:main_thm}, i.e. converge to the distribution $G_\infty$. 
Under $C/N\to c$ we have $p_N \to 1-c$, so we should have $\widehat{Q}^{(-i)}_N(p_N) 
\to Q_\infty(1-c) =\theta^*$ as $N \to \infty$. 

This leads us to consider the approximate formula 
\begin{align*}
m_{\mathcal{C}^*}^{(N)} &\approx \sum_{i=1}^N \lambda_i^{(N)} \p^0_{\Phi_i}\left(X_i^{(N)}\leqslant 
\theta^*\right) =  \sum_{i=1}^N \lambda_i^{(N)}G_0\left(\theta^*/\lambda_i^{(N)}\right); 
\end{align*}
in the last equality we have invoked the distribution of $X_i^{(N)}$, \emph{synchronized} with the arrival under the Palm probability $\p^0_{\Phi_i}$. Of course, the approximation is not a rigorous step, we have taken the limit in the quantile but not in the rest of the formula \eqref{eq:miss rate N}; nevertheless it will help us arrive at the right conjecture. 

Invoking the distribution of intensities in \eqref{eq:def_LN} we may express the approximation as
\begin{equation}\label{eq:miss_rate_approximation}
\frac{m_{\mathcal{C}^*}^{(N)}}{N} 
\approx   \frac{1}{N} \sum_{i=1}^N \lambda_i^{(N)}
G_0\left(\theta^*/\lambda_i^{(N)}\right)
= \int_0^\infty \lambda G_0(\theta^*/\lambda) L_N(d\lambda).
\end{equation}
Under Assumption \ref{ass:weak_convergence}, $L_N\Rightarrow L$. The integrand function 
$\lambda G_0(\theta^*/\lambda)$ is bounded by $\theta^*$, invoking \eqref{eq:lipschitz}; if it is continuous, the  right-hand side will converge to the corresponding integral with the limit distribution $L(\lambda)$. This formula for the asymptotic optimal miss rate is the main conjecture we will prove in Theorem \ref{thm:asymp_miss_rate} below, making rigorous the preceding approximate reasoning. 

We will require some additional technical conditions: 
\begin{itemize}
\item For continuity of $\lambda G_0(\theta^*/\lambda)$, it will suffice to guarantee it at the atoms of the measure $L(d\lambda)$, i.e. the discontinuities of $L(\lambda)$; denote this set by $D_L$, and recall by  Assumption \ref{ass:weak_convergence} that $0\not\in D_L$. Analogously, denote by $D_{G_0}$ the set of discontinuities of $G_0$, and by $D_L\cdot D_{G_0}:=\left\{\lambda x: \lambda\in D_L,x\in D_{G_0}\right\}$; we will require for our limit theorem that $\theta^*\not\in $$D_L\cdot D_{G_0}$.
\item To carry out the limit result for $m_{\mathcal{C}^*}$ in \eqref{eq:miss rate N}, we will treat separately two cases in regard to the \emph{uniform integrability} of the family of popularity distributions $\left\{L_N\right\}_{N\geqslant 1}$.
The family is uniformly integrable if:
\begin{equation*}
\forall\,\epsilon>0,\,\exists\,K\geqslant 0,\text{ such that } \sup_{N\geqslant 1}\frac{1}{N}\sum_{\lambda_i^{(N)}\geqslant K}\lambda_i^{(N)}\leqslant\epsilon.
\end{equation*}
\end{itemize}
Uniform integrability is a standard assumption (see e.g. \cite{billingsley1999}, Section 1.3) that connects convergence in distribution with convergence of the first moment. In particular, if 
 $L_N\Rightarrow L$ and the $\left\{L_N\right\}$ are uniformly integrable, then 
\begin{equation}\label{eq:unif int}
\int_0^\infty \lambda L_N(d\lambda) 
\mathop{\longrightarrow}_{N\to\infty} \int_0^\infty \lambda L(d\lambda) < \infty. 
\end{equation}
A partial converse also holds; if  $L_N\Rightarrow L$, with first moments satisfying \eqref{eq:unif int}, the family must be uniformly integrable. 

\begin{example}
    The Zipf family considered in Example \ref{ex:zipf} has continuous limit $L$, and thus satisfies all necessary continuity assumptions. Furthermore, the family $L_N$ is uniformly integrable for $\beta<1$; for $\beta \geq  1$, the limit distribution has infinite mean.  
\end{example} 

We are now ready to state our main performance result.

\begin{theorem}[Asymptotic miss rate - Uniformly integrable case]\label{thm:asymp_miss_rate}
Suppose $\left\{L_N\right\}_{N\geqslant 1}$ is uniformly integrable. Let $c\in(0,1]$ be such that $1-c$ is a continuity point of $Q_\infty$, and that $\theta^*:=Q_\infty(1-c)\notin D_L\cdot D_{G_0}$. If $C/N\to c$ as $N\to\infty$, then
\begin{equation}\label{eq:asymp_miss_rate}
\lim_{N\to\infty} \frac{m_{\mathcal{C}^*}^{(N)}}{N} 
=\int_0^\infty\lambda G_0\left(\frac{\theta^*}{\lambda}\right)L(d\lambda).
\end{equation}
\end{theorem}
The proof is given in the Appendix. 

The result states that for a uniformly integrable scaling of popularity distributions, under some regularity assumptions, the miss rate of the optimal policy scales linearly with $N$, with a proportionality constant given by the integral in \eqref{eq:asymp_miss_rate}, where the distribution function $G_0$ of the OHR under synchronous sampling appears explicitly. On the other hand, the non-synchronous distribution $G$ of OHR also influences the formula, since it determines the distribution $G_\infty$ in \eqref{eq:def_Ginf}, whose quantile is
the asymptotic threshold $\theta^*$.

\begin{corollary}[Asymptotic miss probability - Uniformly integrable case]\label{cor:asymp_miss_prob}
Under the same hypothesis as Theorem \ref{thm:asymp_miss_rate}, the optimal miss probability satisfies:
\begin{equation}\label{eq:miss prob asympt}
\lim_{N\to\infty} M_{\mathcal{C^*}}^{(N)} =  \frac{\int_0^\infty \lambda G_0\left(\frac{\theta^*}{\lambda}\right)L(d\lambda)}{\int_0^\infty \lambda L(d\lambda)}
\end{equation}
\end{corollary}
\proof{Proof.} Note that 
\[
M_{\mathcal{C^*}}^{(N)}  =\frac{m_{\mathcal{C^*}}^{(N)} }{\lambda^N} = \frac{m_{\mathcal{C^*}}^{(N)} /N}{\lambda^N/N}; 
\]
the limit of the numerator is given in Theorem \ref{thm:asymp_miss_rate}. For the denominator 
we use \eqref{eq:unif int} since $\frac{\lambda^N}{N}= \int_0^\infty \lambda L_N(d\lambda)$ and we have uniform integrability. 
\hfill\halmos\endproof

Note, regarding the formula \eqref{eq:miss prob asympt}, that the numerator is bounded by $\theta^*$; indeed we have $\lambda G_0(\theta^*/\lambda)\leq \theta^*$ from \eqref{eq:lipschitz}, and $L$ has unit mass. As the first moment of this distribution (in the denominator) becomes larger, the miss probability is smaller. This suggests that for $L$ with infinite first moment the miss probability will be zero. This is indeed true, but we must provide a separate argument since uniform integrability will not hold in this case. Our strategy will be to show that a suboptimal policy, the static one with $\mathcal{C}^s=\{1,\ldots,C\}$, already achieves vanishing asymptotic miss probability.

\begin{theorem}[Asymptotic miss probability - Non integrable case]\label{thm:asymp_miss_prob_infty}
Suppose that $\int_0^\infty \lambda L(d\lambda)=\infty$. Assume $C/N\to c\in (0,1]$ as $N\to\infty$. Then
\begin{equation}\label{eq:zero asymp_miss}
M_{\mathcal{C}^*}^{(N)}  \leqslant M_{\mathcal{C}^s}^{(N)}   \mathop{\longrightarrow}_{N\to\infty}0.
\end{equation}
\end{theorem}
\proof{Proof.} The first inequality is immediate by optimality of ${\mathcal{C}^*}$; so we focus on the static policy, where misses are certain for $i>C$, therefore:
\begin{align}
     M_{\mathcal{C}^s}^{(N)}  = \sum_{i=C+1}^N \frac{\lambda^{(N)}_i}{\lambda^N} = 
\frac{\frac{1}{N}\sum_{i=C+1}^N \lambda^{(N)}_i}{{\lambda^N}/{N}}
= \frac{ \int_{0}^{\alpha_N} \lambda L_N(d\lambda)}{ \int_{0}^{\infty} \lambda L_N(d\lambda)},
\end{align}
where $\alpha_N = L_N^{-1}(p_N)$ is the $p_N$-th quantile of $L_N$, and $p_N=1-\frac{C}{N}$.

Since $L_N\Rightarrow L$ and $p_N\to 1-c < 1$, we have that $A:=\sup_N\alpha_N<\infty$. Then
\begin{align}\label{eq:finite_num}
\limsup_{N\to\infty} \int_{0}^{\alpha_N} \lambda L_N(d\lambda) \leq \lim_{N\to\infty}\int_0^A \lambda L_N(d\lambda) = \int_{0}^{A} \lambda L(d\lambda)<\infty.
\end{align}
For the denominator, we have \cite[Thm. 3.4.]{billingsley1999}
\begin{equation}\label{eq:infinite_den}
\liminf_{N\to\infty}\int_0^\infty \lambda L_N(d\lambda)\geqslant \int_0^\infty \lambda L(d\lambda)=\infty.
\end{equation}
Together, (\ref{eq:finite_num}) and (\ref{eq:infinite_den}) imply \eqref{eq:zero asymp_miss}.
\hfill\halmos\endproof

As a comment on the above results, we note the following: when the distribution of 
item intensities is \emph{not} uniformly integrable, as e.g. in the Zipf case with 
$\beta \geqslant 1$, items with the largest intensities dominate the rest; it is thus natural that the static policy would achieve asymptotic optimality, and perfect performance, regardless of the inter-arrival distribution. Note, however, that such convergence may be slow, we refer to Section \ref{sec:simulations} for examples. 

In the more interesting case of uniform integrability, with less disparate popularities, \emph{any} causal memory management policy will incur some performance penalty; the minimum miss probability, achieved by the optimal policy, can be computed explicitly from the inter-arrival distribution characteristics and the storage fraction $c$.

\section{Threshold policies and their timer-based counterparts}\label{sec:threshold}

In Section \ref{sec:asymptotic} we found that the optimal causal policy could be expressed in terms of a threshold on the stochastic intensity (observed hazard rate). This threshold is stochastically varying, but converges in our large scale regime to a deterministic constant. 

This motivates us to look at policies that are \emph{defined} by a deterministic threshold,  i.e. where one keeps in memory items with stochastic intensity larger than $\theta$. An immediate observation is that in such policies the memory constraint $C$ would not be satisfied in a strong way at all times. Rather, we must replace it with the  soft constraint where the \emph{average memory occupation} is $C$. We now make this formal.

\subsection{Threshold policies and memory occupation}

Consider a local memory system with independent request processes $\{\Phi_i:i=1\ldots N\}$, aggregated natural history $\mathcal{F}_t$ and stochastic intensities $\lambda_i(t)$. Consider the following threshold policy:
\begin{equation}\label{eq:def_threshold_policy}
  \mathcal{C}_\theta(t) = \{i: \lambda_i(t)> \theta \} \subset \mathscr{P}(N),
\end{equation}
i.e. the store in memory all the items that have current stochastic intensity above the threshold $\theta$.

Since the $\Phi_i$ are stationary, the intensities $\lambda_i(t)$ are also stationary, and thus we have a stationary and causal policy. The memory occupation of this policy is now the stationary random process:
\begin{equation}\label{eq:memory_occupation}
  U_\theta(t) = \#\mathcal{C}_\theta(t)  = \sum_{i=1}^N \mathbf{1}_{\{\lambda_i(t)> \theta\}}.
\end{equation}
Its average can be computed as (using $t=0$ as a sampling point):
\begin{equation}\label{eq:avg_memory_occupation}
  \E{U_\theta(0)} = 
  \E{\sum_{i=1}^N \mathbf{1}_{\{\lambda_i(0)> \theta\}}} = \sum_{i=1}^N \p(\lambda_i(0)> \theta).
\end{equation}

Note that $\E{U_\theta(0)}$ is decreasing from $N$ to $0$ as $\theta$ goes from $0$ to $\infty$. In order to have a fair comparison against a fixed memory constraint, define the threshold $\theta_C$ as:
\begin{equation}\label{eq:threshold_def}
  \theta_C := \inf\left\{\theta: \sum_{i=1}^{N} \p(\lambda_i(0)> \theta) \leqslant C\right\}.
\end{equation}
Thus $\theta_C$ is the $1-C/N$ quantile of the average of the distributions of the stochastic intensities observed at time $0$. From the right continuity of the distribution functions, it is easy to check that $\E{U_{\theta_C}(0)}\leqslant C$, with equality if $1-C/N$ is a continuity point of the quantile function. We shall assume this henceforth to simplify exposition.

Since the arrival processes are independent, the random variables $\lambda_i(0)$ are independent, and thus we have the following:
\begin{proposition}\label{prop:concentration}
  Assume the above quantile is exact, in the sense that $E[U_{\theta_C}(0)] = C$. Take $C=cN$ with $0< c \leq 1$. Then for any $\varepsilon>0$:
  \begin{equation*}
    \p\left(\left| \frac{U_{\theta_C}(0)}{N} - c\right|>\varepsilon\right)\leqslant 2e^{-2N\varepsilon^2},
  \end{equation*}
  and therefore $U_{\theta_C}(0)/N\overset{\p}{\longrightarrow} c$ as $N\to\infty$.
\end{proposition}

\proof{Proof.} The proof follows from the Hoeffding's inequality applied to the independent Bernoulli random variables $\{\mathbf{1}_{\{\lambda_i(0)>\theta\}}\}$. In fact, for any $N$, and $\varepsilon>0$:
\begin{equation*}
  \p\left(\left| \frac{U_{\theta_C}(0)}{N} - c\right|>\varepsilon\right) = \p\left(\left| \sum_{i=1}^N \mathbf{1}_{\{\lambda_i(0)> \theta_C\}} - \E{\sum_{i=1}^N \mathbf{1}_{\{\lambda_i(0)> \theta_C\}}}\right| > N\varepsilon \right) \leqslant 2e^{-2N\varepsilon^2}
\end{equation*}
\hfill\halmos\endproof

From Proposition \ref{prop:concentration}, we conclude that under independent request processes, the memory occupation of the threshold policy will not deviate much from a fixed memory policy, for large $N$. We will use this to our advantage in the simulations Section.

\subsection{Asymptotic optimality of threshold policies}
We focus now on the situation studied in Sections \ref{sec:asymptotic} and \ref{sec:universalbound}, where by Assumption \ref{ass:renewal}, the request processes are renewal, from a common scale family. Furthermore, the traffic intensities satisfy Assumption \ref{ass:weak_convergence}.

In that case $\lambda_i(0) = \eta_i(-\tau_0^{(i)}) = X_i$, the observed hazard rate, whose distribution is $G^{(i)}(x)$. Therefore $\theta_C$ is the solution to the equation:
\begin{equation*}
  \sum_{i=1}^{N} \left[1-G^{(i)}(\theta_C)\right] = C.
\end{equation*}
Invoking  eq. \eqref{eq:G_scaling} for the scale family, and rearranging terms we have:
\begin{equation*}
  \overline{G}_N(\theta_C) = \frac{1}{N}\sum_{i=1}^{N} G\left(\frac{\theta_C}{\lambda_i}\right) = 1-\frac{C}{N},
\end{equation*}
with the notation introduced in \eqref{eq:gbarN}; invoking the distribution $L_N$ of traffic intensities $\lambda_i$, we may also write the above equation as:
\begin{equation}\label{eq:threshold_n}
  \int_{0}^\infty G\left(\frac{\theta_C}{\lambda}\right)L_N(d\lambda) = 1-\frac{C}{N}.
\end{equation}

We are now in position to prove the following:
\begin{proposition}\label{prop:asymptotic_threshold}
Consider a family of local memory systems, indexed by $N$, with request processes satisfying Assumptions \ref{ass:renewal} and \ref{ass:weak_convergence}. Choose the memory size of the $N$-th system as $C_N=cN$ with $0< c \leq 1$. If there exists a unique solution to $\theta^*$ satisfying:
\begin{equation*}
  G_\infty(\theta^*) = 1-c
\end{equation*}
then the sequence of thresholds $\theta_N:=\theta_{C_N}$ defined by \eqref{eq:threshold_n} satisfies $\theta_{N}\to\theta^*$.
\end{proposition}
The essential step of the proof is to show that  $\overline{G}_N\Rightarrow G_\infty$; this part is identical to the proof of Theorem \ref{thm:main_thm}. The result then follows from the convergence of quantiles.

Let us now analyze the asymptotic performance: consider an object $i$, its miss probability is given by:
\begin{equation*}
  \p_{\Phi_i}^0 \left(i\notin \mathcal{C}_{\theta_N}(0)\right) = \p_{\Phi_i}^0 \left(X_i \leqslant \theta_N \right) = G_0^{(i)}(\theta_N).
\end{equation*}
Therefore, the total miss probability for system $N$ is given by:
\begin{equation*}
  M_{\theta_{N}} = \sum_{i=1}^N \frac{\lambda_i^{(N)}}{\lambda^{(N)}} \p_{\Phi_i}^0 \left(X_i \leqslant \theta_N \right) = \sum_{i=1}^N \frac{\lambda_i^{(N)}}{\lambda^{(N)}} G_0(\theta_N/\lambda_i^{(N)}) =  \frac{\int_0^\infty \lambda G_0(\theta_N/\lambda) L_N(d\lambda)}{\int_0^\infty \lambda L_N(d\lambda)} 
\end{equation*}
We can then state a result analogous to Theorem \ref{thm:asymp_miss_rate} and Corollary \ref{cor:asymp_miss_prob} for the asymptotic performance:
\begin{proposition}\label{prop:asymptotic_performance_threshold}
Under the conditions of Theorem \ref{thm:asymp_miss_rate}, the  
asymptotic miss probability satisfies:
  \begin{equation*}
    \lim_{N\to\infty} M_{\theta_N} =  \frac{\int_0^\infty \lambda G_0\left(\frac{\theta^*}{\lambda}\right)L(d\lambda)}{\int_0^\infty \lambda L(d\lambda)}.
  \end{equation*} 
\end{proposition}
The proof is substantially easier than in the case of Theorem \ref{thm:asymp_miss_rate}, 
because the threshold $\theta_N$ is deterministic, instead of depending on the state of the remaining processes. 

Indeed, by \eqref{eq:lipschitz} and the fact that $\theta_N \to \theta^*$, we see that the function $\lambda G_0(\theta_N/\lambda)$ is uniformly bounded for all $N$. 
Taking proper account of continuity as in  Theorem \ref{thm:asymp_miss_rate}, the numerator converges as stated. For the denominator we invoke uniform integrability.  

The main conclusion of this analysis is that, under the above assumptions, a deterministic threshold policy with a soft memory constraint is asymptotically equivalent to the optimal causal policy derived in Section \ref{sec:optimal}, both in terms of memory usage and performance. 

We shall show now that there is a strong connection between threshold policies and timer based ones. 

\subsection{Connection with timer based policies: monotone hazard rates}

To end this Section, we would like to highlight a strong connection between threshold policies and \emph{timer based} ones. Timer based or time-to-live (TTL) caching has been a longstanding idea in local memory systems: here, each time a request for item $i$ occurs, a timer is started and the item is kept in memory up to timer expiration. When a new request arrives, the timer is reset. TTL policies were first analyzed in terms of point processes requests in \cite{fofack2014performance}, and the optimal timers under stationary requests were obtained in \cite{ferragut2016caches_sigmetrics}.

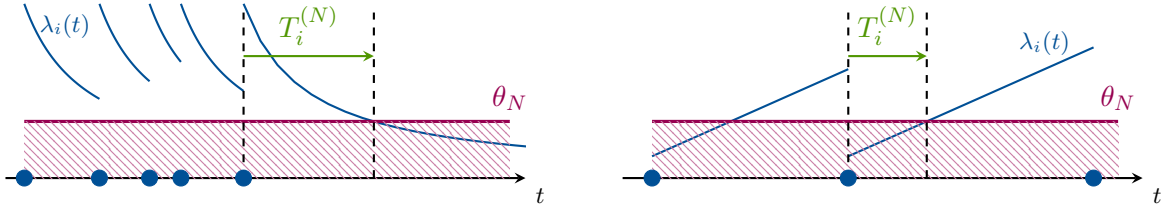
\begin{figure*}
  \begin{center}
  \input{figuras/threshold.tex}
  \end{center}
  \caption{Threshold policies in terms of timers. Right: timer-based caching for decreasing hazard rates. Left: timer-based pre-fetching for increasing hazard  rates.}\label{fig:threshold}
\end{figure*}

Consider the policy defined in \eqref{eq:def_threshold_policy} with deterministic threshold $\theta_N=\theta_{C_N}$ from eq. \eqref{eq:threshold_def}. Assume that the hazard rates are 
strictly \emph{monotonically decreasing} over the distribution domain (as in the Pareto parametric example). Assume for simplicity that the threshold $\theta_N$ lies in the range of the function $\eta_i(\cdot)$ and define: 
\begin{equation*}
  T_i^{(N)} = \eta_i^{-1}(\theta_N).
\end{equation*}
Then, since $\eta_i$ is decreasing, we have:
\begin{equation*}
  \lambda_i(t)>\theta_N \Leftrightarrow t<T_i^{(N)};
\end{equation*}
on the right of Figure \ref{fig:threshold} we illustrate this condition. Consequently, in the case of monotonically decreasing hazard rates, the threshold policy is \emph{exactly} equivalent to TTL caching.\footnote{If $\theta_N$ is not in the range of valid hazard rates, the equivalence also holds by allowing zero (or infinite) timers, i.e. items which are never (or always) stored.} This is consistent with the characterization of optimal timers in \cite{ferragut2016caches_sigmetrics}. In addition, from the asymptotic optimality of threshold policies established above, we further conclude here that \emph{TTL caching is asymptotically optimal}.

An analogous discussion is valid for monotonically \emph{increasing} hazard rates: here, the likelihood of a subsequent request is \emph{smallest} immediately after receiving one, and thus caching is not useful. Instead, it was proposed in \cite{ferragut2024prefetch} to remove the item and \emph{pre-fetch} it at an appropriate later time. 
In fact, this strategy can also be cast as a threshold policy by observing (see Figure \ref{fig:threshold}) that:
\begin{equation*}
  \lambda_i(t)>\theta_N \Leftrightarrow t>T_i^{(N)},
\end{equation*}
so the threshold policy removes the object from memory and recalls it after a time $T_i^{(N)}$; these timers coincide with the optimal timer pre-fetching policy from \cite{ferragut2024prefetch}. And once again, the asymptotical optimality of threshold policies allows us to further conclude that timer based pre-fetching is asymptotically optimal for request processes with increasing hazard rates.

\section{Parametric examples and simulations}\label{sec:simulations}

In this Section, we provide some interesting parametric examples to highlight the power of the results, in particular enabling us to compute sharp estimates for the performance of a local memory system. We must specify the inter-arrival distribution of our scale family of request processes, and the distribution of popularities. We begin with the Pareto-Zipf combination, already introduced as examples above. We then analyze a further example with increasing hazard rates.

\subsection{Pareto inter-arrival times, Zipf popularities.}

\begin{figure}
  \begin{center}
    \begin{tikzpicture}
      \begin{axis}[xlabel=$t$,
        ymin = -.05,
        ymax = 1.3,
        xmin = -.1,
        xmax=6,
        xlabel style = {at={(axis cs:6,0)},anchor=north west},
        y axis line style={->,thick},
        x axis line style={->, thick},
        xtick={2},
        xticklabels={$\frac{\alpha}{\alpha-1}$},
        yticklabels={$1-c_{\alpha,\beta}$, $1$},
        xlabel = {$x$},
        ytick={0.5,1},
        axis x line*=middle,
        axis y line*=middle,
        width=0.4\textwidth,
        height= 0.3\textwidth,
        ]
        \addplot[domain=0:6, dashed] {1};
        \addplot[domain=0:2, dashed] {.5};
        \addplot[mark=none, dashed] coordinates {(2,0) (2, .5)};
        \addplot[azulcito, thick, domain=0:2] {x/4};
        \addplot[azulcito, thick, domain=2:6] {1-1/x};
        \node[above left, azulcito] at (axis cs:6,1) {\footnotesize $G_\infty(x)$};
      \end{axis}
    \end{tikzpicture}\hspace{1em}
    \begin{tikzpicture}
      \begin{axis}[xlabel=$t$,
        ymin = -.05,
        ymax = 10,
        xmin = -.1,
        xmax=1.1,
        xlabel style = {at={(axis cs:1.1,0)},anchor=north west},
        y axis line style={->,thick},
        x axis line style={->, thick},
        xtick={0.5,1},
        xticklabels={$c_{\alpha,\beta}$,$1$},
        yticklabels={$\frac{\alpha}{\alpha-1}$},
        xlabel = {$c$},
        ytick={2},
        axis x line*=middle,
        axis y line*=middle,
        width=0.4\textwidth,
        height= 0.3\textwidth,
        ]
        \addplot[domain=0:.5, dashed] {2};
        \addplot[mark=none, dashed] coordinates {(.5,0) (.5, 2)};
        \addplot[azulcito, thick, domain=.5:1] {4*(1-x)};
        \addplot[azulcito, thick, domain=0:.5] {1/x};
        \node[above right, azulcito] at (axis cs:.5,2) {$\theta^*$};
      \end{axis}
    \end{tikzpicture}
  \end{center}
  \caption{Limit distribution $G_\infty$ and threshold behavior for the Pareto example with Zipf popularities, $\alpha=2$, $\beta=1$.}\label{fig:ginf_pareto_zipf}
\end{figure}
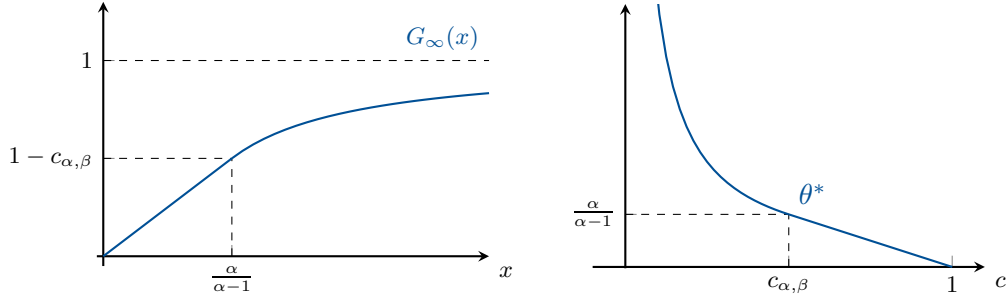

Consider first the Pareto inter-request times introduced in Section \ref{ssec:pareto}, with tail parameter $\alpha$. As we already mentioned, this family represents bursty traffic with decreasing hazard-rates. For the base (unit intensity) process we have the following distribution for  (non-synchronized) observed hazard rates: 
\begin{align*}
  G(x) &= \left(\frac{\alpha-1}{\alpha}x\right)^{\alpha-1}\text{ for } 0\leqslant x\leqslant \frac{\alpha}{\alpha-1}; \quad \quad G(x) = 1 \text{ for }x >\frac{\alpha}{\alpha-1}. 
\end{align*}
Combine it with the Zipf popularities with tail parameter $\beta$ introduced in Example \ref{ex:zipf}, with limit distribution:  
\begin{align*}
L(\lambda) &= 1-\lambda^{-1/\beta}\text{ for }\lambda\geq 1.
\end{align*}
By appropriately integrating eq. \eqref{eq:def_Ginf}, we can obtain the asymptotic distribution:
\begin{equation}\label{eq:Ginf_parteo}
  G_\infty(x) = \begin{cases}
    (1-c_{\alpha,\beta}) \left(\frac{\alpha-1}{\alpha}x\right)^{\alpha-1} & x\leqslant \frac{\alpha}{\alpha-1};\\
    1-c_{\alpha,\beta}\left(\frac{\alpha-1}{\alpha}x\right)^{-1/\beta} & x>\frac{\alpha}{\alpha-1};
  \end{cases} \quad \text{ where } \quad c_{\alpha,\beta} := \frac{(\alpha-1)\beta}{(\alpha-1)\beta + 1},
\end{equation}
a relevant parameter in the following discussion. Note that $G_\infty(x)$ is continuous and takes the value $1-c_{\alpha,\beta}$ at $x=\frac{\alpha}{\alpha-1}$. A depiction of $G_\infty(x)$ is given in Fig. \ref{fig:ginf_pareto_zipf}.

Since $G_\infty$ is strictly increasing, we can now solve for the unique asymptotic threshold $\theta^*$ as a function of the memory size $c$. Imposing $G_\infty(\theta^*)=1-c$, we get the following result:
\begin{equation*}
  \theta^*=\begin{cases}
    \left(\frac{\alpha}{\alpha-1}\right)\left(\frac{c_{\alpha,\beta}}{c}\right)^{\beta} & c\leq c_{\alpha,\beta}, \\
    \left(\frac{\alpha}{\alpha-1}\right)\left(\frac{1-c}{1-c_{\alpha,\beta}}\right)^{\frac{1}{\alpha-1}}& c\geq c_{\alpha,\beta}.
  \end{cases}
\end{equation*}

In the case $\beta<1$, it is easy to see that the measures $\{L_N\}$ are uniformly integrable, and thus we can compute the miss rate estimate from \eqref{eq:asymp_miss_rate}:
\begin{equation*}
  \int_0^\infty \lambda G_0(x/\lambda) L(d\lambda) = \begin{cases}
    (1-c_{\alpha,\beta})\left(\frac{\alpha-1}{\alpha}x\right)^{\alpha} & x\leq\frac{\alpha}{\alpha-1}; \\
    \frac{1}{1-\beta}\left[1-\alpha\beta(1-c_{\alpha,\beta})\left(\frac{\alpha-1}{\alpha}x\right)^{1-1/\beta}\right] & x\geq \frac{\alpha}{\alpha-1}.
  \end{cases}
\end{equation*}
Substituting the appropriate threshold and noting that 
$\int_0^\infty \lambda L(d\lambda) = \frac{1}{1-\beta}$, we reach the following result for the asymptotic miss probability:
\begin{equation}\label{eq:ass_miss_pareto}
  M = \begin{cases}
    1-\alpha\beta(1-c_{\alpha,\beta})\left(\frac{c}{c_{\alpha,\beta}}\right)^{1-\beta} & c\leq c_{\alpha,\beta}; \\
    \\
    (1-\beta)(1-c_{\alpha,\beta})^{-\frac{1}{\alpha-1}}(1-c)^{\frac{\alpha}{\alpha-1}}    & c\geq c_{\alpha,\beta}.
  \end{cases}
\end{equation}

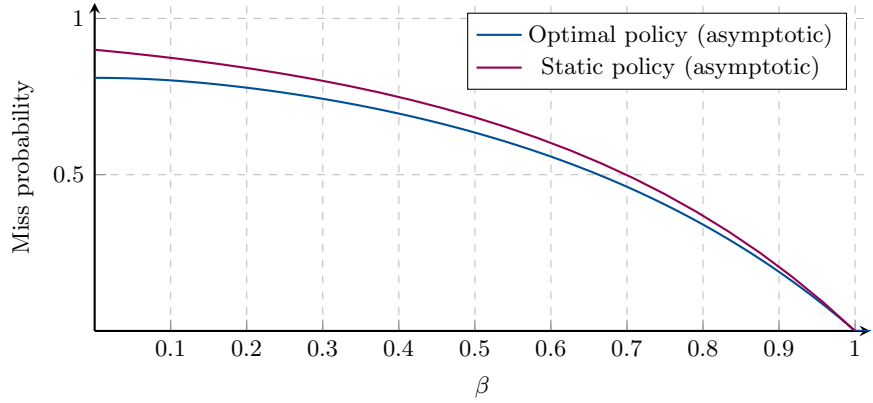
\begin{figure}
  \begin{center}
    \input{figuras/asymptotic_comparison.tex}
  \end{center}
  \caption{Asymptotic miss probability for Zipf popularities with varying parameter $\beta$. Pareto inter-request times with $\alpha=2$, $c=0.1$. Static policy added for comparison.}
  \label{fig:asymptotic_comparison}
\end{figure}

In Figure \ref{fig:asymptotic_comparison}, we plot the asymptotic miss probability of equation \eqref{eq:ass_miss_pareto} for fixed $\alpha$ and $c$, as a function of the popularity tail parameter $\beta$. We also compare it to the asymptotic miss probability of the static policy, which is derived in \cite{ferragut2016caches_sigmetrics} and is given by $1-c^{1-\beta}$ for $0\leqslant \beta \leqslant 1$. The optimal policy has the most advantage when popularities are more homogeneous, in particular as $\beta \to 0$ we have
$M\to (1-c)^\frac{\alpha}{\alpha-1}< 1-c$ for any $\alpha>1$. The gain is larger as $\alpha$ decreases: the optimal policy capitalizes on the burstiness of the incoming traffic, captured by the hazard rates.

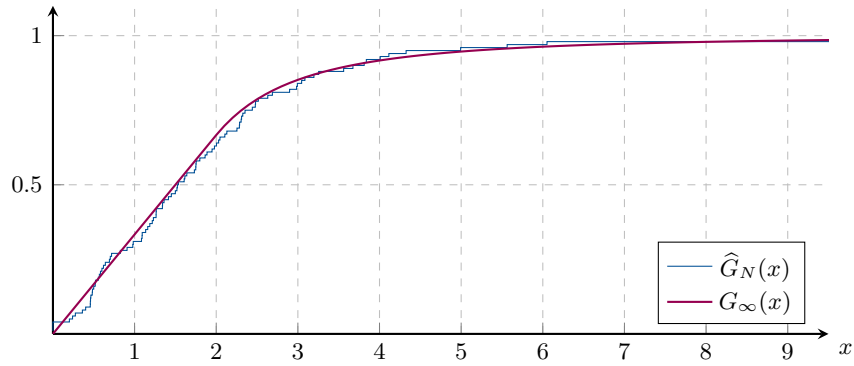
\begin{figure}
  \begin{center}
    \input{figuras/empirical_approximation.tex}
  \end{center}
  \caption{The empirical distribution of the observed hazard rates for $N=100$ and its limit $G_\infty(x)$. Pareto requests with $\alpha=2$ and Zipf popularities with $\beta=0.5$.}\label{fig:empirical_distribution}
\end{figure}

A second observation is that convergence of the empirical distribution of the OHRs is very fast. In Figure \ref{fig:empirical_distribution} we plot a simulated sample in steady state of $\widehat{G}_N(x)$ for $N=100$ items for $\alpha=2$ and $\beta=0.5$ as an example. We see the convergence to $G_\infty(x)$ as established in Theorem \ref{thm:main_thm}. 

\begin{figure}
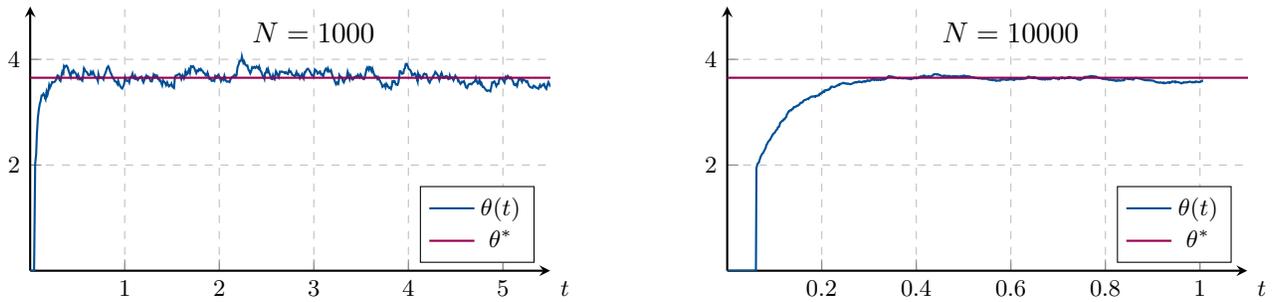

  \begin{center}
  \input{figuras/threshold_process_a2_b05_c01_N1000.tex}\hfill 
  \input{figuras/threshold_process_a2_b05_c01_N10000.tex}
  \end{center}
  \caption{Threshold evolution and large scale limit for Pareto requests $\alpha=2$, $\beta=0.5$, $c=0.1$.}\label{fig:threshold_process}
\end{figure}

Moreover, since convergence to $G_\infty$ is fast, convergence of the random threshold $\widehat{\theta}_N:=\widehat{Q}_N(1-c)\to \theta^*$ is also fast. We plot in Figure \ref{fig:threshold_process} two examples with the same parameters as above, for $N=1000$ and $10000$. The threshold process $\theta_N(t)$ is approximately constant for large $N$ around the value $\theta^*$.

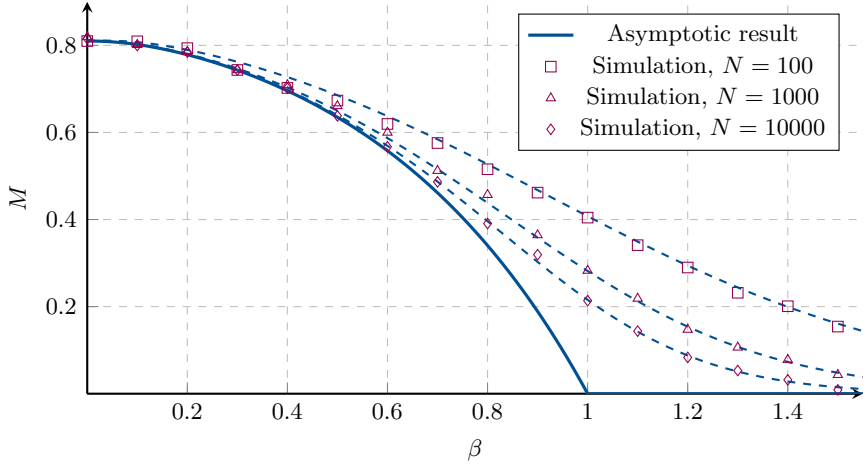
\begin{figure}
  \begin{center}
    \input{figuras/miss_probability.tex}
  \end{center}
  \caption{Miss probability of the optimal policy for Pareto requests, $\alpha=2$, $c=0.1$ and varying $\beta$. The solid curve represents the asymptotic result \eqref{eq:ass_miss_pareto}. Dots represent simulation results for different values of $N$. The dotted lines represent the finite $N$ approximation \eqref{eq:finite_N_miss}.}\label{fig:miss_probability_sim}
\end{figure}

This last observation is crucial for developing finite $N$ approximations for the performance. Due to the fact that, when $\beta \to 1$, the family of distributions $\{L_N\}$ approaches non-uniform integrability, convergence of the miss rate estimate $m_{\mathcal{C}^*}^{(N)}$, the total rate $\lambda^{(N)}$ and the miss probability $M_{\mathcal{C^*}}^{(N)}$ is slow around this value, as depicted in the simulations show in Figure \ref{fig:miss_probability_sim}.

In order to compute finite $N$ approximations, we use the intuition developed in equation \eqref{eq:miss_rate_approximation}, that is, compute the asymptotic threshold $\theta^*$ by solving $G_\infty(\theta^*) = 1-c$, and plug in this estimate in place of the random threshold $\widehat{\theta}_{N}$. Then estimate $M_{\mathcal{C^*}}^{(N)}$ as:
\begin{equation}\label{eq:finite_N_miss}
  M_{\mathcal{C^*}}^{(N)} \approx \frac{\int_0^\infty \lambda G_0(\theta^*/\lambda)L_N(d\lambda)}{\int_0^\infty \lambda L_N(d\lambda)} = \frac{1}{\lambda^{(N)}}\sum_{i=1}^N \lambda_i^{(N)} G_0(\theta^*/\lambda_i^{(N)})
\end{equation}

This turns out to be equivalent to approximate the optimal performance for that of the fixed threshold policy discussed in Section \ref{sec:threshold}, and it is numerically easy to compute, even for distributions which do not have closed form expressions as in this case. Therefore, this procedure provides \emph{sharp estimates} of the maximum achievable performance for any policy and finite $N$. An example of this approximation is depicted in dashed lines in Figure \ref{fig:miss_probability_sim}.

\subsection{Erlang inter-request times and Zipf popularities.}

We now turn our attention to a different example, where the hazard rate of the inter-request distribution is \emph{increasing}. This leads to a totally different behavior: for instance, caching is not a good idea in this setting , since upon receiving a request, a subsequent request becomes \emph{less} likely. It is actually preferable to remove the content from memory and pre-fetch it again closer to request time \cite{ferragut2024prefetch}.

Of course, the optimal memory management policy designed in Section \ref{sec:universalbound} does this automatically by keeping track of the hazard rates, and all the asymptotic results derived for it are still valid in this case.

To illustrate this behavior, we choose the inter-request times to be distributed as an Erlang distribution with $k$ stages, and appropriate means. For $k=1$, this corresponds to the Poisson process, since inter-request times become exponential. As $k$ grows, the process approaches deterministic inter-request times, and thus the traffic pattern becomes more regular.

To be precise, let the base inter-request distribution be Erlang with $k$ stages and mean $1$ (so $\lambda=k$). This corresponds to the following:
\begin{equation*}
  f_0(t) = \frac{1}{(k-1)!}k^k t^{k-1} e^{-kt}, \quad F_0(t) = 1 - \sum_{j=0}^{k-1}\frac{1}{j!}(kt)^j e^{-kt}
   \quad t\geq 0.
\end{equation*}
This yields the following nice formula for the hazard rate function:
\begin{equation*}
  \eta(t) = k \frac{\frac{(kt)^{k-1}}{(k-1)!}}{\sum_{j=0}^{k-1}\frac{1}{j!}(kt)^j} = kB(kt,k-1)
\end{equation*}
where $B(A,C)$ is the classical Erlang-B formula for the blocking probability of telephone systems. For $k>1$, this becomes strictly increasing in $t$, $\eta(0)=0$ and $\eta(t)\uparrow  k$ when $t\to\infty$.

\begin{figure}
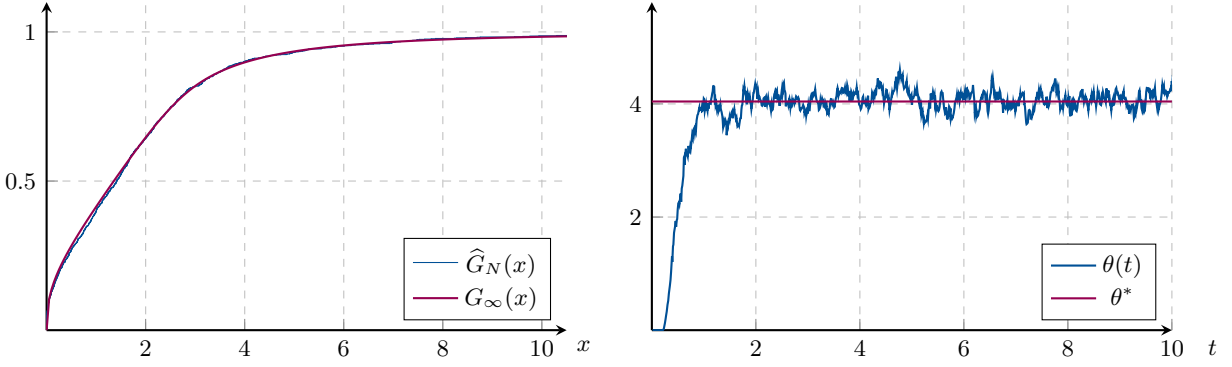

  \begin{center}
    \input{figuras/empirical_approximation_erlang.tex}
    \input{figuras/threshold_process_k4_b05_c01_N1000.tex}
  \end{center}
  \caption{Empirical distribution of the observed hazard rates for $N=1000$ and its limit $G_\infty(x)$ (left), and threshold process evolution (right). Erlang requests with $k=4$ stages, Zipf popularities with $\beta=0.5$ and $c=0.1$.}\label{fig:empirical_distribution_erlang}
\end{figure}

Unfortunately, for this distribution, there is no analytical expression for the age distribution $F$ or the observed hazard rates, and thus $G_\infty$ must be estimated numerically. We do so for the Zipf popularity distribution limit $L(\lambda)$ introduced above. We plot the resulting distribution in Figure \ref{fig:empirical_distribution_erlang}, together with a simulation of the empirical version $\widehat{G}_N(x)$ for $N=1000$, showing good convergence. Also in Figure \ref{fig:empirical_distribution_erlang}, we show the threshold $\widehat{\theta}_N$ evolution over time for the optimal policy, showing convergence to the numerically computed limit $\theta^*$.

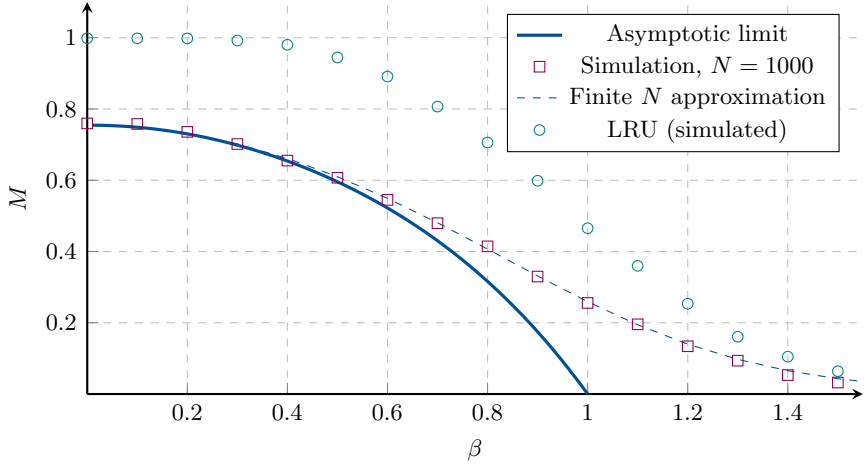
\begin{figure}[!h]
  \begin{center}
    \input{figuras/miss_probability_erlang.tex}
  \end{center}
  \caption{Miss probability of the optimal policy for Erlang requests with $k=4$ stages, $c=0.1$ and varying $\beta$. The solid curve represents the asymptotic result. Squares represent simulation results $N=1000$ and the numerical approximation from \eqref{eq:finite_N_miss} is given in dashed lines. For comparison, LRU simulations are also shown.} \label{fig:miss_probability_sim_erlang}
\end{figure}

Finally, in Figure \ref{fig:miss_probability_sim_erlang}, we plot the asymptotic miss probability of the optimal policy, numerically computed from eq. \eqref{eq:miss prob asympt}, as well as a simulation of the optimal policy for $N=1000$ as a function of the parameter $\beta$. As discussed above, convergence to the performance limit is slow when $\beta$ approaches $1$, so we also compute the performance estimate from \eqref{eq:finite_N_miss}, which also corresponds to the optimal pre-fetching policy from \cite{ferragut2024prefetch}, showing good fit. As a comparison, the classical LRU policy is also simulated for the traffic pattern. As we discussed above, classical caching performs badly, due to the regularity of the request process.

\section{Conclusions}\label{sec:conclusions}

This paper was concerned with the optimal management of local memory systems, analyzed with tools of stationary point processes. As a first contribution, we 
provided a rigorous setting for this problem, formalizing ideas in \cite{panigrahy2022upper} to show that it is optimal (in the sense of maximizing hit rates) to store in memory at every moment the items corresponding to the highest stochastic intensities. 

From this starting point, for the case of renewal request processes we analyzed the limiting behavior of the optimal policy as the catalog size $N\to \infty$, when a fixed fraction $c$ of items can be stored. Assuming item inter-request distributions come from a common scale family, with intensities that have a limit in distribution, we proved that the optimal policy amounts to comparing the stochastic intensity (observed hazard rate) of the process with a fixed \emph{threshold}, defined by the $1-c$ quantile of a certain limiting distribution function. We further characterized the asymptotic performance (miss probability) of this optimal policy.

We also analyzed optimal threshold policies for the finite $N$ case, which satisfy the memory constraint in the mean. We showed that they have the same limiting behavior as the optimal policy. Moreover, we found that for monotonic hazard rates in the inter-request distribution, these are equivalent to \emph{timer-based} policies where caching or prefetching times are determined by a hazard rate threshold. 

We present two detailed examples of our results, for the standard Zipf model of item popularities, and two cases for the inter-request process. For the bursty, decreasing hazard rate case of Pareto inter-arrival times, we obtain closed-form expressions for the optimal asymptotic threshold and the corresponding performance. We also provide sharp estimates of the optimal performance for finite $N$, validated by detailed stochastic simulations. For the regular, increasing hazard rate case of Erlang inter-arrival times, closed-form expressions are not available but we still carry out a numerical validation of the asymptotic threshold and performance. We also use this example to exhibit the significant superiority of the optimal policy in comparison to the popular LRU policy, in contrast with the bursty case where LRU has good behavior. 

For future research, we may explore the application of our machinery to characterize optimality for other kinds of request traffic: Markov-modulated Poisson processes to account for time variations, more general mixtures of traffic from different sources, and possibly heterogeneity in item sizes, as well as different popularity distributions. 

Note also that as presented, the optimal policy serves as a fundamental limit on the achievable performance, but not necessarily a practical algorithm since it requires knowledge of the relevant item intensities and hazard rate distributions. This poses the question of possibly \emph{learning} this information from the actual request data. Alternatively, to design an eviction policy which does not make explicit use of these quantities but approximates optimal performance in practice; to some extent LRU achieves this in the bursty case, but a counterpart for the increasing hazard rate is currently open.  


\input{stochsys24.bbl}

\begin{appendices}

\section{Appendix: Proof of Asymptotic Performance}\label{sec:app proof}

In this section we will provide a proof for our main result on asymptotic performance of the optimal policy, for the uniformly integrable case. Referring back to Section \ref{sec:asymptotic},
we begin by restating equation \eqref{eq:miss rate N}: 
\begin{align}\label{eq:miss rate N-app}
m_{\mathcal{C}^*}^{(N)}  &= \sum_{i=1}^N \lambda_i^{(N)}\p^0_{\Phi_i}\left(\widehat{G}^{(-i)}_N\left(X_i^{(N)}\right)\leqslant p_N\right) = 
\sum_{i=1}^N \lambda_i^{(N)} \p^0_{\Phi_i}\left(X_i^{(N)}\leqslant \widehat{Q}^{(-i)}_N(p_N)\right),
\end{align}
where $p_N=1-\frac{C}{N-1}$, $\widehat{G}^{(-i)}_N$ was defined in \eqref{eq:ghat minus i} to be the empirical distribution of the OHRs of items different from $i$, and $\widehat{Q}^{(-i)}_N$ its corresponding quantile function. 

To compute each term on the right we must invoke the distribution of $X^{(N)}_i$ under the Palm probability; it would be more convenient if the right-hand side of the inequality does not depend on $i$. For this purpose we will bound the quantiles obtained from 
$\widehat{G}^{(-i)}_N$ with those of $\widehat{G}_N$ defined in \eqref{eq:def_empirical_hazrates}, which includes all items; as usual, $\widehat{Q}_N$ is the corresponding quantile function. The key observation is that the contribution of the $i$-th item to the empirical distribution $\widehat{G}_N$ of all OHRs is at most $\frac{1}{N}$, and therefore negligible in the limit as $N$ tends to infinity. More precisely, we have:
\begin{lemma}\label{lemma:1 over N}
$\max_{1\leqslant i\leqslant N}
\left\|
\widehat{G}_N - \widehat{G}_N^{(-i)}
\right\|_\infty
\leqslant
\frac{1}{N}.$
\end{lemma}
\proof{Proof.} Recalling the definitions \eqref{eq:def_empirical_hazrates} and 
\eqref{eq:ghat minus i}, we can write:
\begin{align*}
    N \widehat{G}_N(x) -(N-1) \widehat{G}_N^{(-i)}(x) &= \sum_{j=1}^N \mathbf{1}_{\{X^{(N)}_j\leq x\}} -  \sum_{j\neq i} \mathbf{1}_{\{X^{(N)}_j\leq x\}}  = \mathbf{1}_{\{X^{(N)}_i\leq x\}};\\
    \Longrightarrow N\left(\widehat{G}_N(x)- \widehat{G}_N^{(-i)}(x)\right) &=  \mathbf{1}_{\{X^{(N)}_i\leq x\}} - \widehat{G}_N^{(-i)}(x).
\end{align*}
In the preceding identity, the right-hand side is the difference of two numbers, both in $[0,1]$; therefore
$|\widehat{G}_N(x) - \widehat{G}_N^{(-i)}(x) | \leqslant \frac{1}{N}$.  
This holds for every $x$, and for each $i$.  
\hfill\halmos\endproof

We will now make use of this Lemma to bound the distribution of 
$\widehat{G}^{(-i)}_N (X_i^{(N)})$, as required for \eqref{eq:miss rate N-app}, with 
that of $\widehat{G}_N (X_i^{(N)})$. However, this change brings a new difficulty, since 
$\widehat{G}_N(\cdot)$ and $X_i^{(N)}$ would not be independent. For this reason, to calculate the distribution we will use
a coupling argument and consider independent versions of the corresponding random variables.

\begin{proposition}[Miss rate bounds for fixed $N$]\label{prop:miss_rate_fixed_N}
Let $p_N^{\pm}=p_N\pm\frac{1}{N}$ where $p_N=1-\frac{C}{N-1}$. Let $\varphi_N^{\pm}$ be the probability distribution function of the random variable $\widehat{\theta}_N^{\pm}:=\widehat{Q}_N\left(p_N^{\pm}\right)$ under the probability $\p$. Then the miss rate is bounded from above and below by
\begin{equation}\label{eq:miss_rate_bounds}
\int_0^\infty \int_0^\infty \lambda G_0\left(\frac{\theta}{\lambda}\right)\varphi_N^-(d\theta) L_N(d\lambda)
\leqslant \frac{m_{\mathcal{C}^*}^{(N)} }{N}\leqslant \int_0^\infty \int_0^\infty \lambda G_0\left(\frac{\theta}{\lambda}\right)\varphi_N^+(d\theta) L_N(d\lambda).
\end{equation}
\end{proposition}
\proof{Proof.}
The proof is based on a coupling argument where we consider independent copies of the OHRs. The distributions of these copies mimic the distributions of the OHRs at an arbitrary time $t=0$ (that is, under the probability measure $\p$) and upon an arrival at time $t=0$ (that is, under the Palm distribution $\p_{\Phi}^0$). More precisely, consider two independent, row independent triangular arrays
\[
\left\{Y_i^N:N\geqslant 1,\,1\leqslant i\leqslant N\right\}
\text{ and }
\left\{Z_i^N:N\geqslant 1,\,1\leqslant i\leqslant N\right\},
\]
defined on a common probability space $(\tilde{\Omega},\tilde{\mathcal{F}},\tilde{\p})$ with distributions
\[
Y_i^N\sim G\left(\cdot/\lambda_i^N\right) \text{ and } Z_i^N\sim G_0\left(\cdot/\lambda_i^N\right) \text{ for all }1\leqslant i\leqslant N\text{ and }N\geqslant 1.
\]
Let
\[
\widehat{H}_N(x)=\frac{1}{N}\sum_{j=1}^N 1_{\left\{Y_j^N\leqslant x\right\}}
\text{ and }
\widehat{H}_N^{(-i)}(x)=\frac{1}{N-1}\sum_{j\neq i}1_{\left\{Y_j^N\leqslant x\right\}}.
\]
Analogously to Lemma \ref{lemma:1 over N} we have:
\begin{equation}\label{eq:infinity_norm_G}
\max_{1\leqslant i\leqslant N}
\left\|
\widehat{H}_N^{(-i)}-\widehat{H}_N
\right\|_\infty
\leqslant
\frac{1}{N}.
\end{equation}
For each $N\geqslant 1$ and $1\leqslant i\leqslant N$ the random variables
\[
Y_1^N,\ldots,Y_{i-1}^N, Z_i^N, Y_{i+1}^N,\ldots, Y_N^N
\]
are independent, so if we let $U_i^N=\widehat{H}_N^{(-i)}\left(Z_i^N\right)$, its distribution under $\tilde{\p}$ will coincide with that of $\widehat{G}_N^{(-i)}\left(X_i^{(N)}\right)$ under the Palm probability  $\p_{\Phi_i}^0$, i.e.:
\begin{equation}
    \label{eq:palm equals tilde}
\p^0_{\Phi_i}\left(\widehat{G}^{(-i)}_N(X_i^{(N)})\leqslant p_N\right)
=\tilde{\p}\left(U_i^N\leq p_N\right).
\end{equation}
To bound the probability on the right-hand side, observe that from \eqref{eq:infinity_norm_G} we have:
\begin{equation}\label{eq:in_between}
\left\{\widehat{H}_N\left(Z_i^N\right)\leqslant p_N^{-} \right\}
\subset
\left\{U_i^N\leqslant p_N\right\}
\subset
\left\{\widehat{H}_N\left(Z_i^N\right)\leqslant p_N^{+}\right\};
\end{equation}
the proof proceeds by computing the probabilities of the extreme sets in \eqref{eq:in_between}, noting that under the constructed $\tilde{\p}$, 
$Z^N_i$ is independent of $\widehat{H}_N$ (a function of $Y_1^N,\ldots Y_N^N$).

Focusing momentarily on the upper bound, we can write:
\begin{equation}\label{eq:quantile H}
\tilde{\p}\left(\widehat{H}_N\left(Z_i^N\right)\leq p_N^+\right) = 
\tilde{\p}\left(Z_i^N \leqslant \widehat{Q}_N^H (p_N^+)\right),
\end{equation}
where $\widehat{Q}_N^H$ denotes the quantile function of $\widehat{H}_N$; by construction, 
the distribution of $\widehat{Q}_N^H(p_N^+)$ under $\tilde{\p}$ coincides with that of 
$\widehat{\theta}^+_N = \widehat{Q}_N(p_N^+)$ under $\p$, which by hypothesis is denoted by $\varphi_N^+(\theta)$. 
Since $Z^N_i$ is independent and has distribution $G_0(\cdot/\lambda^N_i)$, we can compute the right hand-side of \eqref{eq:quantile H} to be
\[
\int_0^\infty \tilde{\p}\left(Z_i^N \leqslant \theta \right) \varphi_N^+(d\theta) = 
\int_0^\infty G_0\left(\frac{\theta}{\lambda_i^N}\right)\varphi_N^+(d\theta). 
\]
Multiplying by $\lambda_i^N$ and averaging over $i$ we obtain from \eqref{eq:in_between}:
\begin{align*}
\frac{1}{N} \sum_{i=1}^n \lambda_i^N\tilde{P}(U_i^N\leqslant p_N) &\leqslant
\frac{1}{N} \sum_{i=1}^n \tilde{\p}\left(\widehat{H}_N\left(Z_i^N\right)\leq p_N^+\right) \\
& =  \frac{1}{N} \sum_{i=1}^n \int_0^\infty \lambda_i^N G_0\left(\frac{\theta}
{\lambda_i^N}\right)\varphi_N^+(d\theta) 
= \int_0^\infty \int_0^\infty \lambda G_0\left(\frac{\theta}{\lambda}\right)\varphi_N^+(d\theta) L_N(d\lambda),
\end{align*}
where the last step invokes the definition of $L_N$. An analogous calculation provides the lower bound 
\[
 \int_0^\infty \int_0^\infty \lambda G_0\left(\frac{\theta}{\lambda}\right)\varphi_N^-(d\theta) L_N(d\lambda) \leqslant
\frac{1}{N} \sum_{i=1}^n \lambda_i^N\tilde{P}(U_i^N\leqslant p_N).
 \]
 Invoking \eqref{eq:palm equals tilde} and the expression \eqref{eq:miss rate N-app} for the miss rate, we conclude the proof. 
\hfill\halmos\endproof

\subsection{Convergence of quantiles}\label{sec:quantiles}

Our next step is compute the limit as $N\to\infty$, under suitable conditions, of the lower and upper bounds in (\ref{eq:miss_rate_bounds}). We will use the convergence $L_N \Rightarrow L$ of Assumption \ref{ass:weak_convergence}, but we also need to address the convergence of $\varphi_N^\pm$, distribution of quantiles for the random function $\widehat{G}_N$. 

Recall that by Theorem \ref{thm:main_thm}, we have $\widehat{G}_N\Rightarrow G_\infty$ with $\p$ probability one, where $G_\infty$ is defined by (\ref{eq:def_Ginf}). 
A standard result we have already invoked (e.g \cite[Lemma 21.2]{vaart1998asymptotic}, is that convergence in distribution implies the convergence of quantiles at a fixed point $p$, provided the limit quantile function is continuous at $p$.  We will need a slight generalization of this property, for the case where quantiles of the sequence are evaluated at a variable point $p_N$, convergent to $p$. This is stated and proved next.

\begin{proposition}[Convergence of quantiles]\label{prop:conv_quant}
Let $\{F_N\}_{N\geqslant 1}$ be a sequence of cdfs such that $F_N\Rightarrow F$ for some cdf $F$. Let $Q_N$ and $Q$ be the quantile function of $F_N$ and $F$ respectively. Let $p\in[0,1]$ be a continuity point of $Q$ and $\left\{p_N\right\}_{N\geqslant 1}$ a sequence in $[0,1]$ such that $p_N\to p$. Then $Q_N(p_N)\to Q(p)$.    
\end{proposition}
\proof{Proof.}
In the proof we use the characterization of weak convergence in terms of the Lévy distance
\[
d_L(F, G):=\inf \{\epsilon>0 \mid F(x-\epsilon)-\epsilon \leqslant G(x) \leqslant F(x+\epsilon)+\epsilon, \forall x \in \mathbb{R}\},
\]
between cdfs $F, G: \mathbb{R} \to [0,1]$. Convergence in the Lévy distance $d_L(F_N,F)\to 0$ is equivalent to weak convergence $F_N\Rightarrow F$.

Let $\epsilon>0$. By continuity of $Q$ at $p$, there exists $\delta>0$ such that $|q-p|<\delta$ implies $|Q(q)-Q(p)|<\epsilon/2$. Let $\delta_N:=d_L(F_N,F)$. Choose $N_0$ such that $\left|p_N-p\right|<\delta/2$ and $\delta_N\leqslant \min\{\epsilon,\delta\}/2$ for all $N\geqslant N_0$.

We claim that
\begin{equation}\label{eq:in_between_Q}
Q\left(p_N-\delta_N\right)-\delta_N\leqslant Q_N(p_N)\leqslant Q\left(p_N+\delta_N\right)+\delta_N.
\end{equation}
Indeed, for the left inequality of (\ref{eq:in_between_Q}), let $x\in\mathbb{R}$ be such that $F_N(x)\geqslant p_N$. Then
\[
F(x+\delta_N)\geqslant F_N(x)-\delta_N \geqslant p_N-\delta_N,
\]
which implies $x+\delta_N\geqslant Q(p_N-\delta_N)$. Taking infimum over such $x$ we get $Q_N(p_N)\geqslant Q(p_N-\delta_N) -\delta_N$.

For the right inequality of (\ref{eq:in_between_Q}), let $x\in\mathbb{R}$ be such that $F(x-\delta_N)\geqslant p_N+\delta_N$. Then
\[
F_N(x)\geqslant F(x-\delta_N)-\delta_N\geqslant p_N,
\]
which implies $x\geqslant Q_N(p_N)$. Taking infimum over such $x$ we get $Q(p_N+\delta_N)+\delta_N\geqslant Q_N(p_N)$. This proves the claim.

Let $N\geqslant N_0$. Then
\[
Q(p)-\epsilon \leqslant Q(p) - \frac{\epsilon}{2} - \delta_N < Q(p_N-\delta_N) -\delta_N \leqslant Q_N(p_N)
\]
and
\[
Q_N(p_N)\leqslant Q(p_N+\delta_N)+\delta_N < Q(p) + \frac{\epsilon}{2} + \delta_N \leqslant Q(p) + \epsilon.
\]
Thus $\left|Q_N(p_N)-Q(p)\right| < \epsilon$ for all $N\geqslant N_0$.
\halmos\hfill\endproof

\subsection{Proof of Theorem \ref{thm:asymp_miss_rate}}

We are now in a position to complete the proof of our main result on performance. 

We first establish that $\widehat{\theta}_N^\pm:=\widehat{Q}_N(p_N^\pm) \mathop{\rightarrow}_{N\to\infty} \theta^*:=Q_\infty(1-c)$, almost surely in $\p$. This follows by invoking Proposition \ref{prop:conv_quant} at each $\omega$ in the set where $\widehat{G}_N\Rightarrow G_\infty$, which has unit probability by  Theorem \ref{thm:main_thm}; note that $p_N^\pm\to 1-c$, which by hypothesis is a point of continuity of $Q_\infty$. 

Almost sure convergence implies convergence in probability, and in distribution. Therefore 
the distributions $\varphi_N^\pm$ of the random variables $\widehat{\theta}_N^\pm$ both satisfy $\varphi_N^\pm \Rightarrow \delta_{\theta^*}$, a unit mass at the limit threshold. 
As a consequence, under Assumption \ref{ass:weak_convergence} we have convergence of the
product measures $\varphi_N^{\pm}\otimes L_N \Rightarrow \delta_{\theta^*}\otimes L$.  
We will use this fact to show that both bounds in (\ref{eq:miss_rate_bounds}) converge to the same limit, and it coincides with the right-hand side of (\ref{eq:asymp_miss_rate}).

Define for this purpose, the function $h:\mathbb{R}_+^2 \to \mathbb{R}_+$ by $h(\lambda,\theta) = \lambda G_0\left(\frac{\theta}{\lambda}\right)$. Our next claim is that 
\begin{align}\label{eq:limit hinv}
\left(\varphi_N^{\pm}\otimes L_N\right)h^{-1} \Rightarrow \left(\delta_{\theta^*}\otimes L\right)h^{-1}. 
\end{align}
This claim is established in the mapping theorem \cite[Thm. 2.7.]{billingsley1999} provided that $h$ is measurable and its set of discontinuities $D_h$ has zero measure in the limit, i.e.  $(\delta_{\theta^*}\otimes L) (D_h) = 0$.

Since the limit measure in $\theta$ is a point mass, we have: 
\begin{align}\label{eq:limit disc}
    (\delta_{\theta^*}\otimes L) (D_h) = L\left(\left\{\lambda: \frac{\theta^*}{\lambda}\in D_{G_0}\right\}\right).
\end{align}
Note also that, $G_0$ being a distribution function, its set of discontinuities is countable. 
Therefore the measure on the right of \eqref{eq:limit disc} can only be positive if there exists $\lambda_0 \in D_L$ (an atom of the measure $L$) such that $\frac{\theta^*}{\lambda_0}\in D_{G_0}$. This would imply that $\theta^* \in D_L\cdot D_{G_0}$ and is ruled out by hypothesis. Therefore the claim holds. 

To finish the theorem we must go from the convergence in distribution given in \eqref{eq:limit hinv} to a first moment condition. Here is where uniform integrability will come into play. 
It is easiest to express this via independent random variables  $\widehat{\theta}_N^{\pm}$,  
$\Lambda_N$ with respective distributions $\varphi_N^{\pm}$, $L_N$ with
$\widehat{\theta}_N^{\pm}\stackrel{d}{\rightarrow} \theta^*$ and 
$\Lambda_N \stackrel{d}{\rightarrow} \Lambda$, the latter with  distribution $L$. Condition \eqref{eq:limit hinv} is equivalent to the statement
\[
h\left(\Lambda_N,\widehat{\theta}_N^{\pm}\right) =\Lambda_N G_0\left(\frac{\widehat{\theta}_N^{\pm}}{\Lambda_N}\right) \stackrel{d}{\longrightarrow} \Lambda G_0\left(\frac{\theta^*}{\Lambda}\right)= h\left(\Lambda,\theta^{*}\right).
\]
Since $G_0\leqslant 1$ we have $h(\Lambda_N,\widehat{\theta}_N^{\pm}) \leq \Lambda_N$; now 
$\{\Lambda_N\}_{N\geqslant 1}$ is uniformly integrable by hypothesis, so the same happens with 
$\{h(\Lambda_N,\widehat{\theta}_N^{\pm})\}_{N\geqslant 1}$. Now invoke \cite[Thm. 3.5.]{billingsley1999} to obtain:
\begin{align*}
\int_0^\infty\hspace{-2.8pt}\int_0^\infty \lambda G_0\left(\frac{\theta}{\lambda}\right)L_N(d\lambda)\varphi_N^{\pm}(d\theta)
=\E{\Lambda_N G_0\left(\frac{\widehat{\theta}_N^{\pm}}{\Lambda_N}\right)} 
\mathop{\longrightarrow}_{N\to\infty} 
\E{\Lambda G_0\left(\frac{\theta^*}{\Lambda}\right)}
=\int_0^\infty \lambda G_0\left(\frac{\theta^*}{\lambda}\right)L(d\lambda).
\end{align*}
Thus, both upper and lower bounds in  (\ref{eq:miss_rate_bounds}) converge to the given limit, which concludes the proof.
\hfill\halmos\endproof

\end{appendices}

\end{document}

%% file: figuras/threshold.tex
  \begin{tikzpicture}
    \begin{axis}[xlabel=$t$,
      ymin = -.3,
      ymax = 2,
      xmin = -.3,
      xmax=8,
      xlabel style = {at={(axis cs:8,0)},anchor=north west},
      y axis line style={draw=none},
      x axis line style={->, thick},
      ticks=none,
      axis x line*=middle,
      width=0.5\columnwidth,
      height=0.25\columnwidth,
      ]
      \addplot[azulcito,thick, domain=0:1.2] {2/(1+x)};
      \addplot[azulcito,thick, domain=1.2:2] {2/(1+x-1.2)};
      \addplot[azulcito,thick, domain=2:2.5] {2/(1+x-2)};
      \addplot[azulcito,thick, domain=2.5:3.5] {2/(1+x-2.5)};
      \addplot[azulcito,thick, domain=3.5:10] {2/(1+x-3.5)};
      \node[below, azulcito] at (axis cs:0.65,2) {\footnotesize $\lambda_i(t)$};
      
      \addplot+[azulcito, solid, mark options={fill=azulcito, scale=1.5}, mark=*,only marks] coordinates {
        (0,0)
        (1.2,0)
         (2,0)
         (2.5,0)
         (3.5,0)
      };

      \draw[rojito, very thick] (axis cs:0,.65) -- (axis cs:7.75,.65) node[above]{$\theta_N$};
      \path[pattern=north west lines, pattern color=rojito!50!white] (axis cs:0,0.65) rectangle (7.75,0);
      \draw[dashed,thick] (axis cs:3.5,0) -- (axis cs:3.5,1.9);
      \draw[dashed,thick] (axis cs:5.58,0) -- (axis cs:5.58,1.9);
      \draw[thick, ->, verdecito] (axis cs:3.5,1.4) -- (axis cs:5.58,1.4) node[midway, above] {$T_i^{(N)}$};
    \end{axis}
  \end{tikzpicture}\hspace{2em}
  \begin{tikzpicture}
    \begin{axis}[xlabel=$t$,
      ymin = -.3,
      ymax = 2,
      xmin = -.3,
      xmax=5,
      xlabel style = {at={(axis cs:5,0)},anchor=north west},
      y axis line style={draw=none},
      x axis line style={->, thick},
      ticks=none,
      axis x line*=middle,
      width=0.5\columnwidth,
      height=0.25\columnwidth,
      ]
      \addplot[azulcito, thick, domain=0:2] {0.25+0.5*x};
      \addplot[azulcito,thick, domain=2:4.5] {0.25+0.5*(x-2)};
      \node[below, azulcito] at (axis cs:4,1.8) {\footnotesize $\lambda_i(t)$};

      \addplot+[azulcito, solid, mark options={fill=azulcito, scale=1.5}, mark=*,only marks] coordinates {
    (0,0)
    (2,0)
    (4.5,0)
  };
  
      \draw[rojito, very thick] (axis cs:0,.65) -- (axis cs:4.75,.65) node[above]{$\theta_N$};
      \path[pattern=north west lines, pattern color=rojito!50!white] (axis cs:0,0.65) rectangle (4.75,0);

      \draw[dashed,thick] (axis cs:2,0) -- (axis cs:2,1.9);
      \draw[dashed,thick] (axis cs:2.8,0) -- (axis cs:2.8,1.9);
      \draw[thick, ->, verdecito] (axis cs:2,1.4) -- (axis cs:2.8,1.4) node[midway, above] {$T_i^{(N)}$};
    \end{axis}
  \end{tikzpicture}

%% file: figuras/asymptotic_comparison.tex
\begin{tikzpicture}
    \begin{axis}[xlabel=$t$,
        ymin = 0,
        ymax = 1.05,
        xmin = 0,
        xmax = 1.02,
        y axis line style={->,thick},
        x axis line style={->, thick},
        xlabel = {$\beta$},
        ylabel = {Miss probability},
        axis x line*=middle,
        axis y line*=middle,
        width=0.7\textwidth,
        height= 0.35\textwidth,
        legend pos=north east,
        grid
        ]

    \addplot[color=azulcito, solid, thick] table[row sep={\\}]
    {
        \\
        0.0  0.81  \\
        0.01  0.8099190000000001  \\
        0.02  0.8096760000000001  \\
        0.03  0.8092710000000001  \\
        0.04  0.808704  \\
        0.05  0.807975  \\
        0.06  0.807084  \\
        0.07  0.806031  \\
        0.08  0.8048160000000002  \\
        0.09  0.803439  \\
        0.1  0.8019000000000002  \\
        0.11  0.800199  \\
        0.12  0.7983372977053461  \\
        0.13  0.796322791646331  \\
        0.14  0.7941640436943629  \\
        0.15  0.7918679276522128  \\
        0.16  0.7894399793079565  \\
        0.17  0.7868846577764377  \\
        0.18  0.7842055428639393  \\
        0.19  0.7814054866108218  \\
        0.2  0.7784867313503421  \\
        0.21  0.7754510028628414  \\
        0.22  0.7722995847124952  \\
        0.23  0.7690333781640768  \\
        0.24  0.7656529509081398  \\
        0.25  0.7621585769994558  \\
        0.26  0.7585502698238125  \\
        0.27  0.7548278094796816  \\
        0.28  0.7509907656455393  \\
        0.29  0.7470385167681468  \\
        0.3  0.7429702662294577  \\
        0.31  0.7387850560144084  \\
        0.32  0.7344817782976005  \\
        0.33  0.7300591852859202  \\
        0.34  0.725515897590713  \\
        0.35  0.7208504113530634  \\
        0.36  0.716061104305908  \\
        0.37  0.7111462409248186  \\
        0.38  0.7061039767936  \\
        0.39  0.7009323622899948  \\
        0.4  0.6956293456797971  \\
        0.41  0.6901927756937372  \\
        0.42  0.6846204036500054  \\
        0.43  0.67890988517577  \\
        0.44  0.6730587815731204  \\
        0.45  0.6670645608682485  \\
        0.46  0.6609245985771202  \\
        0.47  0.6546361782162052  \\
        0.48  0.64819649158287  \\
        0.49  0.6416026388266665  \\
        0.5  0.6348516283298892  \\
        0.51  0.6279403764133009  \\
        0.52  0.620865706880839  \\
        0.53  0.6136243504152812  \\
        0.54  0.6062129438352988  \\
        0.55  0.5986280292229472  \\
        0.56  0.5908660529294832  \\
        0.57  0.5829233644663496  \\
        0.58  0.5747962152872796  \\
        0.59  0.5664807574666803  \\
        0.6  0.5579730422787526  \\
        0.61  0.5492690186812061  \\
        0.62  0.540364531706868  \\
        0.63  0.531255320766019  \\
        0.64  0.5219370178618485  \\
        0.65  0.5124051457210493  \\
        0.66  0.502655115841217  \\
        0.67  0.4926822264564269  \\
        0.68  0.4824816604220663  \\
        0.69  0.4720484830197663  \\
        0.7  0.461377639683035  \\
        0.71  0.45046395364399705  \\
        0.72  0.4393021235014489  \\
        0.73  0.427886720710269  \\
        0.74  0.416212186992059  \\
        0.75  0.404272831666749  \\
        0.76  0.39206282890475663  \\
        0.77  0.3795762148991675  \\
        0.78  0.36680688495728364  \\
        0.79  0.35374859051076957  \\
        0.8  0.34039493604352933  \\
        0.81  0.3267393759363404  \\
        0.82  0.31277521122717344  \\
        0.83  0.29849558628604667  \\
        0.84  0.2838934854031637  \\
        0.85  0.2689617292890062  \\
        0.86  0.25369297148497516  \\
        0.87  0.2380796946830851  \\
        0.88  0.2221142069531521  \\
        0.89  0.20578863787582802  \\
        0.9  0.18909493457977522  \\
        0.91  0.17202485768118747  \\
        0.92  0.15456997712380638  \\
        0.93  0.13672166791749996  \\
        0.94  0.1184711057734128  \\
        0.95  0.09980926263361012  \\
        0.96  0.08072690209309286  \\
        0.97  0.061214574711965786  \\
        0.98  0.04126261321549052  \\
        0.99  0.020861127579677063  \\
        1.0  0.0  \\
        2.0  0.0  \\
        }
        ;
        \addlegendentry{Optimal policy (asymptotic)}

        \addplot[color=rojito, solid, thick, domain=0:1]
        {1-(.1)^(1-x)};
        \addlegendentry{Static policy (asymptotic)}

    \end{axis}
\end{tikzpicture}

%% file: figuras/empirical_approximation.tex
\begin{tikzpicture}
\begin{axis}[xlabel=$t$,
    ymin = 0,
    ymax = 1.1,
    xmin = 0,
    xmax = 9.5,
    y axis line style={->,thick},
    x axis line style={->, thick},
    xlabel style = {at={(axis cs:9.5,0)},anchor=north west},
    xlabel = {$x$},
    axis x line*=middle,
    axis y line*=middle,
    width=0.7\textwidth,
    height= 0.35\textwidth,
    legend pos=south east,
    grid
    ]
    \addplot[color=azulcito, const plot, solid] table[row sep={\\}]
        {
            \\
            -0.0  0.01  \\
            -0.0  0.02  \\
            -0.0  0.03  \\
            0.0  0.04  \\
            0.20143057821265326  0.05  \\
            0.24012680117912907  0.06  \\
            0.27452473566734165  0.07  \\
            0.3607670380457966  0.08  \\
            0.4007187922000144  0.09  \\
            0.4573955448563032  0.1  \\
            0.4581949030678096  0.11  \\
            0.4604511716670332  0.12  \\
            0.46433273648946877  0.13  \\
            0.48152158386581895  0.14  \\
            0.4840788676646377  0.15  \\
            0.4995858074889992  0.16  \\
            0.5180411286158892  0.17  \\
            0.521225683768705  0.18  \\
            0.5561919762897224  0.19  \\
            0.5684963734476198  0.2  \\
            0.5822442036667047  0.21  \\
            0.5999551059206216  0.22  \\
            0.6190034323897923  0.23  \\
            0.6420416593948215  0.24  \\
            0.6916970579388543  0.25  \\
            0.702662426770556  0.26  \\
            0.7179631434274742  0.27  \\
            0.8254989620932781  0.28  \\
            0.9081091024884566  0.29  \\
            0.9782500228686068  0.3  \\
            0.9835912154177465  0.31  \\
            1.0846327883957447  0.32  \\
            1.0918036860817637  0.33  \\
            1.0935319921966422  0.34  \\
            1.1356696148661372  0.35  \\
            1.1609668527573638  0.36  \\
            1.1864909075569767  0.37  \\
            1.2166043915929392  0.38  \\
            1.2304029609353917  0.39  \\
            1.2633801813298513  0.4  \\
            1.2640257096813154  0.41  \\
            1.2657807349673322  0.42  \\
            1.3409001367206261  0.43  \\
            1.3424809183737691  0.44  \\
            1.3608178871457237  0.45  \\
            1.4148651254252793  0.46  \\
            1.4500700863888376  0.47  \\
            1.4992394080329097  0.48  \\
            1.5186181656595064  0.49  \\
            1.5241584379517195  0.5  \\
            1.5371047729795855  0.51  \\
            1.6075020920650462  0.52  \\
            1.6155634914215855  0.53  \\
            1.6395025293332712  0.54  \\
            1.7324855456170465  0.55  \\
            1.7492645182443336  0.56  \\
            1.7513957674876128  0.57  \\
            1.7541320492629402  0.58  \\
            1.7990510624290976  0.59  \\
            1.8638630195718822  0.6  \\
            1.8870322289869859  0.61  \\
            1.9487274673631916  0.62  \\
            1.9783196413128679  0.63  \\
            2.003063109545158  0.64  \\
            2.0292659738730165  0.65  \\
            2.045549606911596  0.66  \\
            2.1055604163941104  0.67  \\
            2.1298979729643803  0.68  \\
            2.2537395080178593  0.69  \\
            2.2839505997791325  0.7  \\
            2.286277980614691  0.71  \\
            2.3051492737859864  0.72  \\
            2.310391357138644  0.73  \\
            2.325200579452085  0.74  \\
            2.352904702326109  0.75  \\
            2.4423325135836134  0.76  \\
            2.478892783017448  0.77  \\
            2.481366557495483  0.78  \\
            2.5156867201865176  0.79  \\
            2.6313840034681566  0.8  \\
            2.687935071402705  0.81  \\
            2.8986117780276337  0.82  \\
            2.984593462610827  0.83  \\
            2.9937540922706782  0.84  \\
            3.046118696606882  0.85  \\
            3.08635559788715  0.86  \\
            3.1966158629297494  0.87  \\
            3.254855846380252  0.88  \\
            3.5606494837517304  0.89  \\
            3.6744197696440466  0.9  \\
            3.8119352636674777  0.91  \\
            3.8388169002767327  0.92  \\
            4.009063481911412  0.93  \\
            4.114465784558067  0.94  \\
            4.325556470350843  0.95  \\
            4.993387419999878  0.96  \\
            5.563366386271628  0.97  \\
            6.0514872168902585  0.98  \\
            9.841429297819689  0.99  \\
            10.059133521956218  1.0  \\
        }
        ;
    \addlegendentry {$\widehat{G}_N(x)$}
    \addplot[color=rojito, solid, thick] table[row sep={\\}]
        {
            \\
            0.0  0.0  \\
            0.1  0.03333333333333334  \\
            0.2  0.06666666666666668  \\
            0.3  0.1  \\
            0.4  0.13333333333333336  \\
            0.5  0.16666666666666669  \\
            0.6  0.2  \\
            0.7  0.23333333333333334  \\
            0.8  0.2666666666666667  \\
            0.9  0.30000000000000004  \\
            1.0  0.33333333333333337  \\
            1.1  0.36666666666666675  \\
            1.2  0.4  \\
            1.3  0.4333333333333334  \\
            1.4  0.4666666666666667  \\
            1.5  0.5  \\
            1.6  0.5333333333333334  \\
            1.7  0.5666666666666668  \\
            1.8  0.6000000000000001  \\
            1.9  0.6333333333333334  \\
            2.0  0.6666666666666667  \\
            2.1  0.6976568405139834  \\
            2.2  0.7245179063360883  \\
            2.3  0.7479521109010712  \\
            2.4  0.7685185185185185  \\
            2.5  0.7866666666666666  \\
            2.6  0.8027613412228798  \\
            2.7  0.817101051668953  \\
            2.8  0.8299319727891157  \\
            2.9  0.8414585810543005  \\
            3.0  0.8518518518518519  \\
            3.1  0.8612556364897677  \\
            3.2  0.8697916666666667  \\
            3.3  0.8775635139271503  \\
            3.4  0.8846597462514417  \\
            3.5  0.891156462585034  \\
            3.6  0.897119341563786  \\
            3.7  0.9026053080107135  \\
            3.8  0.9076638965835642  \\
            3.9  0.9123383738768354  \\
            4.0  0.9166666666666666  \\
            4.1  0.9206821336506048  \\
            4.2  0.9244142101284959  \\
            4.3  0.9278889489814314  \\
            4.4  0.931129476584022  \\
            4.5  0.934156378600823  \\
            4.6  0.9369880277252678  \\
            4.7  0.9396408631356572  \\
            4.8  0.9421296296296297  \\
            4.9  0.944467582951548  \\
            5.0  0.9466666666666667  \\
            5.1  0.9487376650006407  \\
            5.2  0.9506903353057199  \\
            5.3  0.9525335231992406  \\
            5.4  0.9542752629172382  \\
            5.5  0.9559228650137741  \\
            5.6  0.9574829931972789  \\
            5.7  0.9589617318149174  \\
            5.8  0.9603646452635751  \\
            5.9  0.9616968304127167  \\
            6.0  0.962962962962963  \\
            6.1  0.9641673385290692  \\
            6.2  0.9653139091224419  \\
            6.3  0.9664063156126648  \\
            6.4  0.9674479166666666  \\
            6.5  0.9684418145956607  \\
            6.6  0.9693908784817876  \\
            6.7  0.9702977649068092  \\
            6.8  0.9711649365628604  \\
            6.9  0.9719946789890079  \\
            7.0  0.9727891156462585  \\
            7.1  0.9735502215168947  \\
            7.2  0.9742798353909465  \\
            7.3  0.9749796709826735  \\
            7.4  0.9756513270026783  \\
            7.5  0.9762962962962963  \\
            7.6  0.976915974145891  \\
            7.7  0.9775116658233541  \\
            7.8  0.9780845934692088  \\
            7.9  0.9786359023660738  \\
            8.0  0.9791666666666666  \\
            8.1  0.9796778946298836  \\
            8.2  0.9801705334126511  \\
            8.3  0.9806454734601054  \\
            8.4  0.981103552532124  \\
            8.5  0.9815455594002307  \\
            8.6  0.9819722372453579  \\
            8.7  0.9823842867838112  \\
            8.8  0.9827823691460055  \\
            8.9  0.9831671085300677  \\
            9.0  0.9835390946502057  \\
            9.1  0.9838988849977861  \\
            9.2  0.984247006931317  \\
            9.3  0.9845839596099741  \\
            9.4  0.9849102157839142  \\
            9.5  0.9852262234533703  \\
            9.6  0.9855324074074074  \\
            9.7  0.9858291706522124  \\
            9.8  0.986116895737887  \\
            9.9  0.9863959459919056  \\
            10.0  0.9866666666666667  \\
        }
        ;
    \addlegendentry {$G_\infty(x)$}
\end{axis}
\end{tikzpicture}

%% file: figuras/miss_probability.tex
\begin{tikzpicture}
    \begin{axis}[xlabel=$t$,
        ymin = 0,
        ymax = .9,
        xmin = 0,
        xmax = 1.55,
        y axis line style={->,thick},
        x axis line style={->, thick},
        xlabel = {$\beta$},
        ylabel = {$M$},
        axis x line*=middle,
        axis y line*=middle,
        width=0.7\textwidth,
        height= 0.4\textwidth,
        legend pos=north east,
        grid
        ]

        \addplot[color=azulcito, solid, very thick] table[row sep={\\}]
        {
            \\
            0.0  0.81  \\
            0.01  0.8099190000000001  \\
            0.02  0.8096760000000001  \\
            0.03  0.8092710000000001  \\
            0.04  0.808704  \\
            0.05  0.807975  \\
            0.06  0.807084  \\
            0.07  0.806031  \\
            0.08  0.8048160000000002  \\
            0.09  0.803439  \\
            0.1  0.8019000000000002  \\
            0.11  0.800199  \\
            0.12  0.7983372977053461  \\
            0.13  0.796322791646331  \\
            0.14  0.7941640436943629  \\
            0.15  0.7918679276522128  \\
            0.16  0.7894399793079565  \\
            0.17  0.7868846577764377  \\
            0.18  0.7842055428639393  \\
            0.19  0.7814054866108218  \\
            0.2  0.7784867313503421  \\
            0.21  0.7754510028628414  \\
            0.22  0.7722995847124952  \\
            0.23  0.7690333781640768  \\
            0.24  0.7656529509081398  \\
            0.25  0.7621585769994558  \\
            0.26  0.7585502698238125  \\
            0.27  0.7548278094796816  \\
            0.28  0.7509907656455393  \\
            0.29  0.7470385167681468  \\
            0.3  0.7429702662294577  \\
            0.31  0.7387850560144084  \\
            0.32  0.7344817782976005  \\
            0.33  0.7300591852859202  \\
            0.34  0.725515897590713  \\
            0.35  0.7208504113530634  \\
            0.36  0.716061104305908  \\
            0.37  0.7111462409248186  \\
            0.38  0.7061039767936  \\
            0.39  0.7009323622899948  \\
            0.4  0.6956293456797971  \\
            0.41  0.6901927756937372  \\
            0.42  0.6846204036500054  \\
            0.43  0.67890988517577  \\
            0.44  0.6730587815731204  \\
            0.45  0.6670645608682485  \\
            0.46  0.6609245985771202  \\
            0.47  0.6546361782162052  \\
            0.48  0.64819649158287  \\
            0.49  0.6416026388266665  \\
            0.5  0.6348516283298892  \\
            0.51  0.6279403764133009  \\
            0.52  0.620865706880839  \\
            0.53  0.6136243504152812  \\
            0.54  0.6062129438352988  \\
            0.55  0.5986280292229472  \\
            0.56  0.5908660529294832  \\
            0.57  0.5829233644663496  \\
            0.58  0.5747962152872796  \\
            0.59  0.5664807574666803  \\
            0.6  0.5579730422787526  \\
            0.61  0.5492690186812061  \\
            0.62  0.540364531706868  \\
            0.63  0.531255320766019  \\
            0.64  0.5219370178618485  \\
            0.65  0.5124051457210493  \\
            0.66  0.502655115841217  \\
            0.67  0.4926822264564269  \\
            0.68  0.4824816604220663  \\
            0.69  0.4720484830197663  \\
            0.7  0.461377639683035  \\
            0.71  0.45046395364399705  \\
            0.72  0.4393021235014489  \\
            0.73  0.427886720710269  \\
            0.74  0.416212186992059  \\
            0.75  0.404272831666749  \\
            0.76  0.39206282890475663  \\
            0.77  0.3795762148991675  \\
            0.78  0.36680688495728364  \\
            0.79  0.35374859051076957  \\
            0.8  0.34039493604352933  \\
            0.81  0.3267393759363404  \\
            0.82  0.31277521122717344  \\
            0.83  0.29849558628604667  \\
            0.84  0.2838934854031637  \\
            0.85  0.2689617292890062  \\
            0.86  0.25369297148497516  \\
            0.87  0.2380796946830851  \\
            0.88  0.2221142069531521  \\
            0.89  0.20578863787582802  \\
            0.9  0.18909493457977522  \\
            0.91  0.17202485768118747  \\
            0.92  0.15456997712380638  \\
            0.93  0.13672166791749996  \\
            0.94  0.1184711057734128  \\
            0.95  0.09980926263361012  \\
            0.96  0.08072690209309286  \\
            0.97  0.061214574711965786  \\
            0.98  0.04126261321549052  \\
            0.99  0.020861127579677063  \\
            1.0  0.0  \\
            2.0  0.0  \\
            }
            ;
            \addlegendentry{Asymptotic result}
    
        \addplot+[only marks, mark=square, rojito]
        table[row sep={\\}]
        {
            \\
            0.0  0.8099278330407645  \\
            0.1  0.8092373374733165  \\
            0.2  0.7935496519949025  \\
            0.3  0.7434383893749341  \\
            0.4  0.7021645021645022  \\
            0.5  0.673237997957099  \\
            0.6  0.619425960869997  \\
            0.7  0.5756914119359534  \\
            0.8  0.5155153118323329  \\
            0.9  0.46179272797977833  \\
            1.0  0.4044191161767646  \\
            1.1  0.3411374812914262  \\
            1.2  0.28981049707312234  \\
            1.3  0.23182266009852215  \\
            1.4  0.20089721084454848  \\
            1.5  0.15419278996865204  \\
            1.6  0.12458881578947367  \\
            1.7  0.09030326443886416  \\
            1.8  0.08134599838318513  \\
            1.9  0.05896753592336346  \\
            2.0  0.05151722119231139  \\
        }
        ;
    \addlegendentry{Simulation, $N=100$}

\addplot+[only marks, mark=triangle, rojito]
        table[row sep={\\}]
        {
            \\
            0.0  0.8191862732354088  \\
            0.1  0.8012088204038257  \\
            0.2  0.7849310460813993  \\
            0.3  0.7434533826994278  \\
            0.4  0.7088882848600163  \\
            0.5  0.6607154578273712  \\
            0.6  0.5998388071730808  \\
            0.7  0.5119340469383684  \\
            0.8  0.4566977765577821  \\
            0.9  0.3642670157068063  \\
            1.0  0.2826706473062831  \\
            1.1  0.21851608667104394  \\
            1.2  0.1472437784318973  \\
            1.3  0.10632278439371012  \\
            1.4  0.07805176278785408  \\
            1.5  0.04335513451554174  \\
            1.6  0.029930328451585386  \\
            1.7  0.020557117985804174  \\
            1.8  0.011083743842364546  \\
            1.9  0.007332785087719285  \\
            2.0  0.004958246346555284  \\
        }
        ;
    \addlegendentry {Simulation, $N=1000$}

    \addplot+[only marks, mark=diamond, rojito]
    table[row sep={\\}]
        {
            \\
            0.0  0.8153210738519028  \\
            0.1  0.7995609001546828  \\
            0.2  0.784709449719143  \\
            0.3  0.7431066627416347  \\
            0.4  0.7016838027119654  \\
            0.5  0.6381424210578992  \\
            0.6  0.5670691658667206  \\
            0.7  0.4864905358837337  \\
            0.8  0.39079266182069916  \\
            0.9  0.31924041562164096  \\
            1.0  0.21376251930838608  \\
            1.1  0.14375487139516763  \\
            1.2  0.08336202232823486  \\
            1.3  0.05322549399457577  \\
            1.4  0.031798673429574764  \\
            1.5  0.00877365930599372  \\
            1.6  0.0  \\
            1.7  0.0  \\
            1.8  0.0  \\
            1.9  0.0  \\
            2.0  0.0  \\
        }
        ;
    \addlegendentry {Simulation, $N=10000$}

    \addplot[dashed, mark=none, azulcito, thick] table[row sep={\\}]
        {
            \\
            0.0  0.8100000000000003  \\
            0.01  0.8104413519014197  \\
            0.02  0.8107220126033163  \\
            0.03  0.8108429583566212  \\
            0.04  0.8108051707668361  \\
            0.05  0.8106096380472462  \\
            0.06  0.8102573562889164  \\
            0.07  0.8097493307450421  \\
            0.08  0.809086577127324  \\
            0.09  0.8082701229121535  \\
            0.1  0.80730100865446  \\
            0.11  0.8061802893071524  \\
            0.12  0.8048355755388108  \\
            0.13  0.8033441518504958  \\
            0.14  0.8017168291295953  \\
            0.15  0.7999693449144051  \\
            0.16  0.7981053560640791  \\
            0.17  0.7961346457707134  \\
            0.18  0.7940551076106934  \\
            0.19  0.7918783448324881  \\
            0.2  0.7896121870558146  \\
            0.21  0.7872393497028901  \\
            0.22  0.7847833581473161  \\
            0.23  0.7822404490360111  \\
            0.24  0.7796137169427122  \\
            0.25  0.7769129674938254  \\
            0.26  0.7741094927389321  \\
            0.27  0.7712486516059179  \\
            0.28  0.7682896387655244  \\
            0.29  0.7652596352020866  \\
            0.3  0.7621579358924474  \\
            0.31  0.7589816972157583  \\
            0.32  0.7557297508149824  \\
            0.33  0.7524025598263164  \\
            0.34  0.7490021908535497  \\
            0.35  0.7455323000095327  \\
            0.36  0.7419981317772011  \\
            0.37  0.7384014430206056  \\
            0.38  0.7347032578938009  \\
            0.39  0.7309523525035746  \\
            0.4  0.7271595210800706  \\
            0.41  0.7232479322712502  \\
            0.42  0.7193048004686986  \\
            0.43  0.715281176218058  \\
            0.44  0.7111945041658191  \\
            0.45  0.7070470157358228  \\
            0.46  0.7028282602712408  \\
            0.47  0.6985471355166795  \\
            0.48  0.6942124074278389  \\
            0.49  0.6897898031613304  \\
            0.5  0.6853598715087175  \\
            0.51  0.6807897745008263  \\
            0.52  0.6762086897124995  \\
            0.53  0.6715682563136353  \\
            0.54  0.6668323130726067  \\
            0.55  0.6620820341408066  \\
            0.56  0.6572600640701415  \\
            0.57  0.6523504947737623  \\
            0.58  0.6474242050520772  \\
            0.59  0.6424589044642494  \\
            0.6  0.6373716786871484  \\
            0.61  0.632266262395825  \\
            0.62  0.6271429400632855  \\
            0.63  0.6219396169117389  \\
            0.64  0.616655486773725  \\
            0.65  0.6113530756575608  \\
            0.66  0.6060329362365711  \\
            0.67  0.6006549118717449  \\
            0.68  0.5951780069391325  \\
            0.69  0.5896840373834501  \\
            0.7  0.5841738090302283  \\
            0.71  0.5786481702877345  \\
            0.72  0.5730397909064968  \\
            0.73  0.5673658552968123  \\
            0.74  0.5616783416411586  \\
            0.75  0.5559783296245969  \\
            0.76  0.5502669349862105  \\
            0.77  0.544536430600445  \\
            0.78  0.5386899644985303  \\
            0.79  0.532834999601709  \\
            0.8  0.5269728545082439  \\
            0.81  0.5211048757994047  \\
            0.82  0.5152324360962101  \\
            0.83  0.5093569320749313  \\
            0.84  0.5034059414737886  \\
            0.85  0.4974175513923423  \\
            0.86  0.4914301504960217  \\
            0.87  0.4854452855721169  \\
            0.88  0.47946451712982846  \\
            0.89  0.47348941713603243  \\
            0.9  0.46752156673144196  \\
            0.91  0.46155305561249726  \\
            0.92  0.45550000901927745  \\
            0.93  0.4494590686118417  \\
            0.94  0.44343189924393034  \\
            0.95  0.43742016299833053  \\
            0.96  0.4314255167994827  \\
            0.97  0.42544961004365023  \\
            0.98  0.4194940822525521  \\
            0.99  0.41356056075633546  \\
            1.0  0.4076506584117187  \\
            1.01  0.40166979753357396  \\
            1.02  0.3957161986108145  \\
            1.03  0.38979147929080427  \\
            1.04  0.38389723308040635  \\
            1.05  0.3780350272453736  \\
            1.06  0.3722064007794744  \\
            1.07  0.36641286244820775  \\
            1.08  0.3606558889116734  \\
            1.09  0.3549369229308482  \\
            1.1  0.3492573716612004  \\
            1.11  0.34361860503722  \\
            1.12  0.33794606075334105  \\
            1.13  0.3323082239927451  \\
            1.14  0.32671592282474904  \\
            1.15  0.3211704099032753  \\
            1.16  0.31567289251137626  \\
            1.17  0.3102245315038392  \\
            1.18  0.3048264403589474  \\
            1.19  0.29947968433995864  \\
            1.2  0.29418527976646514  \\
            1.21  0.28894419339539584  \\
            1.22  0.2837573419110512  \\
            1.23  0.2786255915231751  \\
            1.24  0.27354975767171835  \\
            1.25  0.2685306048366035  \\
            1.26  0.263497434060666  \\
            1.27  0.2585234309522325  \\
            1.28  0.25360923434865396  \\
            1.29  0.24875543072916603  \\
            1.3  0.24396255451238968  \\
            1.31  0.23923108844575447  \\
            1.32  0.23456146408320266  \\
            1.33  0.22995406234733734  \\
            1.34  0.22540921417201737  \\
            1.35  0.22092720122126358  \\
            1.36  0.21650825668021476  \\
            1.37  0.21215256611378952  \\
            1.38  0.20786026838862837  \\
            1.39  0.20363145665384624  \\
            1.4  0.19946617937609826  \\
            1.41  0.19536444142444778  \\
            1.42  0.19132620520054613  \\
            1.43  0.1873435526706039  \\
            1.44  0.18337785029977485  \\
            1.45  0.17947649060381357  \\
            1.46  0.17563928367771783  \\
            1.47  0.17186600206672972  \\
            1.48  0.16815638210859213  \\
            1.49  0.1645101252925617  \\
            1.5  0.16092689963122814  \\
            1.51  0.15740634104132337  \\
            1.52  0.15394805472985862  \\
            1.53  0.15055161658207633  \\
            1.54  0.14721657454786785  \\
            1.55  0.14394245002347328  \\
            1.56  0.14072873922544868  \\
            1.57  0.13757491455405693  \\
            1.58  0.13448042594341408  \\
            1.59  0.13144470219589496  \\
            1.6  0.12846715229847788  \\
            1.61  0.1255471667188808  \\
            1.62  0.12268411867951037  \\
            1.63  0.1198773654074192  \\
            1.64  0.11712624935862583  \\
            1.65  0.11443009941532029  \\
            1.66  0.11178823205463136  \\
            1.67  0.10918802978264686  \\
            1.68  0.10661754590363189  \\
            1.69  0.10410024473626271  \\
            1.7  0.10163539407936217  \\
            1.71  0.0992222542061287  \\
            1.72  0.09686007887380753  \\
            1.73  0.09454811629908026  \\
            1.74  0.09228561009888404  \\
            1.75  0.09007180019648524  \\
            1.76  0.08790592369273352  \\
            1.77  0.08578721570252269  \\
            1.78  0.0837149101565723  \\
            1.79  0.08168824056873196  \\
            1.8  0.07970644076908674  \\
            1.81  0.07776874560321577  \\
            1.82  0.07587439159802098  \\
            1.83  0.07402261759460618  \\
            1.84  0.07221266534873895  \\
            1.85  0.07044378009947978  \\
            1.86  0.06871521110660632  \\
            1.87  0.06702621215749986  \\
            1.88  0.06537604204419677  \\
            1.89  0.0637639650113366  \\
            1.9  0.062189251175766444  \\
            1.91  0.06065117691857937  \\
            1.92  0.05914902525038657  \\
            1.93  0.05768208615063406  \\
            1.94  0.05624965688178698  \\
            1.95  0.05485104227921178  \\
            1.96  0.05348555501759107  \\
            1.97  0.05215251585470744  \\
            1.98  0.05085125385343276  \\
            1.99  0.04958110658275531  \\
            2.0  0.04834142029867184  \\
        }
        ;
        \addplot[dashed, mark=none, azulcito, thick] table[row sep={\\}]
        {
            \\
            0.0  0.8100000000000003  \\
            0.01  0.8099898987378631  \\
            0.02  0.8098181320030499  \\
            0.03  0.8094850296538374  \\
            0.04  0.8089909245054001  \\
            0.05  0.8083361537594622  \\
            0.06  0.8075210604470187  \\
            0.07  0.8065459948877098  \\
            0.08  0.8054113161694297  \\
            0.09  0.8041173936517456  \\
            0.1  0.802664608496675  \\
            0.11  0.8010533552303459  \\
            0.12  0.799278146182397  \\
            0.13  0.7973535554389845  \\
            0.14  0.7952895199477222  \\
            0.15  0.7930934539047947  \\
            0.16  0.7907714828215205  \\
            0.17  0.7883285036304896  \\
            0.18  0.785768776102278  \\
            0.19  0.7830957654236719  \\
            0.2  0.7803123857421521  \\
            0.21  0.7774208695567749  \\
            0.22  0.7744233277609667  \\
            0.23  0.7713214554083674  \\
            0.24  0.7681164997859035  \\
            0.25  0.7648097380425747  \\
            0.26  0.7614017424025624  \\
            0.27  0.757893543222916  \\
            0.28  0.7542856549729687  \\
            0.29  0.7505783809649446  \\
            0.3  0.7467720720123525  \\
            0.31  0.7428670257791624  \\
            0.32  0.7388633670531917  \\
            0.33  0.734761340299276  \\
            0.34  0.7305605964598627  \\
            0.35  0.7262613151539039  \\
            0.36  0.7218635948325499  \\
            0.37  0.717367226706865  \\
            0.38  0.7127722777709231  \\
            0.39  0.7080784311694635  \\
            0.4  0.703286157915186  \\
            0.41  0.6983947150491528  \\
            0.42  0.693404519416595  \\
            0.43  0.6883152927351475  \\
            0.44  0.6831275237201613  \\
            0.45  0.6778410366693044  \\
            0.46  0.6724558821201156  \\
            0.47  0.6669725516238514  \\
            0.48  0.6613910054700094  \\
            0.49  0.6557120076263484  \\
            0.5  0.6499355980168885  \\
            0.51  0.64406214754008  \\
            0.52  0.6380923566306741  \\
            0.53  0.6320272540057742  \\
            0.54  0.6258675397369763  \\
            0.55  0.6196138743047893  \\
            0.56  0.6132671972127723  \\
            0.57  0.6068289957128683  \\
            0.58  0.6003003607839117  \\
            0.59  0.5936826429416046  \\
            0.6  0.5869774520454727  \\
            0.61  0.580186656770792  \\
            0.62  0.5733111068192986  \\
            0.63  0.5663537250206623  \\
            0.64  0.5593162691184684  \\
            0.65  0.5522010851214114  \\
            0.66  0.5450100291645923  \\
            0.67  0.5377464928846333  \\
            0.68  0.5304129343979724  \\
            0.69  0.5230116582400296  \\
            0.7  0.5155459199550816  \\
            0.71  0.5080191176874492  \\
            0.72  0.5004345613370343  \\
            0.73  0.4927957626661093  \\
            0.74  0.48510567757642137  \\
            0.75  0.4773685767412772  \\
            0.76  0.46958861648242955  \\
            0.77  0.4617700543619502  \\
            0.78  0.4539161068130499  \\
            0.79  0.44603150636533456  \\
            0.8  0.4381217414645051  \\
            0.81  0.43018889279605993  \\
            0.82  0.4222404504996737  \\
            0.83  0.4142785245630676  \\
            0.84  0.40631044556250456  \\
            0.85  0.3983386126159258  \\
            0.86  0.39036985307352406  \\
            0.87  0.38240888789793953  \\
            0.88  0.3744593222926512  \\
            0.89  0.366527836335801  \\
            0.9  0.358619356753857  \\
            0.91  0.35073810430847463  \\
            0.92  0.34288918614950187  \\
            0.93  0.3350782659915247  \\
            0.94  0.3273102535043608  \\
            0.95  0.31959001843264323  \\
            0.96  0.31192237631936803  \\
            0.97  0.3043120743884994  \\
            0.98  0.29676377766969186  \\
            0.99  0.2892820554471101  \\
            1.0  0.2818713681124094  \\
            1.01  0.2745360544991827  \\
            1.02  0.26728031977261796  \\
            1.03  0.26010822394378413  \\
            1.04  0.2530236710729342  \\
            1.05  0.24603039922054343  \\
            1.06  0.23913197119856794  \\
            1.07  0.23233176616768977  \\
            1.08  0.22563297211921593  \\
            1.09  0.21903857927289125  \\
            1.1  0.2125513744142994  \\
            1.11  0.20617336313450094  \\
            1.12  0.19990682735500523  \\
            1.13  0.19375450760406018  \\
            1.14  0.18771834566730633  \\
            1.15  0.18180006594646542  \\
            1.16  0.1759999170403738  \\
            1.17  0.17032013381734284  \\
            1.18  0.16476206525590148  \\
            1.19  0.15932638360273882  \\
            1.2  0.15401240516836767  \\
            1.21  0.1488221034331936  \\
            1.22  0.1437557407005063  \\
            1.23  0.13881170402067533  \\
            1.24  0.1339913573156015  \\
            1.25  0.1292945884936781  \\
            1.26  0.12471897225434446  \\
            1.27  0.12026582387715176  \\
            1.28  0.11593369994262069  \\
            1.29  0.11172094539361303  \\
            1.3  0.10762783768220019  \\
            1.31  0.10365149868975564  \\
            1.32  0.09979186959915039  \\
            1.33  0.09604697635664419  \\
            1.34  0.09241499443711407  \\
            1.35  0.0888950806657174  \\
            1.36  0.08548410584599649  \\
            1.37  0.08218195518503663  \\
            1.38  0.0789847018579961  \\
            1.39  0.07589236180784276  \\
            1.4  0.0729011404374057  \\
            1.41  0.07001050868294788  \\
            1.42  0.06721688423506975  \\
            1.43  0.06451948804240673  \\
            1.44  0.06191473125279286  \\
            1.45  0.059401837423489826  \\
            1.46  0.056977028913853395  \\
            1.47  0.05463964875197462  \\
            1.48  0.05238586978251434  \\
            1.49  0.05021485069627208  \\
            1.5  0.04812326785761411  \\
            1.51  0.046109612687562744  \\
            1.52  0.04417131496925086  \\
            1.53  0.04230612450367791  \\
            1.54  0.04051231732556376  \\
            1.55  0.038786850528927046  \\
            1.56  0.0371286237585641  \\
            1.57  0.035534637559519744  \\
            1.58  0.034003390910066364  \\
            1.59  0.032532804573536564  \\
            1.6  0.031120399013911376  \\
            1.61  0.02976494048725283  \\
            1.62  0.02846381071127068  \\
            1.63  0.02721552281352392  \\
            1.64  0.02601841388065573  \\
            1.65  0.024869947169117466  \\
            1.66  0.023769143374469455  \\
            1.67  0.022714006825861235  \\
            1.68  0.021702624041184657  \\
            1.69  0.020733854064872735  \\
            1.7  0.019805752632536096  \\
            1.71  0.01891680437239193  \\
            1.72  0.018065827039492433  \\
            1.73  0.017251007497413602  \\
            1.74  0.01647105338391678  \\
            1.75  0.015724831200056934  \\
            1.76  0.01501071439181541  \\
            1.77  0.014327500757639654  \\
            1.78  0.013674170317550336  \\
            1.79  0.013049318759165394  \\
            1.8  0.01245173642934848  \\
            1.81  0.011880560903480417  \\
            1.82  0.011334629899337008  \\
            1.83  0.01081265770397419  \\
            1.84  0.010313966731481824  \\
            1.85  0.009837579640803173  \\
            1.86  0.00938228278716583  \\
            1.87  0.008947404206841827  \\
            1.88  0.00853213547922233  \\
            1.89  0.008135543197297089  \\
            1.9  0.007756730209258287  \\
            1.91  0.007395126656324799  \\
            1.92  0.007049986239838592  \\
            1.93  0.006720421686574926  \\
            1.94  0.006405855813198365  \\
            1.95  0.006105705941597341  \\
            1.96  0.005819337332143645  \\
            1.97  0.00554595064916666  \\
            1.98  0.00528515188145737  \\
            1.99  0.00503639723423244  \\
            2.0  0.004799149337884596  \\
        }
        ;
        \addplot[dashed, mark=none, azulcito, thick] table[row sep={\\}]
        {
            \\
            0.0  0.8099999999999997  \\
            0.01  0.809927961783599  \\
            0.02  0.8096939997826045  \\
            0.03  0.8092981907584706  \\
            0.04  0.8087406127641898  \\
            0.05  0.8080213457997589  \\
            0.06  0.8071404724813354  \\
            0.07  0.8060980787287241  \\
            0.08  0.8048942544758916  \\
            0.09  0.8035290944093362  \\
            0.1  0.8020026987392528  \\
            0.11  0.8003151740085774  \\
            0.12  0.798467217248701  \\
            0.13  0.7964674723098342  \\
            0.14  0.794324728502686  \\
            0.15  0.7920460069876204  \\
            0.16  0.7896370028228633  \\
            0.17  0.7871023442022718  \\
            0.18  0.7844457964720347  \\
            0.19  0.7816704112001361  \\
            0.2  0.778778647081368  \\
            0.21  0.7757724645243904  \\
            0.22  0.7726534034683359  \\
            0.23  0.7694226451624666  \\
            0.24  0.7660810592836683  \\
            0.25  0.7626292526365103  \\
            0.26  0.759067596191789  \\
            0.27  0.7553962645142368  \\
            0.28  0.7516152552267658  \\
            0.29  0.7477244132948956  \\
            0.3  0.7437234491272424  \\
            0.31  0.7396119587309187  \\
            0.32  0.7353894365019265  \\
            0.33  0.7310552866453482  \\
            0.34  0.7266088472424891  \\
            0.35  0.722049387362171  \\
            0.36  0.7173761270761666  \\
            0.37  0.7125882470904152  \\
            0.38  0.707684897095931  \\
            0.39  0.7026652084335621  \\
            0.4  0.697528297664275  \\
            0.41  0.692273281992344  \\
            0.42  0.6868992891239123  \\
            0.43  0.6814054617515906  \\
            0.44  0.6757909776886933  \\
            0.45  0.6700550448237185  \\
            0.46  0.6641969233907768  \\
            0.47  0.6582159269490264  \\
            0.48  0.6521114476603732  \\
            0.49  0.645882951865168  \\
            0.5  0.6395299986487646  \\
            0.51  0.6330522439390274  \\
            0.52  0.6264494741185214  \\
            0.53  0.6197215866351661  \\
            0.54  0.6128686316647721  \\
            0.55  0.6058908038066437  \\
            0.56  0.5987884674491947  \\
            0.57  0.5915621670512812  \\
            0.58  0.5842126416463691  \\
            0.59  0.5767408301668795  \\
            0.6  0.569147893920149  \\
            0.61  0.5614352371534395  \\
            0.62  0.5536045054331267  \\
            0.63  0.5456575965076903  \\
            0.64  0.5375966979004118  \\
            0.65  0.529424275439778  \\
            0.66  0.521143096023917  \\
            0.67  0.5127562260205861  \\
            0.68  0.5042670683266405  \\
            0.69  0.49567934445774564  \\
            0.7  0.48699711204423773  \\
            0.71  0.4782247786224779  \\
            0.72  0.46936710014425503  \\
            0.73  0.46042917383085924  \\
            0.74  0.45141646030874644  \\
            0.75  0.44233477648368663  \\
            0.76  0.4331902822480567  \\
            0.77  0.4239894797734861  \\
            0.78  0.41473922007943453  \\
            0.79  0.4054466600246758  \\
            0.8  0.39611929556678627  \\
            0.81  0.3867648890580419  \\
            0.82  0.37739151583013547  \\
            0.83  0.368007487130755  \\
            0.84  0.358621341704747  \\
            0.85  0.3492418406002404  \\
            0.86  0.3398779115192949  \\
            0.87  0.3305386192160753  \\
            0.88  0.3212331513197413  \\
            0.89  0.31197075766067106  \\
            0.9  0.3027607302913471  \\
            0.91  0.29361234306919676  \\
            0.92  0.2845348361233869  \\
            0.93  0.27553735038215954  \\
            0.94  0.266628914793054  \\
            0.95  0.2578183655348353  \\
            0.96  0.24911434204976488  \\
            0.97  0.24052522438011675  \\
            0.98  0.23205908600097103  \\
            0.99  0.22372369965757744  \\
            1.0  0.2155264448390733  \\
            1.01  0.20747430938169162  \\
            1.02  0.1995738583400214  \\
            1.03  0.19183120600742828  \\
            1.04  0.1842519654002521  \\
            1.05  0.17684127356867244  \\
            1.06  0.16960373741840226  \\
            1.07  0.16254343996935552  \\
            1.08  0.15566393106292037  \\
            1.09  0.14896820309803482  \\
            1.1  0.1424587334283719  \\
            1.11  0.1361374336386681  \\
            1.12  0.1300056969501631  \\
            1.13  0.1240643935932868  \\
            1.14  0.11831387536464011  \\
            1.15  0.11275401140238683  \\
            1.16  0.10738419660757824  \\
            1.17  0.10220335694583012  \\
            1.18  0.09721000782879761  \\
            1.19  0.09240225081782572  \\
            1.2  0.0877778084144319  \\
            1.21  0.08333405379001982  \\
            1.22  0.07906803876550092  \\
            1.23  0.07497652205185056  \\
            1.24  0.07105598599566149  \\
            1.25  0.06730270046297346  \\
            1.26  0.06371270017733043  \\
            1.27  0.060281864037682306  \\
            1.28  0.057005906860338254  \\
            1.29  0.05388041867580076  \\
            1.3  0.05090088979048741  \\
            1.31  0.048062734673255794  \\
            1.32  0.04536131456542062  \\
            1.33  0.04279195873594421  \\
            1.34  0.04034997886084955  \\
            1.35  0.03803070025525113  \\
            1.36  0.03582946674360543  \\
            1.37  0.0337416636258862  \\
            1.38  0.03176272320129468  \\
            1.39  0.029888151376582  \\
            1.4  0.02811352564591663  \\
            1.41  0.026434510401891886  \\
            1.42  0.02484686544071135  \\
            1.43  0.023346453249038255  \\
            1.44  0.021929245139200547  \\
            1.45  0.020591326302788492  \\
            1.46  0.019328899854809158  \\
            1.47  0.018138286753914903  \\
            1.48  0.017015936422244228  \\
            1.49  0.015958420434914463  \\
            1.5  0.014962434206329272  \\
            1.51  0.014024802227114571  \\
            1.52  0.013142471884308235  \\
            1.53  0.012312513731303208  \\
            1.54  0.011532124164334377  \\
            1.55  0.010798618055772395  \\
            1.56  0.010109428943476781  \\
            1.57  0.009462106054934345  \\
            1.58  0.008854312978038598  \\
            1.59  0.008283823084200418  \\
            1.6  0.007748516662789245  \\
            1.61  0.007246377098460101  \\
            1.62  0.006775489376745825  \\
            1.63  0.006334034929869772  \\
            1.64  0.00592028856998349  \\
            1.65  0.005532614647900743  \\
            1.66  0.005169463087440946  \\
            1.67  0.004829368871594011  \\
            1.68  0.004510944001740669  \\
            1.69  0.0042128774508761796  \\
            1.7  0.003933930493985131  \\
            1.71  0.0036729333717333786  \\
            1.72  0.003428782207409729  \\
            1.73  0.003200436673817493  \\
            1.74  0.002986914454627598  \\
            1.75  0.0027872918948846992  \\
            1.76  0.0026006983791849216  \\
            1.77  0.002426314272425111  \\
            1.78  0.002263368233541539  \\
            1.79  0.002111135613956879  \\
            1.8  0.0019689347438206833  \\
            1.81  0.0018361250744845716  \\
            1.82  0.0017121049238190066  \\
            1.83  0.0015963093218779381  \\
            1.84  0.0014882079549926018  \\
            1.85  0.0013873032060825203  \\
            1.86  0.0012931282887310735  \\
            1.87  0.001205245472379448  \\
            1.88  0.0011232443958412409  \\
            1.89  0.0010467403310345043  \\
            1.9  0.0009753729387139241  \\
            1.91  0.0009088048115645589  \\
            1.92  0.0008467198679527109  \\
            1.93  0.0007888220385620887  \\
            1.94  0.0007348342993278083  \\
            1.95  0.0006844976395363845  \\
            1.96  0.0006375693704486988  \\
            1.97  0.0005938226862236691  \\
            1.98  0.0005530455845546514  \\
            1.99  0.0005150394497441857  \\
            2.0  0.0004796191913761153  \\
        }
        ;
    \end{axis}
\end{tikzpicture}

%% file: figuras/miss_probability_erlang.tex
\begin{tikzpicture}
    \begin{axis}[xlabel=$t$,
        ymin = 0,
        ymax = 1.1,
        xmin = 0,
        xmax = 1.55,
        y axis line style={->,thick},
        x axis line style={->, thick},
        xlabel = {$\beta$},
        ylabel = {$M$},
        axis x line*=middle,
        axis y line*=middle,
        width=0.7\textwidth,
        height= 0.4\textwidth,
        legend pos=north east,
        grid
        ]
        \addplot[color=azulcito, solid, very thick] table[row sep={\\}]
        {
            \\
            0.01  0.7546120356766864  \\
            0.02  0.7544449639259481  \\
            0.03  0.7541620172005645  \\
            0.04  0.7537603170084277  \\
            0.05  0.7532375385500749  \\
            0.06  0.7525918009010072  \\
            0.07  0.7518215793800226  \\
            0.08  0.750925636329515  \\
            0.09  0.7499029665896877  \\
            0.1  0.7487527545076006  \\
            0.11  0.7474743399493853  \\
            0.12  0.7460671913307623  \\
            0.13  0.7445308841309286  \\
            0.14  0.7428650836982826  \\
            0.15  0.7410695314262966  \\
            0.16  0.7391440335733321  \\
            0.17  0.7370884521551226  \\
            0.18  0.7349026974448676  \\
            0.19  0.7325867216995979  \\
            0.2  0.73014051378672  \\
            0.21  0.727564094423971  \\
            0.22  0.7248575117644761  \\
            0.23  0.7220208370715951  \\
            0.24  0.7190541602220539  \\
            0.25  0.7159575847695894  \\
            0.26  0.7127312222848312  \\
            0.27  0.7093751856794924  \\
            0.28  0.7058895812196058  \\
            0.29  0.702274498952457  \\
            0.3  0.6985300013269845  \\
            0.31  0.6946561098941826  \\
            0.32  0.6906527901520202  \\
            0.33  0.6865199348500141  \\
            0.34  0.6822573463929076  \\
            0.35  0.6778647193258355  \\
            0.36  0.6733416241772204  \\
            0.37  0.6686874940480714  \\
            0.38  0.6639016151692615  \\
            0.39  0.658983122128292  \\
            0.4  0.6539309976872497  \\
            0.41  0.6487440762910045  \\
            0.42  0.6434210498475854  \\
            0.43  0.6379604780040697  \\
            0.44  0.6323607761733796  \\
            0.45  0.6266202580040164  \\
            0.46  0.6207370923616259  \\
            0.47  0.6147093748996769  \\
            0.48  0.6085350843383782  \\
            0.49  0.602212068144531  \\
            0.5  0.5957378725474007  \\
            0.51  0.5891104076780485  \\
            0.52  0.5823269178517987  \\
            0.53  0.5753850478749293  \\
            0.54  0.5682826330870805  \\
            0.55  0.5610161190481429  \\
            0.56  0.5535829649004337  \\
            0.57  0.5459802251134221  \\
            0.58  0.5382048582173048  \\
            0.59  0.530253161540818  \\
            0.6  0.5221236298718127  \\
            0.61  0.5138081027530722  \\
            0.62  0.5053129278625056  \\
            0.63  0.4966254359808239  \\
            0.64  0.48774504222685744  \\
            0.65  0.47866802477927456  \\
            0.66  0.46939056166332566  \\
            0.67  0.4599087300204547  \\
            0.68  0.4502185048351951  \\
            0.69  0.44031575730975053  \\
            0.7  0.43019625349475665  \\
            0.71  0.41984242158641305  \\
            0.72  0.4092801672498363  \\
            0.73  0.39849326035955557  \\
            0.74  0.38746224068258356  \\
            0.75  0.37619166639423873  \\
            0.76  0.3646766398247572  \\
            0.77  0.3529121461496177  \\
            0.78  0.3408930513105602  \\
            0.79  0.32861409991897345  \\
            0.8  0.3160699130242943  \\
            0.81  0.30325498578909493  \\
            0.82  0.2901636851061187  \\
            0.83  0.2767902471337054  \\
            0.84  0.2631231799526635  \\
            0.85  0.2491732348989005  \\
            0.86  0.23491687206940803  \\
            0.87  0.2203545213850704  \\
            0.88  0.20547917146978814  \\
            0.89  0.19028483035074695  \\
            0.9  0.17476350720591255  \\
            0.91  0.15890890662038734  \\
            0.92  0.14271396408433812  \\
            0.93  0.12617086674219122  \\
            0.94  0.10927351017538205  \\
            0.95  0.09201360261775997  \\
            0.96  0.07438382166971147  \\
            0.97  0.05637548363817878  \\
            0.98  0.03798169820808137  \\
            0.99  0.019192335171863206  \\
            1.0  0.0  \\
    }
        ;
    \addlegendentry {Asymptotic limit}

    \addplot[only marks, mark=square, rojito] table[row sep={\\}]
        {
            \\
            0.0  0.7595111373940469  \\
            0.1  0.7579771222155328  \\
            0.2  0.7358378056052475  \\
            0.3  0.7016629933880986  \\
            0.4  0.6551792828685259  \\
            0.5  0.6069523997629864  \\
            0.6  0.5449392712550607  \\
            0.7  0.47955537911869794  \\
            0.8  0.4144869215291751  \\
            0.9  0.32970576380491734  \\
            1.0  0.2556346381969158  \\
            1.1  0.19564775269473256  \\
            1.2  0.13407258064516125  \\
            1.3  0.09363672976238424  \\
            1.4  0.053230394148719995  \\
            1.5  0.032258064516129004  \\
            1.6  0.01786439301664633  \\
            1.7  0.011580658268996347  \\
            1.8  0.006485610052695545  \\
            1.9  0.0008001600320064473  \\
            2.0  0.0  \\
        }
        ;
        \addlegendentry {Simulation, $N=1000$}

        \addplot[color=azulcito, dashed] table[row sep={\\}]
        {
            \\
            0.1  0.7494054287048043  \\
            0.2  0.7318196570286837  \\
            0.3  0.7020888231876847  \\
            0.4  0.6611162374898669  \\
            0.5  0.6098823436601344  \\
            0.6  0.5492549204556438  \\
            0.7  0.4806984960718207  \\
            0.8  0.40680892686458836  \\
            0.9  0.3314361434936408  \\
            1.0  0.2591837253222325  \\
            1.1  0.1943663312122577  \\
            1.2  0.14002447534926882  \\
            1.3  0.09725651261369964  \\
            1.4  0.06546245502825915  \\
            1.5  0.04293109889844912  \\
            1.6  0.02757781996556431  \\
            1.7  0.01743280939084646  \\
            1.8  0.01088398605247158  \\
            1.9  0.006736354521423234  \\
            2.0  0.004135846548032875  \\
        }
        ;
        \addlegendentry{Finite $N$ approximation}

        \addplot[color=teal, mark=o, only marks] table[row sep={\\}]
        {
            \\
            0.0  0.9985480190831778  \\
            0.1  0.9982297198792044  \\
            0.2  0.9981269510926118  \\
            0.3  0.9922560660815695  \\
            0.4  0.9801345059493016  \\
            0.5  0.9445253955037469  \\
            0.6  0.8911046147421173  \\
            0.7  0.8067754077791719  \\
            0.8  0.7059795517918  \\
            0.9  0.5985514743921365  \\
            1.0  0.4653589315525877  \\
            1.1  0.35972193401120567  \\
            1.2  0.25339836048562836  \\
            1.3  0.16088360946129  \\
            1.4  0.1049204052098408  \\
            1.5  0.06399099006859832  \\
            1.6  0.04032837352488461  \\
            1.7  0.02118557751069461  \\
            1.8  0.013681970203709315  \\
            1.9  0.010773452586644972  \\
            2.0  0.008460066471950856  \\
        }
        ;
        \addlegendentry{LRU (simulated)}
    \end{axis}
\end{tikzpicture}